\documentclass[A4j,11pt]{article}
\usepackage{tikz}

\usepackage{amssymb,amsmath,amsthm,amscd,amsfonts}
\usepackage{amsmath,amssymb,amsthm,amscd,ascmac,mathrsfs,bm,type1cm,multirow,array,enumerate,mathrsfs,tabularx}
\usepackage{graphicx}
\usepackage{graphics}
\usepackage{pdfpages}

\newcommand{\ve}{\bm{e}} 

\newcommand{\vv}{{\bm{v}}} 
\newcommand{\vw}{{\bm{w}}} 

\newcommand{\vN}{\bm{N}} 

\newcommand{\vV}{\bm{V}}
\newcommand{\vX}{\bm{X}}
\newcommand{\vY}{\bm{Y}}

\begin{document}

\title{{\bf Isoparametric submanifolds\\
in a Riemannian Hilbert manifold
}}
\author{{\bf Naoyuki Koike}}
\date{}
\maketitle

\begin{abstract}
In this paper, we introduce the notion of a regularizable submanifold in a Riemannian Hilbert manifold.  
This submanifold is defined as a curvature-invariant submanifold such that its shape operators and its normal Jacobi operators 
are regularizable, where ``the operators are regularizable'' means that the operators are compact and that their regularized 
traces and the usual traces of their squares exist.  
Furthermore, we introduce the notion of an isoparametric submanifold in a Riemannian Hilbert manifold.  This submanifold is defined 
as a regularizable submanifold with flat section and trivial normal holonomy group satisfying the constancy of the regularized mean curvatures in the radial direction 
of the parallel submanifolds.  
For a curvature-adapted regularizable submanifold $M$ with trivial normal holonomy group in a locally symmetric Riemannian Hilbert manifold, we prove that if, 
for any parallel normal vecrtor field $\widetilde{\xi}$ of $M$, the shape operaors $A_{\widetilde{\xi}_x}$ and the normal Jacobi operator $\widetilde R(\widetilde{\xi}_x)$ 
are independent of the base point $x(\in M)$ (up to orthogonal equivalent), then it is isoparametric under some additional conditions.  
Also, we define the notion of an equifocal submanifold in a Riemannian Hilbert manifold.  We prove that the principal orbits of a certain kind of Hilbert Lie group action 
on the Riemannian Hilbert manifold $\mathcal A_P^{H^s}$ consisting of all $H^s$-connections of a $G$-bundle $P$ over a compact Riemannian manifold $B$ are equifocal, where 
$G$ is a semi-simple Lie group and $s>\frac{1}{2}\,{\rm dim}\,B-1$.  
\end{abstract}

\section{Introduction}
For each unit normal vector $\xi$ of a general Hilbert submanifold in a (separable) Hilbert space, 
the set of all focal points of the submanifold along the normal geodesic $\gamma_{\xi}$ is possible to have accumulating points 
and the multiplicities of each focal points are possible to be infinite because the tangent space of the submanifold is of infinite dimension.  
Thus the focal structure of a general Hilbert submanifold in a Hilbert space is complicate in comparison with the finite dimensional case.  
For example, the shape operator $A_{\xi}$ for any normal vector $\xi$ of the submanifold is a self-adjoint operator but the set of the eigenvalues of $A_{\xi}$ 
is possible to have accumulating points and hence it is possible to be a complicate operator.  
So, R. S. Palais and C. L. Terng (\cite{PaTe}) introduced the notion of a {\it proper Fredholm submanifold} in a Hilbert space 
as a Hilbert submanifold with good focal structure.  The proper Fredholm submanifold is defined as an immersed (Hilbert) submanifold of finite codimension 
in a Hilbert space satisfying the following conditions:

\vspace{0.15truecm}

\noindent
(P)\ The restriction of $\exp^{\perp}$ to the ball normal bundle of any radius is proper;

\noindent
(F)\ The normal exponential map $\exp^{\perp}$ of the submanifold is a Fredholm map.  

\vspace{0.15truecm}

\noindent
Let $M$ be a proper Fredholm submanifold in a Hilbert space $(V,\langle\,\,,\,\,\rangle)$ immersed by $f$.  
Then, for each unit normal vector $\xi$ of $M$, the set of of all focal points of the submanifold along the normal geodesic $\gamma_{\xi}$ of $\xi$-direction 
has no accumulating point and the multiplicities of each focal points are finite.  In fact, from the condition (P), it follows that the set of of all focal points of $M$ 
along the normal geodesic $\gamma_{\xi}$ of $\xi$-direction has no accumulating point and, furthermore, from the condition (F), 
it follows that the multiplicities of each focal points are finite.  Thus a proper Fredholm submanifold has a good focal structure similar to finite dimensional submanifolds.  
Let $d_p^2:M\to\mathbb R$ be the squared distance function from a fixed point $p$ of $V$ other than the focal set of $M$.  This function $d_p^2$ is defined by 
$d_p^2(x):=d(p,f(x))^2$ ($x\in M$), where $d$ is the distance function of $(V,\langle\,\,,\,\,\rangle)$.  
Since $p$ does not belong to the focal set of $M$, $d_p^2$ is a Morse function.  Furthermore, according to the index theorem for the squared distance function, 
it follows from the above facts for the focal points that the index of $d_p^2$ at each critical point is finite.  Hence, by applying the infinite dimensional Morse theory to 
$d_p^2$, we can investigate the topology of $M$.  On the other hand, since the ambient space is a Hilbert space (hence flat), the focal radii (which mean real nembers $s$ such that 
$\gamma_{\xi}(s)$ is a focal point) of $M$ along the normal geodesic $\gamma_{\xi}$ are equal to the inverse numbers of the nonzero eignenvalues of the shape operator $A_{\xi}$.  
This together with the above fact for the focal structure implies that the shape operator $A_{\xi}$ of $M$ for any unit normal vector $\xi$ is a compact operator.  

Let $M$ be a proper Fredholm submanifold in a Hilbert space $(V,\langle\,\,,\,\,\rangle)$ and $A$ the shape tensor of $M$.  
C. L. Terng (\cite{T1}) called $M$ an {\it isoparametric submanifold} if it satisfies the following conditions:

\vspace{0.15truecm}

\noindent
(TI-i)\ \ The normal holonomy group of $M$ is trivial;

\noindent
(TI-ii)\ \ For any parallel normal vector field $\widetilde{\xi}$ of $M$, $A_{\widetilde{\xi}_x}$ and $A_{\widetilde{\xi}_y}$ are orthogonally 

equivalent to each other for all $x,y\in M$.  

\vspace{0.15truecm}

\noindent
The condition (TI-ii) is equivalent to the following condition:

\vspace{0.15truecm}

\noindent
(TI-ii$'$)\ \ For any parallel normal vector field $\widetilde{\xi}$ of $M$, the set of all eigenvalues of 

$A_{\widetilde{\xi}_x}$ is independent of the choice of $x\in M$ with considering the multiplicities.  

\vspace{0.15truecm}

We recall two definitions of the regularized traces of a compact self-adjoint operator of a Hilbert space in the senses of \cite{KT} and \cite{HLO}.  
Let $\mathcal A$ be a compact self-adjoint operator of a (separable) Hilbert space $(V,\langle\,\,,\,\,\rangle)$.  
C. King and C. L. Terng (\cite{KT}) defined the regularized trace (which is denoted by ${\rm Tr}_{\zeta}\,\mathcal A$ in this paper) of 
a compact self-adjoint operator $A$ by 
$$\begin{array}{c}
\displaystyle{{\rm Tr}_{\zeta}\,\mathcal A:=\lim_{s\downarrow 1}(\sum\limits_{i=1}^{\infty}((\lambda^+_i)^s-(\lambda^-_i)^s)}\\
\displaystyle{(-\lambda^-_1\leq-\lambda^-_2\leq\cdots\leq 0\leq\cdots\leq\lambda^+_2\leq\lambda^+_1\,:\,{\rm the}\,\,{\rm spectrum}\,\,{\rm of}\,\,\mathcal A).}
\end{array}$$
In \cite{HLO}, this trace was called the {\it $\zeta$-trace}.  In this paper, we use this terminology.  
Later, E. Heintze, X. Liu and C. Olmos (\cite{HLO}) defined the regularized trace (which is denoted by ${\rm Tr}_r\,\mathcal A$ in this paper) of $\mathcal A$ by 
$$\begin{array}{c}
\displaystyle{{\rm Tr}_r\,\mathcal A:=\sum\limits_{i=1}^{\infty}(\lambda^+_i-\lambda^-_i)}\\
\displaystyle{(-\lambda^-_1\leq-\lambda^-_2\leq\cdots\leq 0\leq\cdots\leq\lambda^+_2\leq\lambda^+_1\,:\,{\rm the}\,\,{\rm spectrum}\,\,{\rm of}\,\,\mathcal A).}
\end{array}$$
The regularized trace in the sense of \cite{HLO} is easier to handle than one in the sense of \cite{KT}.  In almost all relevant cases, 
these regularized traces coincide.  In this paper, we call the regularized trace in the sense of \cite{HLO} the {\it regularized trace}.  
See \cite{K3} about an example of a compact self-adjoint operator whose $\zeta$-trace and regularized trace are different.  

Let $M$ be a proper Fredholm submanifold in a Hilbert space immersed by $f$ and $A$ the shape tensor.  
Note that, for each normal vector $\xi$ of $M$, the shape operator $A_{\xi}$ is a self-adjoint compact operator as above.  
E. Heintze, X. Liu and C. Olmos (\cite{HLO}) called $M$ a {\it regularizable submanifold} if it satisfies the following condition:

\vspace{0.15truecm}

\noindent
(R)\ \ For any $x\in M$ and any normal vector $\xi$ of $M$ at $x$, there exist the regularized 

trace ${\rm Tr}_r\,(A_x)_{\xi}$ of the shape operator $(A_x)_{\xi}$ of $M$ for $\xi$ and the usual trace 

${\rm Tr}\,(A_x)_{\xi}^2$ of $(A_x)_{\xi}^2$.

\vspace{0.15truecm}

\noindent
In particular, they called $M$ a {\it minimal regularizable submanifold} if it is regularizable and if it satisfies the following condition:

\vspace{0.15truecm}

\noindent
(MR)\ \ For any $x\in M$ and any normal vector $\xi$ of $M$ at $x$, ${\rm Tr}_r\,(A_x)_{\xi}$ vanishes.  

\vspace{0.15truecm}

Let $M$ be a regularizable submanifold in $(V,\langle\,\,\,\,\rangle)$.  They (\cite{HLO}) called $M$ an {\it isoparametric submanifold} if it satisfying 
the following conditions:

\vspace{0.15truecm}

\noindent
(HLOI-i)\ \ The normal holonomy group of $M$ is trivial;

\noindent
(HLOI-ii)\ \ The parallel submanifolds sufficiently close to $M$ are of constant 

regularized mean curvature in the radial direction.

\vspace{0.15truecm}

\noindent
Here the {\it parallel submanifold} sufficiently close to $M$ means the image of the immersion $\eta_{\widetilde{\xi}}$ defined by 
$$\eta_{\widetilde{\xi}}:M\hookrightarrow V\quad\left(\mathop{\Longleftrightarrow}_{\rm def}\,\,\eta_{\widetilde{\xi}}(x):=\exp^{\perp}(\widetilde{\xi}_x)\,\,\,(x\in M)\right)$$
for some parallel normal vector field $\widetilde{\xi}$ (whose norm is sufficiently small) of $M$.  The immersion $\eta_{\widetilde{\xi}}$ is called 
the {\it end-point map for} $\widetilde{\xi}$.  Set $\displaystyle{M_{\widetilde{\xi}}:=\eta_{\widetilde{\xi}}(M)}$.  
Also, ``the parallel submanifold {\it $M_{\widetilde{\xi}}$ is of constant regularized mean curvature in the radial direction}'' means that 
the regularized trace of the shape operator $A_{\gamma_{\widetilde{\xi}_x}'(1)}$ of $M_{\widetilde{\xi}}$ for $\gamma_{\widetilde{\xi}_x}'(1)$ is independent of 
the choice of $x\in M$, where $\gamma_{\widetilde{\xi}_x}$ is the normal geodesic of $M$ with $\gamma_{\widetilde{\xi}_x}'(0)=\widetilde{\xi}_x$.  
It is shown that the conditions (TI-ii) and (HLOI-ii) are equivalent to each other (see Corollary 4.3 of \cite{HLO}).  

Next we recall the notion of an isoparametric submanifold with flat section in a finite dimensional Riemannian manifold.  
Let $M$ be an immersed submanifold in a finite dimensional Riemannian manifold $(\widetilde M,\widetilde g)$.  
E. Heintze, X. Liu and C. Olmos (\cite{HLO}) called $M$ an {\it isoparametric submanifold} if it satisfies the following conditions:

\vspace{0.15truecm}

\noindent
(HLOI'-i)\ \ The normal holonomy group of $M$ is trivial;

\noindent
(HLOI'-ii)\ \ The parallel submanifolds sufficiently close to $M$ are of constant mean 

curvature in the radial direction;

\noindent
(HLOI'-iii)\ \ $M$ has sections.

\vspace{0.15truecm}

\noindent
Here ``$M$ has sections'' means that, for any point $x\in M$, the normal umbrella $\Sigma_x:=\exp^{\perp}(T^{\perp}_xM)$ is totally geodesic.  
Furthermore, assume that $M$ satisfies the following condition:

\vspace{0.15truecm}

\noindent
(HLOI'-iv)\ \  The induced metric on each section of $M$ is flat.

\vspace{0.15truecm}

\noindent
Then they called $M$ an {\it isoparametric submanifold with flat section}.  

In this paper, we shall introduce the notions of a CSJ-submanifold, a regularizable submanifold, an isoparametric submanifold 
(which is a generalized notion of an isoparametric submanifold in the sense of Heintze-Liu-Olmos), a weakly isoparametric submanifold 
(which is a generalized notion of an isoparametric submanifold in the sense of Terng) and a curvature-adapted submaifold 
in a Riemannian Hilbert manifold (see Section 2), where ``CSJ-submanfold'' abbreviate ``a submanifold with compact shape operators and compact normal Jacobi operators''.  

\vspace{0.25truecm}

The first result in this paper is as follows.  

\vspace{0.25truecm}

\noindent
{\bf Theorem A.} {\sl (i)\ \ If $M$ be a curvature-adapted CSJ-submanifold with flat section in a locally symmetric Riemannian Hilbert manifold 
$(\widetilde M,\widetilde g)$, then it is proper Fredholm.  

(ii)\ \ Let $M$ be as in the statement (i).  Furthermore, assume that it is a regularizable submanifold with trivial normal holonomy group and that it is of bounded curvatures.  
If $M$ is weakly isoparametric, then it is isoparametric.  
}

\vspace{0.25truecm}

See Section 4 about standard examples of curvature-adapted isoparametric submanifolds with bounded curvatures in a locally symmetric Riemannian Hilbert manifold.  

Let $\pi:P\to B$ be a (smooth) $G$-bundle over a compact Riemannian manifold $(B,g_B)$, where $G$ is a compact semi-simple Lie group.  
Denote by $\mathcal A_P^{H^s}$ the (affine separable) Hilbert space of all $H^s$-connections of $P$ and $\mathcal G_P^{H^{s+1}}$ the $H^{s+1}$-gauge transformation group of $P$, 
where we assume that $\displaystyle{s>\frac{1}{2}\,{\rm dim}\,B-1}$.  
Then, according to Section 9 of \cite{P}, $\mathcal G_P^{H^{s+1}}$ is a smooth Hilbert Lie group.  
Also, according to Lemma 1.2 of \cite{U}, the gauge action $\mathcal G_P^{H^{s+1}}\curvearrowright\mathcal A_P^{H^s}$ is smooth.  
For $\omega\in\mathcal A_P^{H^s}$, we define an operator $\square_{\omega}:\Omega_{\mathcal T,1}^{H^s}(P,\mathfrak g)\to\Omega_{\mathcal T,1}^{H^{s-2}}(P,\mathfrak g)$ by 
$$\square_{\omega}:=d_{\omega}\circ d_{\omega}^{\ast}+d_{\omega}^{\ast}\circ d_{\omega}+{\rm id},$$
where $\Omega_{\mathcal T,1}^{H^j}(P,\mathfrak g)$ ($j=s,s-2$) denotes the Hilbert space of all $\mathfrak g$-valued tensorial $1$-form of class $H^j$ on $P$, 
$d_{\omega}$ denotes the covariant exterior derivative with respect to $\omega$ and $d_{\omega}^{\ast}$ denotes the adjoint operator of $d_{\omega}$.  
The $L^2_s$-inner product $\langle\,\,,\,\,\rangle_s^{\omega}$ of $\Omega_{\mathcal T,1}^{H^s}(P,\mathfrak g)$ is defined by using $\square_{\omega}$ (see Section 6).  
The gauge transformation group $\mathcal G_P^{H^{s+1}}$ does not act isometrically on the (affine) Hilbert space $(\mathcal A_P^{H^s},\langle\,\,,\,\,\rangle^{\omega_0}_s)$.  
Define a Riemannian metric ${\it g}_s$ on $\mathcal A_P^{H^s}$ by 
$({\it g}_s)_{\omega}:=\langle\,\,,\,\,\,\rangle^{\omega}_s$ ($\omega\in\mathcal A_P^{H^s}$).  
This Riemannian metric ${\it g}_s$ is a non-flat metric.  
Note that the Hilbert space $(\mathcal A_P^{H^s},\langle\,\,,\,\,\rangle^{\omega}_s)(=(\Omega_{\mathcal T,1}^{H^s}(P,\mathfrak g),\langle\,\,,\,\,\rangle_s^{\omega}))$ 
is regarded as the tangent space of $(\mathcal A_P^{H^s},{\it g}_s)$ at $\omega$ equipped with the inner product $({\it g}_s)_{\omega}$.  
In \cite{K2}, for a $C^{\infty}$-loop $c:[0,1]\to B$ of nonzero constant speed, we originally defined the notion of a pull-back connection map 
$\mu_c:\mathcal A_P^{H^s}\to H^s([0,1],\mathfrak g)$ (see \cite{K3} also).  
In this paper, for a general $C^{\infty}$-curve $c:[0,1]\to B$ of nonzero constant speed, we define this notion (see Section 6).  
This map is a bounded linear operator of $(\mathcal A_P^{H^s},g_s)$ onto $H^s([0,1],\mathfrak g)$.  
The additive group $\mu_c^{-1}(\hat{\bf 0})$ acts on $\mathcal A_P^{H^s}$ as parallel translations, 
where $\hat{\bf 0}$ is the constant path at the zero vector ${\bf 0}$ of $\mathfrak g$ and 
$\mu_c^{-1}(\hat{\bf 0})$ is regarded as a closed subgroup of the additive (Hilbert Lie) group 
$\Omega_{\mathcal T,1}^{H^s}(P,\mathfrak g)$ acting on $\mathcal A_P^{H^s}$ as parallel translations.  
Denote by $\widetilde{\mathcal G_P^{H^{s+1}}}^c$ the semi-direct group $\mathcal G_P^{H^{s+1}}\ltimes\mu_c^{-1}(\hat{\bf 0})$ of 
$\mathcal G_P^{H^s}$ and $\mu_c^{-1}(\hat{\bf 0})$.  
In this paper, we shall introduce the notion of a horizontally isometric action on a Riemannian Hilbert manifold (see Section 3).  
The second result in this paper is as follows.  

\vspace{0.25truecm}

\noindent
{\bf Theorem B.} {\sl Let $s>\frac{1}{2}\,{\rm dim}\,B-1$.  Then the action $\widetilde{\mathcal G_P^{H^{s+1}}}^c\curvearrowright(\mathcal A_P^{H^s},g_s)$ is horizontally isometric.  
If $c$ is not a loop, then the action $\widetilde{\mathcal G_P^{H^{s+1}}}^c\curvearrowright(\mathcal A_P^{H^s},g_s)$ is transitive.  
}

\vspace{0.25truecm}



The correspondence $\square_{\bullet}:\omega\mapsto\square_{\omega}$ ($\omega\in\mathcal A_P^{H^s}$) is regarded as a $\Omega_{\mathcal T,1}^{H^{s-2}}(P,\mathfrak g)$-valued 
$1$-form on $\mathcal A_P^{H^s}$ by identifying the source space $\Omega_{\mathcal T,1}^{H^s}(P,\mathfrak g)$ of $\square_{\omega}$ with $T_{\omega}\mathcal A_P^{H^s}$.  
Since $g_s$ is defined by using $\square_{\omega}$'s ($\omega\in\mathcal A_P^{H^s}$), it is very important to investigate whether this 
$\Omega_{\mathcal T,1}^{H^{s-2}}(P,\mathfrak g)$-valued $1$-form $\square_{\bullet}$ on $\mathcal A_P^{H^s}$ is closed.  
The third result in this paper is as follows.  

\vspace{0.25truecm}

\noindent
{\bf Theorem C.} {\sl {\rm (i)}\ \ Define a subspace $W_{\omega}$ of the tangent space $T_{\omega}\mathcal A_P^{H^s}$ by 
$$W_{\omega}:=\{A\in\Omega_{\mathcal T,1}^{H^s}(P,\mathfrak g)\,|\,\nabla^{S,\omega}A={\bf 0}\}\,(\subset\Omega_{\mathcal T,1}^{H^s}(P,\mathfrak g)\approx T_{\omega}\mathcal A_P^{H^s}),$$
where $\nabla^{S,\omega}$ denotes the Riemannian connection of the Sasaki metric $g_S^{\omega}$ of $P$ associated to $\omega$.  
Then, for the $\Omega_{\mathcal T,1}^{H^{s-2}}(P,\mathfrak g)$-valued $1$-form $\square_{\bullet}$ on $\mathcal A_P^{H^s}$, 
$(d\square_{\bullet})_{\omega}|_{W_{\omega}\times W_{\omega}}={\bf 0}$ holds.  

{\rm (ii)}\ \ The $\Omega_{\mathcal T,1}^{H^{s-2}}(P,\mathfrak g)$-valued $1$-form $\square_{\bullet}$ on $\mathcal A_P^{H^s}$ is not closed.}

\vspace{0.25truecm}

\noindent
{\it Remark 1.1.}\ \ The exterior derivative $d\square_{\bullet}$ of the $\Omega_{\mathcal T,1}^{H^s}(P,\mathfrak g)$-valued $1$-form $\square_{\bullet}$ 
on $\mathcal A_P^{H^s}$ is given by 
$$\begin{array}{l}
\displaystyle{(d\square_{\bullet})_{\omega}(A_1,A_2)=
-\sum_{\sigma\in S_2}{\rm sgn}\,\sigma\,{\rm Tr}_{g_S^{\omega}}^{\bullet}\left[(\nabla^{S,\omega}A_{\sigma(1)})\otimes A_{\sigma(2)}\right](\bullet,\bullet,\cdot)}\\
\hspace{3truecm}\displaystyle{
-2\sum_{\sigma\in S_2}{\rm sgn}\,\sigma\,{\rm Tr}_{g_S^{\omega}}^{\bullet}
\left[(\nabla^{S,\omega}A_{\sigma(1)})\otimes A_{\sigma(2)}\right](\cdot,\bullet,\bullet)}\\
\hspace{4.25truecm}\displaystyle{(\omega\in\mathcal A_P^{H^s},\,\,A_1,\,A_2\in T_{\omega}\mathcal A_P^{H^s}(\approx\Omega_{\mathcal T,1}^{H^s}(P,\mathfrak g)))}
\end{array}$$
(see Lemma 11.3), where $S_2$ denotes the symmetric group of degree two and ${\rm sgn}\,\sigma$ denotes the signature of $\sigma$ and 
$\nabla^{S,\omega}$ denotes the Riemannian connection of the Sasaki metric of $P$ defined by the Riemannian metric of $B$, 
the ${\rm Ad}(G)$-invariant inner product of $\mathfrak g$ and the connection $\omega$.  By using this relation, we can derive the non-closedness of $\square_{\bullet}$.  

\vspace{0.25truecm}

In this paper, we shall introduce the notions of a horizontally hyperpolar action, a CSJ-action and a regularizable action on a Riemannian Hilbert manifold (see Section 3).  
Let $K_i$ ($i=1,2$) be symmetric subgroups of $G$.  Give $G/K_i$ ($i=1,2$) the Riemannian metrics $\bar g_i$ induced from the bi-invariant metric of $G$.  
Then $(G/K_i,\bar g_i)$ ($i=1,2$) are Riemannian symmetric spaces of compact type.  
Define a closed subgroup $P(G,K_1\times K_2)$ of $H^{s+1}([0,1],G)$ by 
$$P(G,K_1\times K_2):=\{{\bf g}\in H^{s+1}([0,1],G)\,|\,({\bf g}(0),{\bf g}(1))\in K_1\times K_2\}.$$
Let $\pi:P\to B,\,\,c$ and $\sigma$ be as above.  
Also, define a closed subgroup $(\mathcal G_P^{H^{s+1}})_{K_1,K_2}$ of $\mathcal G_P^{H^{s+1}}$ by 
$$(\mathcal G_P^{H^{s+1}})_{K_1,K_2}:=\{{\bf g}\in\mathcal G_P^{H^{s+1}}\,|\,((\widehat{\bf g}\circ\sigma)(0),(\widehat{\bf g}\circ\sigma)(1))\in K_1\times K_2\},$$
where $\widehat{\bf g}$ is the element of $H^{s+1}(P,G)$ defined by ${\bf g}(u)=u\widehat{\bf g}(u)$ ($u\in P$).  
Also, denote by $(\mathcal G_P^{H^{s+1}})_{K_1,K_2;c}$ the semi-direct product $(\mathcal G_P^{H^{s+1}})_{K_1,K_2}\ltimes\mu_c^{-1}(\hat{\bf 0})$.  
The main result of this paper is as follows.  

\vspace{0.25truecm}

\noindent
{\bf Theorem D.} {\sl 
If $c$ is not a loop, then the action $(\mathcal G_P^{H^{s+1}})_{K_1,K_2;c}\curvearrowright(\mathcal A_P^{H^s},g_s)$ is horizontally hyperpolar.}

\vspace{0.25truecm}

Let $M$ be a compact (embedded) submanifold in a symmetric space.  
C. L. Terng and G. Thorbergsson (\cite{TT}) called $M$ an {\it equifocal submanifold} if it satisfies the following conditions:

\vspace{0.15truecm}

\noindent
(E-i)\ \ The normal holonomy group of $M$ is trivial;

\noindent
(E-ii)\ \ For any parallel normal vector field $\widetilde{\xi}$ of $M$, the end-point map $\eta_{\widetilde{\xi}}$ for $\widetilde{\xi}$ has 

constant rank;

\noindent
(E-iii)\ \ $M$ has flat sections.

\vspace{0.15truecm}

\noindent
Under the assumptions of (E-i) and (E-iii), the condition (E-ii) is equivalent to the following condition:

\vspace{0.15truecm}

\noindent
(E-ii')\ \ For any parallel unit normal vector field $\widetilde{\xi}$ of $M$, the set of all focal radii 

along the normal geodesic $\gamma_{\widetilde{\xi}_x}$ is independent of the choice of $x\in M$.  

\vspace{0.15truecm}

\noindent
Here a {\it focal radius along $\gamma_{\widetilde{\xi}_x}$} means a real number $r$ such that the rank of $(\eta_{r\widetilde{\xi}})_{\ast x}$ is smaller than the dimenison of $M$ 
(equivalently, ${\rm Ker}(\eta_{r\widetilde{\xi}})_{\ast x}\not=\{{\bf 0}\}$).  
In this paper, we define the notion of an equifocal submanifold in a Riemannian Hilbert manifold similarly (see Section 2).  

For a horizontally hyperpolar action, a hyperpolar CSJ-action and a hyperpolar regluarizable action, we can prove the following fact.  

\vspace{0.35truecm}

\noindent
{\bf Theorem E.} {\sl 
(i)\ \ Principal orbits of a horizontally hyperpolar action on a Riemannian Hilbert manifold are equifocal.  

(ii)\ \ Principal orbits of a hyperpolar CSJ-action on a Riemannian Hilbert manifold are weakly isoparametric.  

(iii)\ \ If a principal orbit of a hyperpolar regularizable action on a locally symmetric Riemannian Hilbert manifold is curvature-adapted and of bounded curvature, 
then the principal orbit is isoparametric.  
}

\vspace{0.35truecm}

From Theorem D and the statement (i) of Theorem E, we can derive the following fact.  

\vspace{0.35truecm}

\noindent
{\bf Corollary F.} {\sl The principal orbits of the action $(\mathcal G_P^{H^{s+1}})_{K_1,K_2;c}\curvearrowright(\mathcal A_P^{H^s},g_s)$ are equifocal.}


\section{Isoparametric submanifolds in a Riemannian Hilbert manifold} 
In this section, we shall introduce the notion of an isoparametric submanifold in a Riemannian Hilbert manifold.  
Let $M$ be an immersed submanifold of finite codimension in a Riemannian Hilbert manifold $(\widetilde M,\widetilde g)$.  
We call $M$ a {\it proper submanifold} if it satisfies the condition (P) stated in Introduction.  
Also, we call a {\it proper Fredholm} if it satisies both conditions (P) and (F) stated in Introduction.  
Proper Fredholm submanifolds have good focal structure.  However, since the ambient space is of curvature, the focal structure of $M$ depends not only on the datas 
of the shape operators of $M$ but also on that of the curvature of $(\widetilde M,\widetilde g)$.  Hence, even if $M$ is a proper Fredholm submanifold, 
its shape operators are not necessarily compact.  
Denote by $\widetilde{\nabla}$ and $\widetilde R$ the Riemannian connection and the curvature tensor of $(\widetilde M,\widetilde g)$, respectively.  
Let $M$ be an immersed submanifold in $(\widetilde M,\widetilde g)$.  
For a nonzero normal vector $\xi$ of at $p$, if the operator $\widetilde R(\cdot,\xi)\xi$ preserves the tangent space $T_pM$ invariantly, then 
the restriction $\widetilde R(\xi):=\widetilde R(\cdot,\xi)\xi|_{T_pM}$ is called a {\it normal Jacobi operator for} $\xi$, which is is a self-adjoint operator of $(T_pM,g_p)$, where 
$g$ denotes the induced metric on $M$.  

\vspace{0.25truecm}

\noindent
{\bf Definition 2.1.}\ \ Let $M$ be a proper submanifold in $(\widetilde M,\widetilde g)$.  

{\rm (i)}\ If $\widetilde R(\cdot,\xi)\xi$ preserves the tangent space $T_pM$ invariantly for any $p\in M$ and any $\xi\in T_p^{\perp}M$, then 
we call $M$ a {\it curvature-invariant submanifold}.  

{\rm(ii)}\ If $M$ is curvature-invariant and if the shape operator $A_{\xi}$ and the normal Jacobi operator $\widetilde R(\xi)$ are compact for any normal vector $\xi$ of $M$, 
then we call $M$ a {\it submanifold with compact shape operators and compact normal Jacobi operators} (abbreviately, {\it CSJ-submanifold}).  

{\rm(iii)}\ If $M$ is a CSJ-submanifold and if the shape operator $A_{\xi}$ and the normal Jacobi operator $\widetilde R(\xi)$ are regularizable for any normal vector $\xi$ 
of $M$, then we call $M$ a {\it regularizable submanifold}.  

{\rm(iv)}\ If $M$ is curvature-invariant and if, for any normal vector $\xi$ of $M$, the shape operator $A_{\xi}$ and the normal Jacobi operator $\widetilde R(\xi)$ 
commutes, then we call $M$ a {\it curvature-adapted submanifold}.  

{\rm(v)}\ If $M$ is a CSJ-submanifold and if $\displaystyle{\mathop{\sup}_{\xi\in S^{\perp}M}\,\|A_{\xi}\|_{\rm op}}$ and 
$\displaystyle{\mathop{\sup}_{\xi\in S^{\perp}M}\,\|\widetilde R(\xi)\|_{\rm op}}$ are finite, then we say that $M$ {\it is of bounded curvatures}, 
where $S^{\perp}M$ denotes the spherical normal bundle consisting of unit normal vectors of $M$.  

\vspace{0.35truecm}

\noindent
{\bf Definition 2.2.}\ \ Let $M$ be a CSJ-submanifold in $(\widetilde M,\widetilde g)$.  
We call $M$ a {\it weakly isoparametric submanifold} if it satisfies the following conditions:

(WI-i)\ \ The normal holonomy group of $M$ is trivial;

(WI-ii)\ \ For any parallel normal vector field $\widetilde{\xi}$ of $M$, the shape operators $A_{\widetilde{\xi}_x}$ and $A_{\widetilde{\xi}_y}$ are orthogonally equivalent, 
where $x$ and $y$ are any points of $M$;

(WI-iii)\ \ For any parallel normal vector field $\widetilde{\xi}$ of $M$, the normal Jacobi operators $\widetilde R(\widetilde{\xi}_x)$ and $\widetilde R(\widetilde{\xi}_y)$ 
are orthogonally equivalent, where $x$ and $y$ are any points of $M$;

(WI-iv)\ \ $M$ has flat sections.  

\vspace{0.35truecm}

\noindent
{\bf Definition 2.3.}\ \ Let $M$ be a CSJ-submanifold in $(\widetilde M,\widetilde g)$.  
We call $M$ a {\it isoparametric submanifold} if it satisfies the conditions (WI-i), (WI-iv) and the following condition:

(I)\ \ $M$ is a regularizable submanifold with bounded curvatures and the parallel 

submanifolds sufficiently close to $M$ are regularizable submanifolds of constant 

regularized mean curvature in the radial direction.  

\vspace{0.15truecm}

\noindent
Here we note that there exist parallel submanifolds sufficiently close to $M$ because $M$ is of bounded curvature.  

\vspace{0.35truecm}

\noindent
{\bf Definition 2.4.}\ \ Let $M$ be an immersed submanifold of finite codimension in $(\widetilde M,\widetilde g)$.  
We call $M$ an {\it equifocal submanifold} if it satisfies the following conditions:

(${\rm E}_{\infty}$-i)\ \ The normal holonomy group of $M$ is trivial;

(${\rm E}_{\infty}$-ii)\ \ For any parallel normal vector field $\widetilde{\xi}$ of $M$, ${\rm dim}\,{\rm Ker}(\eta_{\widetilde{\xi}})_{\ast x}$ is finite and 

independent of the choice of $x\in M$;

(${\rm E}_{\infty}$-iii)\ \ $M$ has flat sections.

\vspace{0.35truecm}

For $\xi\in T_x^{\perp}M$, denote by $\gamma_{\xi}$ the normal geodesic of direction $\xi$ (i.e., $\gamma_{\xi}'(0)=\xi$).  
Let $\vY$ be a Jacobi field along $\gamma_{\xi}$.  If $\vY(0)$ is a nonzero vector of $df_x(T_xM)$, then it is called a {\it $M$-Jacobi field along} $\gamma_{\xi}$.  
Furthermore, if $\vY'(0)\in df_x(T_xM)$, then it is called a {\it strongly $M$-Jacobi field along} $\gamma_{\xi}$.  
If there exists a $M$-Jacobi field $\vY$ along $\gamma_{\xi}$ with $\vY(r)={\bf 0}$, then $r$ is called a {\it focal radius along} $\gamma_{\xi}$ and 
$\gamma_{\xi}(r)(=\exp^{\perp}(r\xi))$ is called a {\it focal point along} $\gamma_{\xi}$.  Also, the dimension of the linear subspace 
$$\{\vY(0)\,|\,\vY\,:\,M-{\rm Jacobi}\,\,{\rm field}\,\,{\rm along}\,\,\gamma_{\xi}\,\,{\rm with}\,\,\vY(r)={\bf 0}\}$$
of $df_x(T_xM)$ is called the {\it multiplicity} of the focal radius $r$.  
In the case where $(\widetilde M,\widetilde g)$ is locally symmetric and $M$ has trivial normal holonomy group and flat section, 
$(\eta_{r\widetilde{\xi}})_{\ast x}(\vY(0))=\vY(r)$ holds (see Section 5).  
Hence, in this case, the above condition (${\rm E}_{\infty}$-ii) is equivalent to 
the following condition:

\vspace{0.15truecm}

(${\rm E}_{\infty}$-ii$'$)\ \ For any parallel normal vector field $\widetilde{\xi}$ of $M$ and any $x\in M$, the 

multiplicity of each focal radius along $\gamma_{\widetilde{\xi}_x}$ is finite and the set of all focal radii 

along $\gamma_{\widetilde{\xi}_x}$ is independent of the choice of $x\in M$ with considering the multiplicities.  

\section{Horizontally hyperpolar actions on a Rimannian Hilbert manifold} 
In this section, we shall introduce the notion of a (horizontally) hyperpolar action on a Riemannian Hilbert manifold.  
Let $\mathcal G\curvearrowright(\widetilde M,\widetilde g)$ be a proper action of a Hilbert Lie group $\mathcal G$ on a Riemannian Hilbert manifold $(\widetilde M,\widetilde g)$.  

\vspace{0.25truecm}

\noindent
{\bf Definition 3.1.}\ \ {\rm (i)}\ If $\mathcal G\curvearrowright(\widetilde M,\widetilde g)$ is an isometric action of finite cohomogeneity and if there exists a finite dimensional 
complete totally geodesic submanifold $\Sigma$ (in $(\widetilde M,\widetilde g)$) meeting to all $\mathcal G$-orbits orthogonally, then we call this action a {\it polar action} 
and $\Sigma$ its {\it section}.  Furthermore, if the induced metric on $\Sigma$ is flat, then we call this action a {\it hyperpolar action} and $\Sigma$ its {\it flat section}.  

{\rm (ii)}\ If $\mathcal G\curvearrowright(\widetilde M,\widetilde g)$ is an isometric action of finite cohomogeneity and if the principal orbits of this action are 
CSJ-submanifolds (resp. regularizable submanifolds), then we call this action a {\it CSJ-action} (resp. a {\it regularizable action}).  

{\rm (iii)}\ If, for any ${\bf g}\in\mathcal G$ and any regular point $u(\in\widetilde M)$ of the action $\mathcal G\curvearrowright(\widetilde M,\widetilde g)$, 
$${\rm pr}_{T^{\perp}_{{\bf g}\cdot u}(\mathcal G\cdot u)}\circ{\bf g}_{\ast u}|_{T_u^{\perp}(\mathcal G\cdot u)}$$
is a linear isometry, then we call this action a {\it horizontally isometric action}, where ${\rm pr}_{T^{\perp}_{{\bf g}\cdot u}(\mathcal G\cdot u)}$ denotes 
the orthogonal projection of $T_{{\bf g}\cdot u}\widetilde M$ onto $T^{\perp}_{{\bf g}\cdot u}(\mathcal G\cdot u)$.  

{\rm (iv)}\ If $\mathcal G\curvearrowright(\widetilde M,\widetilde g)$ is a horizontally isometric action of finite cohomogeneity and if there exists a finite dimensional complete 
totally geodesic submanifold $\Sigma$ (in $(\widetilde M,\widetilde g)$) meeting to all $\mathcal G$-orbits orthogonally, then we call this action a {\it horizontally polar action} 
and $\Sigma$ its {\it section}.  Furthermore, if the induced metric on $\Sigma$ is flat, then we call this action a {\it horizontally hyerpolar action} and $\Sigma$ 
its {\it flat section}.  

\vspace{0.25truecm}

\noindent
{\it Remark 3.1.}\ \ (i)\ Let $\mathcal U$ be the set of all regular points of a proper action $\mathcal G\curvearrowright(\widetilde M,\widetilde g)$.  If this action is 
horizontally isometric, then there exists the Riemannian metric $\bar g$ on the orbit space $\mathcal U/\mathcal G$ such that the restriction of the orbit map 
$\pi_{\mathcal G}$ of this action to $\mathcal U$ is a Riemannin submersion of $(\mathcal U,g)$ onto $(\mathcal U/\mathcal G,\bar g)$.  
In particular, if this action is horizontally polar, then the horizontal distribution of the Riemannin submersion 
$\pi_{\mathcal G}|_{\mathcal U}:(\mathcal U,g)\to(\mathcal U/\mathcal G,\bar g)$ is integrable.  From this fact, we can show that, for any $\xi\in T_u\widetilde M$, 
the $\mathcal G$-equivariant normal vector field 
$$\widetilde{\xi}:{\bf g}\cdot u\mapsto({\rm pr}_{T^{\perp}_{{\bf g}\cdot u}(\mathcal G\cdot u)}\circ{\bf g}_{\ast u})(\xi)\quad\,\,({\bf g}\in\mathcal G))$$
of the orbit $\mathcal G\cdot u$ is parallel with respect to the normal connection of $\mathcal G\cdot u$.  This fact follows from the fact that the integral manifolds of 
the horizontal distribution are totally geodesic.  

\vspace{0.25truecm}

\section{Examples.} 
In this section, we shall give standard examples of curvature-adapted isoparametric submanifolds in a locally symmetric Riemannian Hilbert maifold.  

\vspace{0.25truecm}

\noindent{\it Example 4.1.}\ \ (i)\ Let $(V,\langle\,\,,\,\,\rangle)$ be a (separable) Hilbert space and $\{\ve_i\}_{i=1}^{\infty}$ be an orthonormal basis of 
$(V,\langle\,\,,\,\,\rangle)$.  For a bounded sequence $\{m_i\}_{i=1}^{\infty}$ of positive integers and a sequence $\{r_i\}_{i=1}^{\infty}$ of positive numbers 
satisfying $\sum\limits_{i=1}^{\infty}r_i^2<\infty$, define a Hilbert submanifold $\widetilde M(\{m_i\}_{i=1}^{\infty},\{r_i\}_{i=1}^{\infty})$ by 
$$\begin{array}{l}
\hspace{0.5truecm}\displaystyle{\widetilde M(\{m_i\}_{i=1}^{\infty},\{r_i\}_{i=1}^{\infty})}\\
\displaystyle{:=\left\{\left.\sum_{i=1}^{\infty}a_i\ve_{2i}+\sum_{j=1}^{\infty}b_j\ve_{2j-1}\,\right\vert\,
\sum_{i=1}^{m_1}a_i^2=r_1^2,\,\,\,\sum_{i=m_1+\cdots+m_k+1}^{m_1+\cdots+m_{k+1}}a_i^2=r_{k+1}^2\,\,\,(k\geq 1),\right.}\\
\hspace{5truecm}\displaystyle{\left.\sum_{i=1}^{\infty}b_i^2<\infty\right\}.}
\end{array}$$
Denote by $\widetilde g$ the induced metric on $\widetilde M(\{m_i\}_{i=1}^{\infty},\{r_i\}_{i=1}^{\infty})$.  
Then it is easy to show that $(\widetilde M(\{m_i\}_{i=1}^{\infty},\{r_i\}_{i=1}^{\infty}),\,\,\widetilde g)$ is a locally symmetric.  
Take any unit vector $\vv$ of $V$, which is expressed as the Fourier series $\vv=\sum\limits_i\langle\vv,\ve_i\rangle\ve_i$.  
For the simplicity, set $\alpha_i:=\langle\vv,\ve_i\rangle$ and 
$$I_k:=\left\{\begin{array}{ll}
\{1,\cdots,m_1\} & (k=1)\\
\{m_1+\cdots+m_{k-1}+1,\cdots,m_1+\cdots+m_k\} & (k\geq 2).
\end{array}\right.$$
Set $\vv_k:=\sum\limits_{i\in I_k}\alpha_{2i}\ve_{2i}$ ($i\in\mathbb N$).  
Take another $\vw\in V$ and set $\beta_i:=\langle\vw,\ve_i\rangle$ and $\vw_k:=\sum\limits_{i\in I_k}\beta_{2i}\ve_{2i}$ ($i\in\mathbb N$).  
Then we have 
$$\widetilde R(\vw,\vv)\vv=\sum_{k=1}^{\infty}\widetilde R(\vw_k,\vv_k)\vv_k=\sum_{k=1}^{\infty}\frac{1}{r_k^2}\,\|\vv_k\|^2\vw_k.\leqno{(4.1)}$$
Set $I_{\vv}:=\{i\in\vN\,|\,\vv_i\not={\bf 0}\}$.  From $(4.1)$, we see that, if $\sharp\,I_{\vv}<\infty$, then $\widetilde R(\cdot,\vv)\vv$ is a compact operator 
but, if $\sharp\,I_{\vv}=\infty$, then $\widetilde R(\cdot,\vv)\vv$ is not necessarily a compact operator.  
Define a subset $M$ of $\widetilde M(\{m_i\}_{i=1}^{\infty},\{r_i\}_{i=1}^{\infty})$ by 
\begin{align*}
M:=\{\sum_{i=1}^{\infty}&a_i\ve_{2i}+\sum_{j=1}^{\infty}b_j\ve_{2j-1}\,\in\,\widetilde M(\{m_i\}_{i=1}^{\infty},\{r_i\}_{i=1}^{\infty})\,|\,\\
&a_{m_1+\cdots+m_j}=r'_j\,\,(j=1,\cdots,k_1),\,\,b_j=0\,\,(j=1,\cdots,k_2)\},
\end{align*}
where $r'_j$ ($j=1,\cdots,k_1$) are positive numbers smaller than $r_j$.  
Take any $x\in M$ and any $\xi\in T_x^{\perp}M$.  Set 
$$(E_j)_x:=\{\sum_{i\in I_j}a_i\ve_{2j}\,|\,a_i\in\mathbb R\}\cap T_xM\quad\,\,(j=1,\cdots,k_1)$$
and 
$$(E_0)_x:=\left(\mathop{\oplus}_{j=1}^{k_1}(E_j)_x\right)^{\perp}\cap T_xM.$$
It is clear that 
$$T_xM=\left(\mathop{\oplus}_{j=1}^{k_1}(E_j)_x\right)\oplus(E_0)_x$$
holds and that the correspondences $x\mapsto(E_j)_x$ ($j=1,\cdots,k_1$) and $x\mapsto(E_0)_x$ give $C^{\infty}$-distibutios on $M$.  
Since 
$$T_x^{\perp}M\subset\left(\mathop{\oplus}_{j=1}^{k_1}\{\sum_{i\in I_j}a_i\ve_{2_i}\,|\,a_i\in\mathbb R\}\right)\oplus
{\rm Span}\{\ve_1,\ve_3,\cdots,\ve_{2k_2-1}\},$$
we have 
$$\xi=\sum_{j=1}^{k_1}\sum_{i\in I_j}\langle\xi,\ve_{2i}\rangle\ve_{2i}+\sum_{i=1}^{k_2}\langle\xi,\ve_{2i-1}\rangle\ve_{2i-1}.$$
Set $\xi_j:=\sum\limits_{i\in I_j}\langle\xi,\ve_{2i}\rangle\,\ve_{2i}$ ($j=1,\cdots,k_1$) and $\xi_0:=\sum\limits_{i=1}^{k_2}\langle\xi,\ve_{2i-1}\rangle\ve_{2i-1}$. 
Then we obtain 
$$\widetilde R_x(\cdot,\xi)\xi|_{(E_j)_x}=\frac{1}{r_j^2}\cdot\|\xi_j\|^2\cdot{\rm id},\,\,\,\,\,\,\widetilde R_x(\cdot,\xi)\xi|_{(E_0)_x}={\bf 0}$$
and 
$$A_{\xi}|_{(E_j)_x}=\sqrt{\frac{1}{{r'}_j^2}-\frac{1}{r_j^2}}\cdot\|\xi_j\|\cdot{\rm id},\,\,\,\,\,\,A_{\xi}|_{(E_0)_x}={\bf 0}.$$
Thereforethe decomposition 
$$T_xM=\left(\mathop{\oplus}_{j=1}^{k_1}(E_j)_x\right)\oplus(E_0)_x$$
gives the simultaneously eigenspace decomposition of the commutative family 
$$\{\widetilde R(\cdot,\xi)\xi\,|\,\xi\in T_x^{\perp}M\}\cup\{A_{\xi}\,|\,\xi\in T_x^{\perp}M\}$$
of self-adjoint operators.  Thus $M$ is a curvature-adapted CSJ-submanifold.  
Also, it is easy to show that $M$ has trivial normal holonomy group and flat sections.  
Since $A_{\xi}$ is of finite rank, the regularized trace of ${\rm Tr}_rA_{\xi}$ is equal to the usual trace of $A_{\xi}$, which is equal to 
$\displaystyle{\sum_{j=1}^{k_1}\sqrt{\frac{1}{{r'}_j^2}-\frac{1}{r_j^2}}\cdot\|\xi_j\|\cdot(m_j-1)}$.  
Thus $M$ is regularizble.  Furthermore, it is easy to show that $M$ satisfies the condition (I) in Definition 2.3.  
Thus $M$ is isoparametric.  

\section{Proof of Theorem A} 
In this section, we prove Theorem A stated in Introduction.  
Let $M$ be a CSJ-submanifold with flat section in a locally symmetric Riemannian Hilbert manifold $(\widetilde M,\widetilde g)$ immersed by $f$.  
Then the focal structure of $M$ is determined by the shape operators $A_{\xi}$'s ($\xi\in T^{\perp}M$) and 
the normal Jacobi operators $\widetilde R(\xi)$'s ($\xi\in T^{\perp}M$).  
In fact, since $(\widetilde M,\widetilde g)$ is locally symmetric, the $M$-Jacobi field $\vY$ along the normal geodesic $\gamma_{\xi}$ of $\xi$-direction 
($\xi\in T_x^{\perp}M$) is described as 
$$\vY(s)=P_{\gamma_{\xi}|_{[0,s]}}\left(\,\,D^{co}_{s\xi}(\vY(0))+sD^{si}_{s\xi}(\vY'(0))\,\,\right),\leqno{(5.1)}$$
where 
$P_{\gamma_{\xi}|_{[0,s]}}$ denotes the parallel translation along the geodesic segment $\gamma_{\xi}|_{[0,s]}$ with respect to the Riemannian connection $\widetilde{\nabla}$ 
of $\widetilde g$, and $D^{co}_{s\xi}$ and $D^{si}_{s\xi}$ are compact self-adjoint operators of $T_{f(x)}\widetilde M$ onto oneself defined by 
$$D^{co}_{s\xi}:=P_{\gamma_{\xi}|_{[0,s]}}\circ\left(\sum_{k=0}^{\infty}\frac{(-1)^ks^{2k}}{(2k)!}\,\,\widetilde R(\xi)^k\right)\circ P_{\gamma_{\xi}|_{[0,s]}}^{-1}$$
and 
$$D^{si}_{s\xi}:=P_{\gamma_{\xi}|_{[0,s]}}\circ\left(\sum_{k=0}^{\infty}\frac{(-1)^ks^{2k}}{(2k+1)!}\,\,\widetilde R(\xi)^k\right)\circ P_{\gamma_{\xi}|_{[0,s]}}^{-1},$$
respectively.  The relation $(5.1)$ is shown by imitating the proof of Lemma 3.4 in \cite{TT} (see Section 1 of \cite{K1} also).  
Since $M$ has flat section, $D^{co}_{s\xi}$ and $D^{si}_{s\xi}$ preserves $df_x(T_xM)$ and $T^{\perp}_xM$ invariantly, respectively.  
Hence we have 
$$\begin{array}{c}
\displaystyle{(D^{co}_{s\xi}(\vY(0))+sD^{si}_{s\xi}(\vY'(0)))_T=D^{co}_{s\xi}(\vY(0))+sD^{si}_{s\xi}(\vY'(0)_T),}\\
\displaystyle{(D^{co}_{s\xi}(\vY(0))+sD^{si}_{s\xi}(\vY'(0)))_{\perp}=sD^{si}_{s\xi}(\vY'(0)_{\perp})=s\vY'(0)_{\perp},}
\end{array}$$
where $(\cdot)_T$ (resp. $(\cdot)_{\perp}$) denotes is the tangential (resp. the normal) component of $(\cdot)$.  
Hence, if $\vY'(0)_{\perp}\not={\bf 0}$, then there does not exist a nonzero real constant $s$ with $\vY(s)={\bf 0}$.  
On the other hand, if $\vY'(0)_{\perp}={\bf 0}$ (then $\vY$ is called a {\it strongly $M$-Jacobi field along} $\gamma_{\xi}$), then we have 
$\vY'(0)=-df_x(A_{\xi}(df_x^{-1}(\vY(0))))$ and hence 
$$\vY(s)=P_{\gamma_{\xi}|_{[0,s]}}\left(\,\,D^{co}_{s\xi}(\vY(0))-sD^{si}_{s\xi}(df_x(A_{\xi}(df_x^{-1}(\vY(0)))))\,\,\right).\leqno{(5.2)}$$
From thiese facts, we can derive that $r$ is a focal radius along $\gamma_{\xi}$ if and only if there exists a strongly $M$-Jacobi field $\vY$ along $\gamma_{\xi}$ 
with $\vY(r)={\bf 0}$.  
Define a family $\{Q_{\xi}(s)\}_{s\in\mathbb R}$ of linear operators of $T_xM$ into $T_{f(x)}\widetilde M$ by 
$$Q_{\xi}(s):=D^{co}_{s\xi}\circ df_x-sD^{si}_{s\xi}\circ df_x\circ A_{s\xi}.$$
Then we can derive that the set of all focal radii of $M$ along $\gamma_{\xi}$ is equal to 
$$\{s\in\mathbb R\,|\,{\rm Ker}\,Q_{\xi}(s)\not=\{{\bf 0}\}\}.\leqno{(5.3)}$$
Any strongle $M$-Jacobi field along $\gamma_{\xi}$ is given as the variational vector field of the geodesic variation $\delta$ given by $\delta(s,t):=\gamma_{\widetilde{\xi}_t}(s)$, where 
$\widetilde{\xi}_t$ is the parallel normal vector field along a curve $c$ in $M$ with $c'(0)=\vY(0)$ satisfying $\widetilde{\xi}_0=\xi$.  
From this fact, if $M$ has trivial normal holonomy group and flat section, then, for any parallel unit normal vector field $\widetilde{\xi}$ of $M$ 
$(\eta_{s\widetilde{\xi}})_{\ast x}(\vY(0))=\vY(s)$ ($s\in\mathbb R$) holds.  Hence $r$ is a focal radius along $\gamma_{\xi}$ if and only if 
${\rm Ker}(\eta_{r\widetilde{\xi}})_{\ast x}\not=\{{\bf 0}\}$ holds.  

%

\vspace{0.25truecm}

\noindent
{\it Proof of Theorem A.}\ \ First we shall  show the statement (i) of Theorem A.  
Let $M$ be a curvature-adapted CSJ-submanifold with flat section in a locally symmetric Riemannian Hilbert manifold $(\widetilde M,\widetilde g)$.  
Take any unit normal vector $\xi$ of $M$ at $x$.  Since $A_{\xi}$ and $\widetilde R(\xi)$ are compact sef-adjoint operators of $(T_xM,g_x)$, 
the spectrums of $A_{\xi}$ and $\widetilde R(\xi)$ are described as 
$${\rm Spec}\,A_{\xi}=\{0\}\cup\{\lambda^A_i\,\vert\,i\in I_A\}\quad\,\,
(\vert\lambda^A_i\vert>\vert\lambda^A_{i+1}\vert\,\,{\rm or}\,\,\lambda^A_i=-\lambda^A_{i+1}>0)$$
and 
$${\rm Spec}\,\widetilde R(\xi)=\{0\}\cup\{\lambda^R_i\,\vert\,i\in I_R\}\quad\,\,
(\vert\lambda^R_i\vert>\vert\lambda^R_{i+1}\vert\,\,{\rm or}\,\,\lambda^R_i=-\lambda^R_{i+1}>0),$$
respectively, where $I_A$ is equal to $\mathbb N$ or $\{1,\cdots,l_A\}$ for some $l_A\in\mathbb N$, and 
$I_R$ is equal to $\mathbb N$ or $\{1,\cdots,l_R\}$ for some $l_R\in\mathbb N$.  
Also, the eigenspace decompositions of $T_xM$ for $A_{\xi_x}$ and $\widetilde R(\xi_x)$ are described as 
$$T_xM=\overline{{\rm Ker}\,A_{\xi}\oplus\left(\mathop{\oplus}_{i\in I_A}{\rm Ker}(A_{\xi}-\lambda^A_i\,{\rm id})\right)}
\leqno{(5.4)}$$
and 
$$T_xM=\overline{{\rm Ker}\,\widetilde R(\xi)\oplus
\left(\mathop{\oplus}_{i\in I_R}{\rm Ker}(\widetilde R(\xi)-\lambda^R_i\,{\rm id})\right)}.\leqno{(5.5)}$$
For the simplicity, we set 
$$\begin{array}{c}
E^A_0:={\rm Ker}\,A_{\xi},\,\,\,\,E^A_i:={\rm Ker}(A_{\xi}-\lambda^A_i\,{\rm id}),\\
E^R_0:={\rm Ker}\,\widetilde R(\xi),\,\,\,\,E^R_i:={\rm Ker}(\widetilde R(\xi)-\lambda^R_i\,{\rm id}).
\end{array}$$
Since $M$ is curvature-adapted, $A_{\xi}$ and $\widetilde R(\xi)$ are simultaneously diagonalizable.  
Therefore we have their common eigenspace decomposition:
$$\begin{array}{l}
\displaystyle{T_xM=\overline{(E^R_0\cap E^A_0)\oplus\left(\mathop{\oplus}_{i\in I_A}(E^R_0\cap E^A_i)\right)
\oplus\left(\mathop{\oplus}_{i\in I_R}(E^R_i\cap E^A_0)\right)}}\\
\hspace{1.7truecm}\displaystyle{\overline{
\oplus\left(\mathop{\oplus}_{i\in I_R}\mathop{\oplus}_{j\in I_A}(E^R_i\cap E^A_j)\right)}.}
\end{array}\leqno{(5.6)}$$
Set $m_{ij}:={\rm dim}(E^R_i\cap E^A_j)$ ($i\in\{0\}\cup I_R,\,\,j\in\{0\}\cup I_A$).  
Let $\{\ve_{i,j;k}\,\vert\,k=1,\cdots,m_{ij}\}$ ($(i,j)\in((\{0\}\cup I_R)\times(\{0\}\cup I_A))\setminus\{(0,0)\}$) be 
an orthonormal basis of $E^R_i\cap E^A_j$.  Take a curve $\alpha:(-\varepsilon,\varepsilon)\to M$ with $\alpha'(0)=\ve_{i,j;k}$ 
and define a geodesic variation $\delta:[0,r+\varepsilon)\times(-\varepsilon,\varepsilon)\to\widetilde M$ by 
$\delta(s,t):=\gamma_{\widetilde{\xi}_{\alpha(t)}}(s)$, where $\widetilde{\xi}$ denotes the parallel normal vector field along $\alpha$ with $\widetilde{\xi}(0)=\xi$ and 
$\varepsilon$ is a sufficiently small positive number.  
Let $\vY$ be the variational vector field of this geodesic variation $\delta$, that is, 
$\displaystyle{\vY=\left.\frac{\partial\delta}{\partial t}\right|_{t=0}}$.  This variational vector field $\vY$ is a strongly 
$M$-Jacobi field along the normal geodesic $\gamma_{\xi}$ with $\vY(0)=df_x(\ve_{i,j;k})$ and 
$\vY'(0)=-df_x(A_{\widetilde{\xi}_x}(\ve_{i,j;k}))=-\lambda^A_j\,df_x(\ve_{i,j;k})$.  
Hence, from $(5.2)$, we have 
$$\vY(s)=\left(\cos(s\sqrt{\lambda^R_i})-\frac{\sin(s\sqrt{\lambda^R_i})}{\sqrt{\lambda^R_i}}\cdot\lambda^A_j\right)\,
P_{\gamma_{\xi}|_{[0,s]}}\left(f_{\ast x}(\ve_{i,j;k})\right),\leqno{(5.7)}$$
Here, in the case of $\lambda^R_i<0$, $\cos(s\sqrt{\lambda^R_i})$ and $\frac{\sin(s\sqrt{\lambda^R_i})}{\sqrt{\lambda^R_i}}$ means 
$\cosh(s\sqrt{-\lambda^R_i})$ and $\frac{\sinh(s\sqrt{-\lambda^R_i})}{\sqrt{-\lambda^R_i}}$, respectively.  Also, in the case of 
$\lambda^R_i=0$, $\frac{\sin(s\sqrt{\lambda^R_i})}{\sqrt{\lambda^R_i}}$ means $s$.  
From $(5.7)$, we can derive that $\vY(r)={\bf 0}$ holds if and only if 
$$r=\left\{\begin{array}{ll}
\displaystyle{\frac{1}{\sqrt{\lambda^R_i}}\left({\rm arctan}\frac{\sqrt{\lambda^R_i}}{\lambda^A_j}+k\pi\right)\,\,\,\,(k\in\mathbb Z)\,\,\,\,} & (\lambda^R_i>0)\\
\displaystyle{\frac{1}{\sqrt{\lambda^R_i}}\cdot{\rm arctanh}\frac{\sqrt{\lambda^R_i}}{\lambda^A_j}} & (\lambda^R_i<0,\,\,|\lambda^A_j|>|\lambda^R_i|))\\
{\rm does}\,\,{\rm not}\,\,{\rm exist} & (\lambda^R_i<0,\,\,|\lambda^A_j|\leq|\lambda^R_i|))\\
\displaystyle{\frac{1}{\lambda^A_j}} & (\lambda^R_i=0)
\end{array}\right.\leqno{(5.8)}$$
Set 
$$r_{ijk}:=\frac{1}{\sqrt{\lambda^R_i}}\left({\rm arctan}\frac{\sqrt{\lambda^R_i}}{\lambda^A_j}+k\pi\right)\quad\,\,(k\in\mathbb Z)$$
when $\lambda^R_i>0$, 
$$r_{ij}:=\frac{1}{\sqrt{\lambda^R_i}}\cdot{\rm arctanh}\frac{\sqrt{\lambda^R_i}}{\lambda^A_j}$$
when $\lambda^R_i<0$ and $|\lambda^A_j|>|\lambda^R_i|$, and 
$$r_{0j}:=\frac{1}{\lambda^A_j}$$
when $\lambda^R_i=0$.  
Also, define by $I^R_{\pm}$ and $I^A_i$ ($i\in I^R\cup\{0\}$) by 
$$I^R_+:=\{i\in I^R\,|\,\lambda^R_i>0\},\quad\,\,I^R_-:=\{i\in I^R\,|\,\lambda^R_i<0\}$$
and 
$$I^A_i:=\{j\in I^A\,|\,E^R_i\cap E^A_j\not=\{{\bf 0}\}\,\}.$$
Then we see that that the set $\mathcal{FR}_{\xi}$ of all focal radii of $M$ along $\gamma_{\xi}$ is given by 
$$\begin{array}{l}
\displaystyle{\mathcal{FR}_{\xi}=\left(\mathop{\cup}_{i\in I^R_+}\{r_{ijk}\,|\,j\in I^A_i,\,\,k\in\mathbb Z\}\right)}\\
\hspace{1.5truecm}\displaystyle{\cup\left(\mathop{\cup}_{i\in I^R_-}\{r_{ij}\,|\,j\in I^A_i\,\,{\rm s.t.}\,\,|\lambda^A_j|>|\lambda^R_i|\}\right)}\\
\hspace{1.5truecm}\displaystyle{\cup\{r_{0j}\,|\,j\in I^A_0\}.}
\end{array}
\leqno{(5.9)}$$
Since $A_{\xi}$ and $\widetilde R(\xi)$ are compact self-adjoint operators, it follows from $(5.9)$ that the set of all focal radii along $\gamma_{\xi}$ has no accumulating point 
other than $0$ and that the multiplicities of the focal radii is finite.  
Therefore, from the arbitrariness of $\xi$, we see that $M$ is proper Fredholm.  

Next we shall show the statement (ii).  
Let $M$ be a curvature-adapted regularizabe submanifold as in the statement (ii) of Theorem A.  
Take $x\in M$ and a parallel normal vector field $\widetilde{\xi}$ of $M$.  Let $\gamma_{\widetilde{\xi}_x}$ be the normal geodesic of $f(M)$ 
with $\gamma_{\widetilde{\xi}_x}'(0)=\widetilde{\xi}_x$ and set $\widetilde{\xi}^r_x:=\gamma_{\widetilde{\xi}_x}'(r)$.  
Since $M$ is regularizable, $A_{\widetilde{\xi}_x}$ and $\widetilde R(\widetilde{\xi}_x)$ are regularizable compact self-adjoint operators of $(T_xM,g_x)$.  Hence the spectrums of $A_{\widetilde{\xi}_x}$ and $\widetilde R(\widetilde{\xi}_x)$ are described as 
$${\rm Spec}\,A_{\widetilde{\xi}_x}=\{0\}\cup\{\lambda^A_{x,i}\,\vert\,i\in I^x_A\}\quad\,\,
(\vert\lambda^A_{x,i}\vert>\vert\lambda^A_{x,i+1}\vert\,\,{\rm or}\,\,\lambda^A_{x,i}=-\lambda^A_{x,i+1}>0)$$
and 
$${\rm Spec}\,\widetilde R(\widetilde{\xi}_x)=\{0\}\cup\{\lambda^R_{x,i}\,\vert\,i\in I^x_R\}\quad\,\,
(\vert\lambda^R_{x,i}\vert>\vert\lambda^R_{x,i+1}\vert\,\,{\rm or}\,\,\lambda^R_{x,i}=-\lambda^R_{x,i+1}>0),$$
respectively, where $I^x_A$ is equal to $\mathbb N$ or $\{1,\cdots,l_A\}$ for some $l_A\in\mathbb N$, and 
$I^x_R$ is equal to $\mathbb N$ or $\{1,\cdots,l_R\}$ for some $l_R\in\mathbb N$.  
Also, the eigenspace decompositions of $T_xM$ for $A_{\widetilde{\xi}_x}$ and $\widetilde R(\widetilde{\xi}_x)$ are described as 
$$T_xM=\overline{{\rm Ker}\,A_{\widetilde{\xi}_x}\oplus
\left(\mathop{\oplus}_{i\in I^x_A}{\rm Ker}(A_{\widetilde{\xi}_x}-\lambda^A_{x,i}\,{\rm id})\right)}\leqno{(5.10)}$$
and 
$$T_xM=\overline{{\rm Ker}\,\widetilde R(\widetilde{\xi}_x)\oplus
\left(\mathop{\oplus}_{i\in I^x_R}{\rm Ker}(\widetilde R(\widetilde{\xi}_x)-\lambda^R_{x,i}\,{\rm id})\right)}.\leqno{(5.11)}$$
For the simplicity, we set 
$$\begin{array}{c}
E^A_{x,0}:={\rm Ker}\,A_{\widetilde{\xi}_x},\,\,\,\,E^A_{x,i}:={\rm Ker}(A_{\widetilde{\xi}_x}-\lambda^A_{x,i}\,{\rm id}),\\
E^R_{x,0}:={\rm Ker}\,\widetilde R(\widetilde{\xi}_x),\,\,\,\,E^R_{x,i}:={\rm Ker}(\widetilde R(\widetilde{\xi}_x)-\lambda^R_{x,i}\,{\rm id}).
\end{array}$$
Since $M$ is curvature-adapted, $A_{\widetilde{\xi}_x}$ and $\widetilde R(\widetilde{\xi}_x)$ are simultaneously diagonalizable.  
Therefore we have their common eigenspace decomposition:
$$\begin{array}{l}
\displaystyle{T_xM=\overline{(E^R_{x,0}\cap E^A_{x,0})\oplus\left(\mathop{\oplus}_{i\in I^x_A}(E^R_{x,0}\cap E^A_{x,i})\right)
\oplus\left(\mathop{\oplus}_{i\in I^x_R}(E^R_{x,i}\cap E^A_{x,0})\right)}}\\
\hspace{1.7truecm}\displaystyle{\overline{
\oplus\left(\mathop{\oplus}_{i\in I^x_R}\mathop{\oplus}_{j\in I^x_A}(E^R_{x,i}\cap E^A_{x,j})\right)}.}
\end{array}\leqno{(5.12)}$$
For a positive number $r$, denote by $f^r$ the end-point map $\eta_{r\widetilde{\xi}}$ for $r\widetilde{\xi}$ and set $M_{r\widetilde{\xi}}:=f^r(M)$.  
If $r$ is a sufficiently small, then $f^r$ is an immersion because $M$ is of bounded curvatures by the assumption.  
Denote by $A^r$ the shape tensor of the parallel submanifold $M_{r\widetilde{\xi}}$ of $M$.  
Define a normal vector field $\widetilde{\xi}^r$ by $(\widetilde{\xi}^r)_x:=\gamma_{\widetilde{\xi}_x}'(r)$ ($x\in M$).  
Set $m^x_{ij}:={\rm dim}(E^R_{x,i}\cap E^A_{x,j})$ ($i\in\{0\}\cup I^x_R,\,\,j\in\{0\}\cup I^x_A$).  
Let $\{\ve_{i,j;k}\,\vert\,k=1,\cdots,m^x_{ij}\}$ ($(i,j)\in((\{0\}\cup I^x_R)\times(\{0\}\cup I^x_A))\setminus\{(0,0)\}$) be an orthonormal basis of $E^R_{x,i}\cap E^A_{x,j}$.  
Take a curve $\alpha:(-\varepsilon,\varepsilon)\to M$ with $\alpha'(0)=\ve_{i,j;k}$ and define a geodesic variation 
$\delta:[0,r+\varepsilon)\times(-\varepsilon,\varepsilon)\to\widetilde M$ by $\delta(s,t):=\gamma_{\widetilde{\xi}_{\alpha(t)}}(s)$, where $\varepsilon$ is a sufficiently 
small positive number.  Let $\vY$ be the variational vector field of this geodesic variation $\delta$, that is, 
$\displaystyle{\vY=\left.\frac{\partial\delta}{\partial t}\right|_{t=0}}$.  This variational vector field $\vY$ is a strongly $M$-Jacobi field along the normal geodesic 
$\gamma_{\widetilde{\xi}_x}$ with $\vY(0)=df_x(\ve_{i,j;k})$ and $\vY'(0)=-df_x(A_{\widetilde{\xi}_x}(\ve_{i,j;k}))=-\lambda^A_j\,df_x(\ve_{i,j;k})$.  
Hence, from $(5.2)$, we have 
$$\vY(s)=\left(\cos(s\sqrt{\lambda^R_{x,i}})-\frac{\sin(s\sqrt{\lambda^R_{x,i}})}{\sqrt{\lambda^R_{x,i}}}\cdot\lambda^A_{x,j}\right)\,
P_{\gamma_{\widetilde{\xi}_x}|_{[0,s]}}\left(f_{\ast x}(\ve_{i,j;k})\right),\leqno{(5.13)}$$
Here, in the case of $\lambda^R_{x,i}<0$, $\cos(s\sqrt{\lambda^R_{x,i}})$ and $\frac{\sin(s\sqrt{\lambda^R_{x,i}})}{\sqrt{\lambda^R_{x,i}}}$ means 
$\cosh(s\sqrt{-\lambda^R_{x,i}})$ and $\frac{\sinh(s\sqrt{-\lambda^R_{x,i}})}{\sqrt{-\lambda^R_{x,i}}}$, respectively.  Also, in the case of $\lambda^R_{x,i}=0$, 
$\frac{\sin(s\sqrt{\lambda^R_{x,i}})}{\sqrt{\lambda^R_{x,i}}}$ means $s$.  
From $(5.13)$, we have 
$$df^r_x(\ve_{i,j;k})=\vY(r)=\left(\cos(r\sqrt{\lambda^R_{x,i}})-\frac{\sin(r\sqrt{\lambda^R_{x,i}})}{\sqrt{\lambda^R_{x,i}}}\cdot\lambda^A_{x,j}\right)\,
P_{\gamma_{\widetilde{\xi}_x}|_{[0,r]}}(df_x(\ve_{i,j;k}))\leqno{(5.14)}$$
and 
$$\begin{array}{l}
\hspace{0.5truecm}\displaystyle{df^r_x((A^r)_{\widetilde{\xi}^r_x}(\ve_{i,j;k}))=-\vY'(r)}\\
\displaystyle{=\left(\sqrt{\lambda^R_{x,i}}\sin(r\sqrt{\lambda^R_{x,i}})+\lambda^A_j\cos(r\sqrt{\lambda^R_{x,i}})\right)
P_{\gamma_{\widetilde{\xi}_x}\vert_{[0,r]}}(df_x(\ve_{i,j;k})).}
\end{array}\leqno{(5.15)}$$
Hence we obtain 
$$(A^r)_{(\widetilde{\xi}^r)_x}(\ve_{i,j;k})=\frac{\sqrt{\lambda^R_{x,i}}\tan(r\sqrt{\lambda^R_{x,i}})+\lambda^A_{x,j}}
{1-\lambda^A_{x,j}\tan(r\sqrt{\lambda^R_{x,i}})/\sqrt{\lambda^R_{x,i}}}\,\ve_{i,j;k}.\leqno{(5.16)}$$
From $(5.12)$ and $(5.16)$, we can show that $(A^r)_{\widetilde{\xi}^r_x}$ 
is regularizable.  Since $M$ is weakly isoparametric by the assumption, (WI-ii) and (WI-iii) holds.  These together with $(5.12)$ and $(5.16)$ implies that 
${\rm Tr}_r(A^r)_{\widetilde{\xi}^r_x}={\rm Tr}_r(A^r)_{\widetilde{\xi}^r_y}$ holds for any $x,y\in M$, that is, 
$M_{r\widetilde{\xi}}$ is of constant regularized mean curvature in the radial direction.  
This together with the arbitariness of $\widetilde{\xi}$ implies that $M$ is isoparametric.  \qed

\section{The holonomy element map and the pull-back connection map} 
In this section, we shall define the notions of the holonomy element map and the pull-back connection map for a $C^{\infty}$-curve of constant speed in the base manifold 
of a principal bundle and derive some facts for these maps.  
Note that, in \cite {K3}, these maps had already been defined for a $C^{\infty}$-loop of constant speed in base manifold.  
Let $\pi:P\to B$ be a $G$-bundle of class $C^{\infty}$ over a compact Riemannian manifold $(B,g_B)$, where $G$ is a compact semi-simple Lie group.  
Fix an ${\rm Ad}(G)$-invariant inner product $\langle\,\,,\,\,\rangle_{\mathfrak g}$ (for example, the $(-1)$-multiple of the Killing form of $\mathfrak g$) of 
the Lie algebra $\mathfrak g$ of $G$, where ${\rm Ad}$ denotes the adjoint representation of $G$.  Denote by $g_G$ the bi-invariant metric of $G$ induced from 
$\langle\,\,,\,\,\rangle_{\mathfrak g}$.  Let $\mathcal A_P^{\infty}$ be the affine Hilbert space of all $C^{\infty}$-connections of $P$.  Set 
\begin{align*}
\Omega^{\infty}_{\mathcal T,i}(P,\mathfrak g):=\{A\in\Omega_i^{\infty}(P,\mathfrak g)\,\vert\,\,&R^{\ast}_gA={\rm Ad}(g^{-1})\circ A\,\,\,(\forall\,g\in G)\\
&A(\cdots,\vv,\cdots)=0\,\,\,\,(\forall\,\vv\in\mathcal V)\},
\end{align*}
where $\Omega_i^{\infty}(P,\mathfrak g)$ is the space of all $\mathfrak g$-valued $i$-forms of class $C^{\infty}$ on $P$, 
$R_g$ denotes the right translation by $g$ and $\mathcal V$ denotes the vertical distribution of the bundle $P$.  
Each element of $\Omega^{\infty}_{\mathcal T,i}(P,\mathfrak g)$ is called a {\it $\mathfrak g$-valued tensorial $i$-form of class $C^{\infty}$} on $P$.  
Also, let 
$$\Omega_i^{\infty}(B,{\rm Ad}(P))(=\Gamma^{\infty}((\wedge^iT^{\ast}B)\otimes{\rm Ad}(P)))$$
be the space of all ${\rm Ad}(P)$-valued $i$-forms of class $C^{\infty}$ over $B$, where ${\rm Ad}(P)$ denotes the adjoint bundle 
$P\times_{\rm Ad}\mathfrak g$.  The space ${\mathcal A}_P^{\infty}$ is the affine space having $\Omega^{\infty}_{\mathcal T,1}(P,\mathfrak g)$ as the associated vector space.  
Furthermore, $\Omega^{\infty}_{\mathcal T,i}(P,\mathfrak g)$ is identified with $\Omega_i^{\infty}(B,{\rm Ad}(P))$ by identifying 
$A\in\Omega_{\mathcal T,i}^{\infty}(P,\mathfrak g)$ with $\widehat A\in\Omega_i^{\infty}(B,{\rm Ad}(P))$ defined by 
{\small
$$\left\{\begin{array}{ll}
\displaystyle{\widehat A_{\pi(u)}(\pi_{\ast}(\vv_1),\cdots,\pi_{\ast}(\vv_i))=u\cdot A_u(\vv_1,\cdots,\vv_i)\quad(u\in P,\,\,\vv_1,\cdots,\vv_i\in T_uP)} & (i\geq 1)\\
\displaystyle{\widehat A(\pi(u))=u\cdot A(u)\quad\,\,(u\in P)} & (i=0)
\end{array}\right.$$
}
According to \cite{GP2}, 
We shall define the $H^s$-completion of $\Omega_{\mathcal T,i}^{\infty}(P,\mathfrak g)$, where $s\geq 0$.  
Denote by $\mathcal A_P^{w,s}$ the space of all $s$-times weak differentiable connections of $P$ and 
$\Omega_{\mathcal T,i}^{w,s}(P,\mathfrak g)$ the space of all $s$-times weak differentiable tensorial $i$-forms on $P$.  
Fix a $C^{\infty}$-connection $\omega_0$ of $P$ as the base point of $\mathcal A_P^{w,s}$.  
Set 
$$\Delta_{\omega_0}:=\left\{
\begin{array}{ll}
d_{\omega_0}\circ d_{\omega_0}^{\ast}+d_{\omega_0}^{\ast}\circ d_{\omega_0} & (i\geq 1)\\
d_{\omega_0}^{\ast}\circ d_{\omega_0} & (i=0)
\end{array}\right.$$
(the Hodge-de Rham Laplace operator with respect to $\omega_0$), 
where $d_{\omega_0}$ denotes the covariant exterior derivative with respect to $\omega_0$ and 
$d_{\omega_0}^{\ast}$ denotes the adjoint operator of $d_{\omega_0}$ with respect to the $L^2$-inner products of 
$\Omega_i^{w,j}(B,{\rm Ad}(P))$ ($i\geq 0,\,\,\,j\geq 1$).  
Define an operator $\square_{\omega_0}:\Omega_{\mathcal T,i}^{w,s}(P,\mathfrak g)\to\Omega_{\mathcal T,i}^{w,s-2}(P,\mathfrak g)$ by 
$$
\square_{\omega_0}:=\left\{\begin{array}{ll}
d_{\omega_0}\circ d_{\omega_0}^{\ast}+d_{\omega_0}^{\ast}\circ d_{\omega_0}+{\rm id} & (i\geq 1)\\
d_{\omega_0}^{\ast}\circ d_{\omega_0}+{\rm id} & (i=0).
\end{array}\right.$$
The $L^2_s$-inner product $\langle\,\,,\,\,\rangle^{\omega_0}_s$ of 
$$T_{\omega}\mathcal A_P^{w,s}(\approx\Omega_{\mathcal T,1}^{w,s}(P,\mathfrak g)\approx\Omega_1^{w,s}(B,{\rm Ad}(P))\approx\Gamma^{w,s}(T^{\ast}B\otimes{\rm Ad}(P)))$$
is defined by 
$$\begin{array}{r}
\displaystyle{\langle A_1,A_2\rangle^{\omega_0}_s
:=\int_{x\in B}\langle(\widehat A_1)_x,(\widehat{\square_{\omega_0}^s(A_2)})_x\rangle_{B,\mathfrak g}\,dv_B}\\
\displaystyle{(A_1,A_2\in\Omega_{\mathcal T,1}^{w,s}(P,\mathfrak g)),}
\end{array}\leqno{(6.1)}$$
where $\widehat{\square_{\omega_0}^s(A_2)}$ denotes the ${\rm Ad}(P)$-valued $1$-form over $B$ corresponding to 
$\square_{\omega_0}^s(A_2)$, $\langle\,\,,\,\,\rangle_{B,\mathfrak g}$ denotes the fibre metric of 
$T^{\ast}B\otimes{\rm Ad}(P)$ defined naturally from $g_B$ and $\langle\,\,,\,\,\rangle_{\mathfrak g}$ and $dv_B$ denotes the volume element of $g_B$.  
Let $\Omega_{\mathcal T,1}^{H^s}(P,\mathfrak g)$ be the completion of $\Omega_{\mathcal T,1}^{\infty}(P,\mathfrak g)$ with respect to 
$\langle\,\,,\,\,\rangle^{\omega_0}_s$, that is, 
$$\Omega_{\mathcal T,1}^{H^s}(P,\mathfrak g):=\{A\in\Omega_{\mathcal T,1}^{w,s}(P,\mathfrak g)\,\vert\,\langle A,A\rangle^{\omega_0}_s<\infty\}.$$
Also, set 
$${\mathcal A}_P^{H^s}:=\{\omega_0+A\,\vert\,A\in\Omega_{\mathcal T,1}^{H^s}(P,\mathfrak g)\}.$$
Let $\Omega_1^{H^s}(B,{\rm Ad}(P))$ be the completion of 
$\Omega_1^{\infty}(B,{\rm Ad}(P))$ corresponding to $\Omega_{\mathcal T,1}^{H^s}(P,\mathfrak g)$.  


Let $\mathcal G_P^{\infty}$ be the group of all $C^{\infty}$-gauge transformations ${\bf g}$'s of $P$ with $\pi\circ{\bf g}=\pi$.  
The gauge action $\mathcal G_P^{\infty}\curvearrowright\mathcal A_P^{\infty}$ is given by 
$$\begin{array}{r}
\displaystyle{({\bf g}\cdot\omega)_u={\rm Ad}(\widehat{\bf g}(u))\circ\omega_u-(R_{\widehat{\bf g}(u)})_{\ast}^{-1}\circ\widehat{\bf g}_{\ast u}}\\
\displaystyle{({\bf g}\in\mathcal G_P^{\infty},\,\,\,\omega\in\mathcal A_P^{\infty}).}
\end{array}\leqno{(6.2)}$$
For each ${\bf g}\in\mathcal G_P^{\infty}$, 
$\widehat{\bf g}\in C^{\infty}(P,G)$ is defined by ${\bf g}(u)=u\widehat{\bf g}(u)\,\,\,(u\in P)$.  
This element $\widehat{\bf g}$ satisfies 
$$\widehat{\bf g}(ug)={\rm Ad}(g^{-1})(\widehat{\bf g}(u))\,\,\,\,(\forall u\in P,\,\,\forall g\in G),$$
where ${\rm Ad}$ denotes the homomorphism of $G$ to ${\rm Aut}(G)$ defined by 
${\rm Ad}(g_1)(g_2):=g_1\cdot g_2\cdot g_1^{-1}$ ($g_1,g_2\in G$).  
Under the correspondence ${\bf g}\leftrightarrow\widehat{\bf g}$, $\mathcal G_P^{\infty}$ is identified with 
$$\widehat{\mathcal G}_P^{\infty}:=\{\widehat{\bf g}\in C^{\infty}(P,G)\,|\,
\widehat{\bf g}(ug)={\rm Ad}(g^{-1})(\widehat{\bf g}(u))\,\,\,(\forall\,u\in P,\,\,\forall\,g\in G)\}.$$
For $\widehat{\bf g}\in\widehat{\mathcal G}_P^{\infty}$, 
the $C^{\infty}$-section $\breve{\bf g}$ of the associated $G$-bundle $P\times_{{\rm Ad}}G$ is 
defined by 
$\breve{\bf g}(x):=u\cdot\widehat{\bf g}(u)\,\,\,(x\in B)$, where $u$ is any element of $\pi^{-1}(x)$.  
Under the correspondence $\widehat{\bf g}\leftrightarrow\breve{\bf g}$, 
$\widehat{\mathcal G}_P^{\infty}(=\mathcal G_P^{\infty})$ is identified with the space 
$\Gamma^{\infty}(P\times_{{\rm Ad}}G)$ of all $C^{\infty}$-sections of $P\times_{{\rm Ad}}G$.  
The $H^{s+1}$-completion of $\Gamma^{\infty}(P\times_{{\rm Ad}}G)$ was defined in \cite{GP1} (see Section 1 (P668)).  
Denote by $\Gamma^{H^{s+1}}(P\times_{{\rm Ad}}G)$ this $H^{s+1}$-completion.  
Also, denote by $\mathcal G_P^{H^{s+1}}$ (resp. $\widehat{\mathcal G}_P^{H^{s+1}}$) the $H^{s+1}$-completion of 
$\mathcal G_P^{\infty}$ (resp. $\widehat{\mathcal G}_P^{\infty}$) corresponding to 
$\Gamma^{H^{s+1}}(P\times_{{\rm Ad}}G)$.  
If $\displaystyle{s>\frac{1}{2}\,{\rm dim}\,B-1}$, then the $H^{s+1}$-gauge transformation group $\mathcal G_P^{H^{s+1}}$ of $P$ is a smooth Hilbert Lie group 
and the gauge action $\mathcal G_P^{H^{s+1}}\curvearrowright\mathcal A_P^{H^s}$ is smooth as stated in Introduction.  
However, by this action, $\mathcal G_P^{H^{s+1}}$ does not act isometrically on the Hilbert space 
$(\mathcal A_P^{H^s},\langle\,\,,\,\,\rangle^{\omega_0}_s)$.  
Define a Riemannian metric ${\it g}_s$ on $\mathcal A_P^{H^s}$ by 
$({\it g}_s)_{\omega}:=\langle\,\,\,\,\,\rangle^{\omega}_s$ ($\omega\in\mathcal A_P^{H^s}$), where 
$\langle\,\,,\,\,\rangle_s^{\omega}$ is the $L^2_s$-inner product defined as in $(6.1)$ by using $\omega$ instead of $\omega_0$.  
This Riemannian metric ${\it g}_s$ is non-flat.  
The gauge transformation group $\mathcal G_P^{H^{s+1}}$ acts isometrically on 
the Riemannian Hilbert manifold $(\mathcal A_P^{H^s},{\it g}_s)$ (see Section 7 about this proof), where we note that 
the Hilbert space $(\mathcal A_P^{H^s},\langle\,\,,\,\,\rangle^{\omega_0}_s)(=(\Omega_{\mathcal T,1}^{H^s}(P,\mathfrak g),\langle\,\,,\,\,\rangle_s^{\omega_0}))$ 
is regarded as the tangent space of $(\mathcal A_P^{H^s},{\it g}_s)$ at $\omega_0$.  
Hence we can give the moduli space $\mathcal M_P^{H^s}:=\mathcal A_P^{H^s}/\mathcal G_P^{H^{s+1}}$ the Riemannian orbimetric 
$\overline{\it g}_s$ such that the orbit map 
$\pi_{\mathcal M_P}:(\mathcal A_P^{H^s},{\it g}_s)\to(\mathcal M_P^{H^s},\overline{\it g}_s)$ is a Riemannian orbisubmersion.  

In \cite{K3}, we defined the notion of a holonomy map along a $C^{\infty}$-loop of constant speed in $(B,g_B)$.  
In more general, we shall define the notion of a holonomy element map along a $C^{\infty}$-curve of constant speed in $(B,g_B)$ as follows.  
Fix $x_0\in B$ and $u_0\in\pi^{-1}(x_0)$.  
Take a $C^{\infty}$-curve $c:[0,1]\to B$ of constant speed $a(>0)$ starting from $x_0$.  Let $\sigma$ be the horizontal lift of $c$ starting from $u_0$ 
with respect to $\omega_0$.  

\vspace{0.25truecm}

\noindent
{\bf Definition 6.1.}\ \ We define a map ${\rm hol}_c:{\mathcal A}_P^{H^s}\to G$ by 
$$\mathcal P_c^{\omega}(u_0)=P_c^{\omega_0}(u_0)\cdot{\rm hol}_c(\omega).\leqno{(1.1)}$$
We call this map the {\it holonomy element map along} $c$, where $\mathcal P_c^{\omega}$ (resp. $\mathcal P_c^{\omega_0}$) denotes the parallel translation along $c$ 
with respect to $\omega$ (resp. $\omega_0$).  

\vspace{0.25truecm}

\noindent
{\it Remark 6.1.}\ \ In \cite{K3}, we defined this map ${\rm hol}_c$ for a $C^{\infty}$-loop $c:[0,1]\to B$ and called it the holonomy map along $c$.  

\vspace{0.25truecm}

\vspace{0.25truecm}

\noindent
{\bf Definition 6.2.}\ \ Define a map $\mu_c:\mathcal A_P^{H^s}\to H^s([0,1],\mathfrak g)$ by 
$$(\mu_c(\omega))(t):=-\omega_{\sigma(t)}(\sigma'(t))\,(=-(\sigma^{\ast}\omega)_t\left(\frac{d}{dt}\right))\quad\,\,
(t\in[0,1],\,\,\,\,\omega\in\mathcal A_P^{H^s}).\leqno{(6.3)}$$
We call this map $\mu_c$ the {\it pull-back connection map by $c$}.  

\vspace{0.25truecm}

\noindent
{\it Remark 6.2.}\ \ In \cite{K3}, for a $C^{\infty}$-loop $c:[0,1]\to B$, we defined this map $\mu_c$ by 
$$(\mu_c(\omega))(t):=\omega_{\sigma(t)}(\sigma'(t))\quad\,\,(t\in[0,1],\,\,\,\,\omega\in\mathcal A_P^{H^s}).$$
However, in the right-hand side of this definition, ``$-$'' was lacking.  

\vspace{0.25truecm}

It is clear that $\mu_c$ is a bounded linear operator of $(\mathcal A_P^{H^s},\langle\,\,,\,\,\rangle_s^{\omega})$ onto $H^s([0,1],\mathfrak g)$ for any 
$\omega\in\mathcal A_P^{H^s}$.  Define a map $\phi:H^s([0,1],\mathfrak g)\to G$ by 
$$\phi(u):={\bf g}_u(1)\quad\,\,(u\in H^s([0,1],\mathfrak g)),$$
where ${\bf g}_u$ is the element of $H^{s+1}([0,1],G)$ with ${\bf g}_u(0)=e$ ($e\,:\,$ the identity element of $G$) and 
$(R_{{\bf g}_u(t)})_{\ast}^{-1}({\bf g}_u'(t))=u(t)$ $(t\in[0,1])$.  
Here $R_{{\bf g}_u(t)}$ denotes the right translation by ${\bf g}_u(t)$.  
This map $\phi$ is called the {\it parallel transport map for} $G$.  

According to the proof of Lemma 2.1 in \cite{K3}, we can derive the following fact.  

\vspace{0.25truecm}

\noindent
{\bf Proposition 6.1.} {\sl Among ${\rm hol}_c,\,\phi$ and $\mu_c$, the relation ${\rm hol}_c=\phi\circ\mu_c$ holds.}

\vspace{0.25truecm}

According to the proof of Proposition 3.2 in \cite{K3}, we can derive the following fact.  

\vspace{0.25truecm}

\noindent
{\bf Proposition 6.2.} {\sl For any nonnegative integer $s$, 
$$\mu_c:(\mathcal A_P^{H^s},{\rm g}_s)\to(H^s([0,1],\mathfrak g),\langle\,\,,\,\,\rangle_s^o)$$
is a homothetic submersion of coefficient $a$ onto $H^s([0,1],\mathfrak g)$, where $\langle\,\,,\,\,\rangle_s^o$ is the $L^2_s$-inner product defined by 
$\langle\,\,,\,\,\rangle_{\mathfrak g}$.}

\vspace{0.25truecm}6

According to the proof of Proposition 3.3 in \cite{K3}, the following fact holds.  

\vspace{0.25truecm}

\noindent
{\bf Proposition 6.3.} {\sl For any nonnegative integer $s$, 
$$\phi:(H^s([0,1],\mathfrak g),{\rm g}_s^o)\to(G,{\rm g}_G)$$
is a Riemannian submersion, where $g_G$ is the bi-invarint metric 
defined by $\langle\,\,,\,\,\rangle_{\mathfrak g}$.}

\vspace{0.25truecm}

From Propositions 6.1, 6.2 and 6.3, we can derive the following fact.  

\vspace{0.25truecm}

\noindent
{\bf Theorem 6.4.} 
{\sl The holonomy element map ${\rm hol}_c:(\mathcal A_P^{H^s},g_s)\to(G,g_G)$ is a homothetic submersion of coefficient $a$.}

\vspace{0.25truecm}

Set 
$$\begin{array}{c}
\hspace{1.5truecm}P^{H^{s+1}}(G,e\times G):=\{{\bf g}\in H^{s+1}([0,1],G)\,\vert\,{\bf g}(0)=e\},\\
{\rm and}\quad\,\,\Lambda_e^{H^{s+1}}(G):=\{{\bf g}\in H^{s+1}([0,1],G)\,\vert\,{\bf g}(0)={\bf g}(1)=e\}.\\
\end{array}$$
For $g_0\in G$, define a subgroup $\Delta_{{\rm Ad}(g_0)}(G)$ of $G\times G$ by 
$$\Delta_{{\rm Ad}(g_0)}(G):=\{(g,{\rm Ad}(g_0)(g))\,|\,g\in G\}.$$
Set 
$$P^{H^{s+1}}(G,\Delta_{{\rm Ad}(g_0)}(G)):=\{{\bf g}\in H^{s+1}([0,1],G)\,|\,({\bf g}(0),{\bf g}(1))\in\Delta_{{\rm Ad}(g_0)}(G)\}.$$
Also, set 
$${\rm Ad}(G)_{g_0}:=\{R_{g_0}^{-1}\circ{\rm Ad}(g)\circ R_{g_0}\,|\,g\in G\},$$
where ${\rm Ad}$ is the adjoint action of $G$ on oneself.  
According to the proof of Theorem 6.1 of \cite{TT}, we have 
$$\phi\circ{\bf g}=L_{{\bf g}(0)}\circ R_{{\bf g}(1)}^{-1}\circ\phi\leqno{(6.4)}$$
for any ${\bf g}\in H^{s+1}([0,1],G)$, where ${\bf g}$ in the left-hand side means the diffeomorphism of $H^s([0,1],\mathfrak g)$ onto oneself defined by the action of 
$H^{s+1}([0,1],G)$ on $H^0([0,1],\mathfrak g)$ and $L_{{\bf g}(0)}$ and $R_{{\bf g}(1)}$ are the left translation by ${\bf g}(0)$ and the right translation by ${\bf g}(1)$, 
respectively.  
From $(6.4)$, we can derive the following facts.  

\vspace{0.15truecm}

(i)\ \ The subaction $P^{H^{s+1}}(G,e\times G)\curvearrowright H^s([0,1],\mathfrak g)$ of the gauge action $H^{s+1}([0,1],G)$\newline
$\curvearrowright H^s([0,1],\mathfrak g)$ is transitive (that is, the orbit space $H^s([0,1],\mathfrak g)/P^{H^{s+1}}(G,e\times G)$ is the one-point set.  

(ii)\ \ The orbit space $H^s([0,1],\mathfrak g)/P^{H^{s+1}}(G,\Delta_{{\rm Ad}(g_0)}(G))$ of the subaction \newline
$P^{H^{s+1}}(G,\Delta_{{\rm Ad}(g_0)}(G))\curvearrowright H^s([0,1],\mathfrak g)$ is identified with $G/{\rm Ad}(G)_{g_0}$;

(iii)\ \ The orbit map of the subaction $\Lambda^{H^{s+1}}_e(G)\curvearrowright H^s([0,1],\mathfrak g)$ coincides with the parallel transport map $\phi$ for $G$.  
(hence $H^s([0,1],\mathfrak g)/\Lambda^{H^{s+1}}_e(G)=G$).  

\vspace{0.35truecm}

\noindent
{\bf Definition 6.3.}\ \ (i)\ Define a map $\lambda_c:\mathcal G_P^{H^{s+1}}\to H^{s+1}([0,1],G)$ by 
$$\lambda_c({\bf g})(t):=\widehat{\bf g}\circ\sigma\quad\,\,({\bf g}\in\mathcal G_P^{H^{s+1}}).\leqno{(6.5)}$$

(ii)\ Define $\lambda_G:H^{s+1}([0,1],G)\to G\times G$ by 
$$\lambda_G({\bf g}):=({\bf g}(0),{\bf g}(1))\quad\,\,({\bf g}\in H^{s+1}([0,1],G)).\leqno{(6.6)}$$

\vspace{0.25truecm}

The {\it based gauge transformation group} $(\mathcal G_P^{H^{s+1}})_x$ {\it at} $x\in B$ is defined by 
$$(\mathcal G_P^{H^{s+1}})_x:=\{{\bf g}\in\mathcal G_P^{H^{s+1}}\,\vert\,\widehat{\bf g}(\pi^{-1}(x))=\{e\}\}.$$
Denote by $\mathcal M_P^{H^s}$ the muduli space $\mathcal A_P^{H^s}/\mathcal G_P^{H^{s+1}}$ and 
$\pi_{\mathcal M_P}$ the orbit map of the action $\mathcal G_P^{H^{s+1}}\curvearrowright\mathcal A_P^{H^s}$.  
Also, denote by $(\mathcal M_P^{H^s})_x$ the muduli space $\mathcal A_P^{H^s}/(\mathcal G_P^{H^{s+1}})_x$ and 
$\pi_{(\mathcal M_P)_x}$ the orbit map of the action $(\mathcal G_P^{H^{s+1}})_x\curvearrowright\mathcal A_P^{H^s}$.  
Throughout $\lambda_c$, $\mathcal G_P^{H^{s+1}}$ acts on $H^s([0,1],\mathfrak g)$.  
Also, throughout $\lambda_G\circ\lambda_c$, $\mathcal G_P^{H^{s+1}}$ acts on $G$ by 
$${\bf g}\cdot g:=((\lambda_G\circ\lambda_c)({\bf g}))\cdot g\quad\,\,({\bf g}\in\mathcal G_P^{H^{s+1}},\,\,\,\,g\in G).$$
According to the proof of Lemma 2.2 in\cite{K3}, we can derive the following facts from Proposition 6.1, $(6.2),\,(6.4)$ and the defintions of $\mu_c$ and $\lambda_c$.  

\vspace{0.25truecm}

\noindent
{\bf Lemma 6.5.} {\sl 
{\rm (i)} The pull-back connection map $\mu_c$ is $\mathcal G_P^{H^{s+1}}$-equivariant, that is, 
the following relation holds:
$$\mu_c({\bf g}\cdot\omega)=\lambda_c({\bf g})\cdot\mu_c(\omega)\quad\,\,
({\bf g}\in\mathcal G_P^{H^{s+1}},\,\,\omega\in\mathcal A_P^{H^s}).\leqno{(6.7)}$$

{\rm (ii)} The holonomy element map ${\rm hol}_c$ is $\mathcal G_P^{H^{s+1}}$-equivariant, that is, 
the following relation holds:
$${\rm hol}_c({\bf g}\cdot\omega)=(\lambda_G\circ\lambda_c)({\bf g})\cdot({\rm hol}_c(\omega))\quad\,\,
({\bf g}\in\mathcal G_P^{H^{s+1}},\,\,\omega\in\mathcal A_P^{H^s}).\leqno{(6.8)}$$

{\rm(iii)} Set $x_0:=c(0)$.  Then $\mu_c$ maps each $(\mathcal G_P^{H^{s+1}})_{x_0}$-orbit onto $H^s([0,1],\mathfrak g)$.  
If $c(0)=c(1)$, then $\mu_c$ maps $(\mathcal G_P^{H^{s+1}})_{x_0}$-orbits in $\mathcal A_P^{H^s}$ to $\Lambda_e^{H^{s+1}}(G)$-orbits in $H^s([0,1],\mathfrak g)$ 
and hence there uniquely exists the map $\overline{\mu}_c$ of the based moduli space $(\mathcal M_P^{H^s})_{x_0}$ onto $H^s([0,1],\mathfrak g)/\Lambda_e^{H^{s+1}}(G)(\approx G)$ 
satisfying 
$$\overline{\mu}_c\circ\pi_{(\mathcal M_P)_{x_0}}=\phi\circ\mu_c(={\rm hol}_c)\leqno{(6.9)}$$
(see Diagram 6.1).  

{\rm(iv)} Assume that $c(0)=c(1)$ and let $\sigma(1)=u_0\cdot g_0$ ($g_0\in G$).  Then the pull-back connection map $\mu_c$ maps $\mathcal G_P^{H^{s+1}}$-orbits 
in $\mathcal A_P^{H^s}$ to $P^{H^{s+1}}(G,\Delta_{{\rm Ad}(g_0^{-1})}(G))$-orbits in $H^s([0,1],\mathfrak g)$ 
and hence there uniquely exists the map $\overline{\overline{\mu}}_c$ of the moduli sapce $\mathcal M_P^{H^s}$ onto 
$H^s([0,1],\mathfrak g)/P^{H^{s+1}}(G,\Delta_{{\rm Ad}(g_0^{-1})}(G))(\approx G/{\rm Ad}(G)_{g_0^{-1}})$ satisfying 
$$\overline{\overline{\mu}}_c\circ\pi_{\mathcal M_P}=\pi_{{\rm Ad}_{g_0^{-1}}}\circ\phi\circ\mu_c(=\pi_{{\rm Ad}_{g_0^{-1}}}\circ{\rm hol}_c=:\overline{\rm hol}_c)\leqno{(6.10)}$$
(see Diagram 6.2), where $\pi_{{\rm Ad}_{g_0^{-1}}}$ denotes the natural projection of $G$ onto $G/{\rm Ad}(G)_{g_0^{-1}}$.  
}

\vspace{0.25truecm}


{\small 
\centerline{
\unitlength 0.1in
\begin{picture}( 52.7000, 15.0500)(  5.4000,-24.3500)
\put(40.7000,-16.6000){\makebox(0,0)[rt]{$\mathcal A_P^{H^s}$}}%
\put(53.5000,-16.6000){\makebox(0,0)[lt]{$H^s([0,1],\mathfrak g)$}}%
\put(34.2000,-12.2000){\makebox(0,0)[rb]{$(\mathcal G_P^{H^{s+1}})_{x_0}$}}%
\put(52.2000,-12.2000){\makebox(0,0)[rb]{$\Lambda_e^{H^{s+1}}(G)$}}%
%
\special{pn 8}%
\special{pa 3620 1160}%
\special{pa 4510 1160}%
\special{fp}%
\special{sh 1}%
\special{pa 4510 1160}%
\special{pa 4444 1140}%
\special{pa 4458 1160}%
\special{pa 4444 1180}%
\special{pa 4510 1160}%
\special{fp}%
\put(39.6000,-11.0000){\makebox(0,0)[lb]{$\lambda_c$}}%
%
\special{pn 8}%
\special{pa 4230 1740}%
\special{pa 5180 1740}%
\special{fp}%
\special{sh 1}%
\special{pa 5180 1740}%
\special{pa 5114 1720}%
\special{pa 5128 1740}%
\special{pa 5114 1760}%
\special{pa 5180 1740}%
\special{fp}%
\put(45.5000,-16.9000){\makebox(0,0)[lb]{$\mu_c$}}%
%
\special{pn 13}%
\special{ar 3280 1620 590 450  5.0767644 5.9136541}%
%
\special{pn 13}%
\special{pa 3840 1460}%
\special{pa 3880 1540}%
\special{fp}%
\special{sh 1}%
\special{pa 3880 1540}%
\special{pa 3868 1472}%
\special{pa 3856 1492}%
\special{pa 3832 1490}%
\special{pa 3880 1540}%
\special{fp}%
%
\special{pn 13}%
\special{ar 5080 1630 590 450  5.0767644 5.9136541}%
%
\special{pn 13}%
\special{pa 5640 1470}%
\special{pa 5680 1550}%
\special{fp}%
\special{sh 1}%
\special{pa 5680 1550}%
\special{pa 5668 1482}%
\special{pa 5656 1502}%
\special{pa 5632 1500}%
\special{pa 5680 1550}%
\special{fp}%
%
\special{pn 8}%
\special{pa 3890 1920}%
\special{pa 3890 2290}%
\special{fp}%
\special{sh 1}%
\special{pa 3890 2290}%
\special{pa 3910 2224}%
\special{pa 3890 2238}%
\special{pa 3870 2224}%
\special{pa 3890 2290}%
\special{fp}%
%
\special{pn 8}%
\special{pa 5730 1890}%
\special{pa 5730 2260}%
\special{fp}%
\special{sh 1}%
\special{pa 5730 2260}%
\special{pa 5750 2194}%
\special{pa 5730 2208}%
\special{pa 5710 2194}%
\special{pa 5730 2260}%
\special{fp}%
\put(41.4000,-23.4000){\makebox(0,0)[rt]{$(\mathcal M_P^{H^s})_{x_0}$}}%
\put(56.7000,-23.5000){\makebox(0,0)[lt]{$G$}}%
%
\special{pn 8}%
\special{pa 4330 2430}%
\special{pa 5480 2430}%
\special{fp}%
\special{sh 1}%
\special{pa 5480 2430}%
\special{pa 5414 2410}%
\special{pa 5428 2430}%
\special{pa 5414 2450}%
\special{pa 5480 2430}%
\special{fp}%
\put(47.6000,-24.0000){\makebox(0,0)[lb]{$\overline{\mu}_c$}}%
%
\special{pn 8}%
\special{pa 4200 1870}%
\special{pa 5450 2280}%
\special{fp}%
\special{sh 1}%
\special{pa 5450 2280}%
\special{pa 5394 2240}%
\special{pa 5400 2264}%
\special{pa 5380 2278}%
\special{pa 5450 2280}%
\special{fp}%
\put(47.9000,-20.4000){\makebox(0,0)[lb]{{\small ${\rm hol}_c$}}}%
\put(58.1000,-20.1000){\makebox(0,0)[lt]{$\phi$}}%
%
\special{pn 8}%
\special{ar 4312 2142 72 72  0.6528466 5.7088805}%
%
\special{pn 8}%
\special{pa 4372 2102}%
\special{pa 4392 2152}%
\special{fp}%
\special{sh 1}%
\special{pa 4392 2152}%
\special{pa 4386 2084}%
\special{pa 4372 2102}%
\special{pa 4350 2098}%
\special{pa 4392 2152}%
\special{fp}%
\put(37.9000,-21.5000){\makebox(0,0)[rb]{$\pi_{(\mathcal M_P)_{x_0}}$}}%
%
\special{pn 8}%
\special{ar 5412 2002 72 72  0.6528466 5.7088805}%
%
\special{pn 8}%
\special{pa 5472 1962}%
\special{pa 5492 2012}%
\special{fp}%
\special{sh 1}%
\special{pa 5492 2012}%
\special{pa 5486 1944}%
\special{pa 5472 1962}%
\special{pa 5450 1958}%
\special{pa 5492 2012}%
\special{fp}%
\end{picture}%
\hspace{8truecm}}
}

\vspace{0.25truecm}

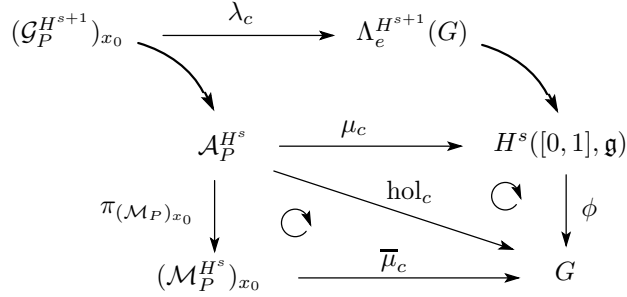
\centerline{{\bf Diagram 6.1$\,\,:\,\,$ The map $\overline{\mu}_c$ induced from $\mu_c$ ($c(0)=c(1)$-case)}}

\vspace{0.35truecm}

{\small 
\centerline{
\unitlength 0.1in
\begin{picture}( 45.5000, 15.0500)(  3.8000,-23.1500)
\put(31.9000,-15.4000){\makebox(0,0)[rt]{$\mathcal A_P^{H^s}$}}%
\put(44.7000,-15.4000){\makebox(0,0)[lt]{$H^s([0,1],\mathfrak g)$}}%
\put(25.4000,-11.0000){\makebox(0,0)[rb]{$\mathcal G_P^{H^{s+1}}$}}%
\put(38.0000,-10.9000){\makebox(0,0)[lb]{$P^{H^{s+1}}(G,\Delta_{{\rm Ad}(g_0)}(G))$}}%
%
\special{pn 8}%
\special{pa 2740 1040}%
\special{pa 3630 1040}%
\special{fp}%
\special{sh 1}%
\special{pa 3630 1040}%
\special{pa 3564 1020}%
\special{pa 3578 1040}%
\special{pa 3564 1060}%
\special{pa 3630 1040}%
\special{fp}%
\put(30.8000,-9.8000){\makebox(0,0)[lb]{$\lambda_c$}}%
%
\special{pn 8}%
\special{pa 3350 1620}%
\special{pa 4300 1620}%
\special{fp}%
\special{sh 1}%
\special{pa 4300 1620}%
\special{pa 4234 1600}%
\special{pa 4248 1620}%
\special{pa 4234 1640}%
\special{pa 4300 1620}%
\special{fp}%
\put(36.7000,-15.7000){\makebox(0,0)[lb]{$\mu_c$}}%
%
\special{pn 13}%
\special{ar 2400 1500 590 450  5.0767644 5.9136541}%
%
\special{pn 13}%
\special{pa 2960 1340}%
\special{pa 3000 1420}%
\special{fp}%
\special{sh 1}%
\special{pa 3000 1420}%
\special{pa 2988 1352}%
\special{pa 2976 1372}%
\special{pa 2952 1370}%
\special{pa 3000 1420}%
\special{fp}%
%
\special{pn 13}%
\special{ar 4302 1496 560 370  5.0752732 5.9119355}%
%
\special{pn 13}%
\special{pa 4832 1364}%
\special{pa 4870 1430}%
\special{fp}%
\special{sh 1}%
\special{pa 4870 1430}%
\special{pa 4854 1362}%
\special{pa 4844 1384}%
\special{pa 4820 1382}%
\special{pa 4870 1430}%
\special{fp}%
%
\special{pn 8}%
\special{pa 3010 1800}%
\special{pa 3010 2170}%
\special{fp}%
\special{sh 1}%
\special{pa 3010 2170}%
\special{pa 3030 2104}%
\special{pa 3010 2118}%
\special{pa 2990 2104}%
\special{pa 3010 2170}%
\special{fp}%
%
\special{pn 8}%
\special{pa 4850 1770}%
\special{pa 4850 2140}%
\special{fp}%
\special{sh 1}%
\special{pa 4850 2140}%
\special{pa 4870 2074}%
\special{pa 4850 2088}%
\special{pa 4830 2074}%
\special{pa 4850 2140}%
\special{fp}%
\put(32.6000,-22.2000){\makebox(0,0)[rt]{$\mathcal M_P^{H^s}$}}%
\put(47.9000,-22.3000){\makebox(0,0)[lt]{$G/{\rm Ad}(G)_{g_0^{-1}}$}}%
%
\special{pn 8}%
\special{pa 3450 2310}%
\special{pa 4600 2310}%
\special{fp}%
\special{sh 1}%
\special{pa 4600 2310}%
\special{pa 4534 2290}%
\special{pa 4548 2310}%
\special{pa 4534 2330}%
\special{pa 4600 2310}%
\special{fp}%
\put(38.8000,-22.8000){\makebox(0,0)[lb]{$\overline{\overline{\mu}}_c$}}%
%
\special{pn 8}%
\special{pa 3320 1750}%
\special{pa 4570 2160}%
\special{fp}%
\special{sh 1}%
\special{pa 4570 2160}%
\special{pa 4514 2120}%
\special{pa 4520 2144}%
\special{pa 4500 2158}%
\special{pa 4570 2160}%
\special{fp}%
\put(39.1000,-19.2000){\makebox(0,0)[lb]{{\small $\overline{{\rm hol}}_c$}}}%
\put(49.3000,-18.9000){\makebox(0,0)[lt]{$\pi_{\rm Ad}\circ\phi$}}%
%
\special{pn 8}%
\special{ar 3432 2022 72 72  0.6528466 5.7088805}%
%
\special{pn 8}%
\special{pa 3492 1982}%
\special{pa 3512 2032}%
\special{fp}%
\special{sh 1}%
\special{pa 3512 2032}%
\special{pa 3506 1964}%
\special{pa 3492 1982}%
\special{pa 3470 1978}%
\special{pa 3512 2032}%
\special{fp}%
\put(29.1000,-20.3000){\makebox(0,0)[rb]{$\pi_{\mathcal M_P}$}}%
%
\special{pn 8}%
\special{ar 4532 1882 72 72  0.6528466 5.7088805}%
%
\special{pn 8}%
\special{pa 4592 1842}%
\special{pa 4612 1892}%
\special{fp}%
\special{sh 1}%
\special{pa 4612 1892}%
\special{pa 4606 1824}%
\special{pa 4592 1842}%
\special{pa 4570 1838}%
\special{pa 4612 1892}%
\special{fp}%
\end{picture}%
\hspace{6truecm}}
}

\vspace{0.25truecm}

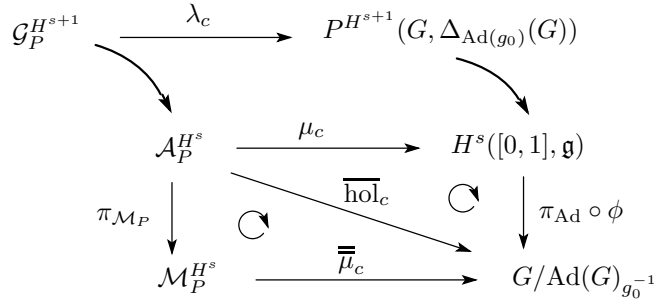
\centerline{{\bf Diagram 6.2$\,\,:\,\,$ The map $\overline{\overline{\mu}}_c$ induced from $\mu_c$ ($c(0)=c(1)$-case)}}

\section{The Green function} 
In this section, we recall the definition of Green function of an invertible linear differential operator from the space of sections of a vector bundle to oneself.  
Let $\pi:E\to B$ be a $C^{\infty}$-vector bundle (equipped with a fibre metric $g_E$) over a compact Riemannian manifold $B$ and 
$\mathcal L$ an invertible linear differential operator from the space $\Gamma^{H^s}(E)$ to oneself, where 
$\Gamma^{H^s}(E)$ denotes the $H^s$-completion of the sapce $\Gamma^{\infty}(E)$ of $C^{\infty}$-sections of $E$.  
We consider the linear differential equation 
$$\mathcal L(\sigma)=\rho,\leqno{(7.1)}$$
where $\rho$ is a $H^s$-section of $E$.  
Take $y\in B$ and $\vv\in E_y$.  Define the section $\delta_{\vv}$ of $E$ by 
$$\langle\delta_{\vv},\sigma\rangle_0=(g_E)_y(\vv,\sigma(y))\quad\,\,(\forall\,\sigma\in\Gamma(E))$$
Let $G_{\mathcal L}^{\vv}\in\Gamma(E)$ ($\vv\in E_y$) be the solution of $\mathcal L(\sigma)=\delta_{\vv}$.  
Let $pr_i:B\times B\to B$ ($i=1,2$) be the projection of the $i$-th component.  
We consider the vector bundle ${\rm pr_1}^{\ast}E\otimes{\rm pr}_2^{\ast}E^{\ast}$ over $B\times B$ 
(i.e., $\displaystyle{{\rm pr_1}^{\ast}E\otimes{\rm pr}_2^{\ast}E^{\ast}:=\mathop{\amalg}_{(x,y)\in B\times B}(E_x\otimes E_y^{\ast})}$).  
Let $G_{\mathcal L}$ be the section of ${\rm pr_1}^{\ast}E\otimes{\rm pr}_2^{\ast}E^{\ast}$ defined by 
$$G_{\mathcal L}(x,y)(\vv):=G_{\mathcal L}^{\vv}(x)\quad\,\,((x,y)\in B\times B,\,\,\,\,\vv\in E_y).$$
Here we note that $G_{\mathcal L}^{\vv}$ exists uniquely because $\mathcal L$ is invertible.  
This section $G_{\mathcal L}$ is called the {\it Green function} of $\mathcal L$.  
Here we note that 
$$G_{\mathcal L}(x,y)(\vv)=\sum_i\lambda_i^{-1}\cdot(g_E)_y(\vv,\eta_i(y))\,\delta_{y,1}(x)\eta_i(x),$$
where $(\eta_i)_i$ is a $H^s$-orthonormal basis of $\Gamma^{H^s}(E)$ with $\mathcal L(\eta_i)=\lambda_i\eta_i$ ($\lambda_i\in\mathbb R$) and $\delta_{y,1}$ denotes 
the delta function defined by $\int_B\delta_{y,1}(x)f(x)\,dx=f(y)$ ($\forall\,f\in C^{\infty}(B)$).  
In fact, by operating $\mathcal L$ to the right-hand side of the above relation and using the linearity of $\mathcal L$, it is equal to $\delta_{\vv}$.  
The {\it Green operator} $G^o_{\mathcal L}:\Gamma^{H^s}(E)\to\Gamma^{H^s}(E)$ of $\mathcal L$ is defined by 
$$G^o_{\mathcal L}(\psi)(y):=\sum_{i=1}^m\int_{x\in B}g_E(\psi(x),G_{\mathcal L}(x,y)(\ve_i))\,dv_B\cdot\ve_i\quad\,\,(\psi\in\Gamma^{H^s}(E),\,\,\,y\in B),$$
where $(\ve_1,\cdots,\ve_m)$ is an orthonormal basis of $E_y$.  
Then we have 
\begin{align*}
\mathcal L(G^o_{\mathcal L}(\psi))(y)&=\sum_{i=1}^m\int_{x\in B}g_E(\psi(x),\mathcal L(G_{\mathcal L}(x,y)(\ve_i))\,dv_B\cdot\ve_i\\
&=\sum_{i=1}^m\int_{x\in B}g_E(\psi(x),\delta_{\ve_i}(x))\,dv_B\cdot\ve_i\\
&=\sum_{i=1}^mg_E(\psi(y),\ve_i)\cdot\ve_i=\psi(y).
\end{align*}
Thus $G^o_{\mathcal L}(\psi)$ is a solution of $\mathcal L(\sigma)=\psi$.  

\section{The orbit geometry of the action $\mathcal G_P^{H^{s+1}}\curvearrowright(\mathcal A_P^{H^s},g_s)$} 
In this section, we shall recall known facts for the orbit geometry of the action $\mathcal G_P^{H^{s+1}}\curvearrowright(\mathcal A_P^{H^s},g_s)$.  
We shall use the notations in Section 6.  
Take $\omega\in\mathcal A_P^{H^s}$ and set $A:=\omega-\omega_0\in\Omega_{\mathcal T,1}^{H^s}(P,\mathfrak g)=\Omega_1^{H^s}(B,{\rm Ad}(P))$.  
Denote by $\mathcal V^P$ the vertical distribution of $\pi:P\to B$ and $\mathcal H^{\omega}$ the horizontal distribution with respect to $\omega$ 
(i.e., $(\mathcal H^{\omega})_u={\rm Ker}\,\omega_u$ ($u\in P$)).  
Set $\Delta_{\omega}^0:=d_{\omega}^{\ast}\circ d_{\omega}$.  

\vspace{0.25truecm}

\noindent
{\bf Proposition 8.1.} {\sl The connection $\omega$ is irreducible if and only if $\Delta_{\omega}^0$ is invertible.}  

\vspace{0.25truecm}

\noindent
{\it Proof.}\ Assume that $\omega$ is irreducible.  Then the holonomy group of $\omega$ coincides with the whole of $G$.  
Take $f\in{\rm Ker}\,d_{\omega}$.  Then we have $df|_{\mathcal H^{\omega}}=0$.  
That is, $f$ is constant along any $\omega$-horizontal curve.  Hence $f$ is constant on each holonomy bundle with respect to the connection $\omega$.  
Since the holonomy group of $\omega$ coincides with the whole of $G$, the holonomy bundle coincides with the whole of $P$.  
Therefore $f$ is constant on the whole of $P$.  
On the other hand, since $f\in\Omega_{\mathcal T,0}(P,\mathfrak g)$, we have 
$$f(ua)={\rm Ad}(a^{-1})(f(u))\quad(u\in P,\,\,\,a\in G),$$
which together with the constancy of $f$ implies 
$$f(u)={\rm Ad}(a^{-1})(f(u))\quad(u\in P,\,\,\,a\in G),$$
that is, $f(u)\in{\rm Fix}({\rm Ad}(G))=\{{\bf 0}\}$.  Thus we obtain $f\equiv{\bf 0}$.  
Therefore we can conclude ${\rm Ker}\,d_{\omega}=\{{\bf 0}\}$, that is, $d_{\omega}$ is invertible.  
Hence so is $\Delta_{\omega}^0$ also.  The converse also is shown by the similar discussion.  \qed

\vspace{0.25truecm}

In the case where $\omega$ is irreducible, the Green function $G_{\Delta_{\omega}^0}$ of $\Delta_{\omega}^0$ exists uniquely as stated in Section 7 because 
$\Delta_{\omega}^0$ is invertible.  
Note that, in the case where $\omega$ is reducible, for each $\vv\in{\rm Ad}(P)_x$, we can find $G_{\Delta_{\omega}^0}^{\vv}\in\Gamma^{H^s}({\rm Ad}(P))$ satisfying 
$$\Delta_{\omega}^0(G_{\Delta_{\omega}^0}^{\vv})=\delta_{\vv}-\sum_ig(\vv,\eta_i(x))\eta_i,$$
where $(\eta_i)_i$ is an orthonormal basis of ${\rm Ker}\,\Delta_A^0$.  
The quotient group $\mathcal G_P^{H^{s+1}}/(\mathcal G_P^{H^{s+1}})_{x_0}$ is isomorphic to $G/Z$, where $Z$ is the center of $G$.  
For example, if $G=SU(2)$, then we have $Z=S^1$ and hence $G/Z=SO(3)$.  
The group $\mathcal G_P^{H^{s+1}}/(\mathcal G_P^{H^{s+1}})_{x_0}(\cong G/Z)$ acts on $(\mathcal M_P^{H^s})_{x_0}$ and the orbit space 
$(\mathcal M_P^{H^s})_{x_0}/(G/Z)$ is identified with $\mathcal M_P^{H^s}$.  
It is shown that $\pi_{(\mathcal M_P)_{x_0}}$ is a submersion and that $\pi_{\mathcal M_P}$ is an orbisubmersion.  

For the action $\mathcal G_P^{H^{s+1}}\curvearrowright(\mathcal A_P^{H^s},g_s)$, we have the following fact.  

\vspace{0.35truecm}

\noindent
{\bf Proposition 8.2.} {\sl The action $\mathcal G_P^{H^{s+1}}\curvearrowright(\mathcal A_P^{H^s},g_s)$ is isometric.}

\vspace{0.35truecm}

\noindent
{\it Proof.}\ \ 
Take ${\bf g}\in\mathcal G_P^{H^{s+1}}$ and $A_i\in T_{\omega}\mathcal A_P^{H^s}(\approx\Omega_{\mathcal T,1}^{H^s}(P,\mathfrak g))$ ($i=1,2$).  
Easily we can show 
$${\bf g}_{\ast\omega}(A_i)={\rm Ad}(\widehat{\bf g})(A_i)$$
and 
$$\square_{{\bf g}\cdot\omega}\circ{\rm Ad}(\widehat{\bf g})={\rm Ad}(\widehat{\bf g})\circ\square_{\omega},$$
where ${\rm Ad}(\widehat{\bf g})$ is a linear map defined by 
$$({\rm Ad}(\widehat{\bf g})(A))(u)={\rm Ad}(\widehat{\bf g}(u))(A(u))\quad\,\,(u\in P).$$
From these relations, we can derive 
\begin{align*}
(g_s)_{{\bf g}\cdot\omega}({\bf g}_{\ast\omega}(A_1),{\bf g}_{\ast\omega}(A_2))=&(g_s)_{{\bf g}\cdot\omega}({\rm Ad}(\widehat{\bf g})(A_1),{\rm Ad}(\widehat{\bf g})(A_2))\\
=&\langle{\rm Ad}(\widehat{\bf g})(A_1),{\rm Ad}(\widehat{\bf g})(A_2))\rangle_s^{{\bf g}\cdot\omega}\\
=&\langle{\rm Ad}(\widehat{\bf g})(A_1),\square_{{\bf g}\cdot\omega}^s({\rm Ad}(\widehat{\bf g})(A_2))\rangle_0\\
=&\langle{\rm Ad}(\widehat{\bf g})(A_1),{\rm Ad}(\widehat{\bf g})(\square_{\omega}^s(A_2))\rangle_0\\
=&\int_{x\in B}\langle(\widehat{{\rm Ad}(\widehat{\bf g})(A_1))})_x,(\widehat{{\rm Ad}(\widehat{\bf g})(\square_{\omega}^s(A_2))})_x\rangle_{B,\mathfrak g}\,dv_B\\
=&\int_{x\in B}\langle(\widehat A_1)_x,(\widehat{\square_{\omega}^s(A_2)})_x\rangle_{B,\mathfrak g}\,dv_B\\
=&\langle A_1,\square_{\omega}^s(A_2)\rangle_0=\langle A_1,A_2\rangle_s^{\omega}=(g_s)_{\omega}(A_1,A_2).
\end{align*}
Thus ${\bf g}:(\mathcal A_P^{H^s},g_s)\to(\mathcal A_P^{H^s},g_s)$ is an isometry.  Therefore $\mathcal G_P^{H^{s+1}}\curvearrowright(\mathcal A_P^{H^s},g_s)$ is isometric.  \qed

\vspace{0.35truecm}

Denote by $\nabla^s$ the Levi-Civita connection of $g_s$.  
Since $\mathcal G_P^{H^{s+1}}$ acts on $(\mathcal A_P^{H^s},g_s)$ isometrically, there exists the Riemannian metric $\overline g_s$ on 
$\mathcal M_P^{H^s}$ such that the orbit map $\pi_{\mathcal M_P}:(\mathcal A_P^{H^s},g_s)\to(\mathcal M_P^{H^s},\overline g_s)$ is a Riemannian orbisubmersion.  
Similarly, there exists the Riemannian metric $\overline g_s^{x_0}$ on 
$(\mathcal M_P^{H^s})_{x_0}$ such that the orbit map $\pi_{(\mathcal M_P)_{x_0}}:(\mathcal A_P^{H^s},g_s)\to((\mathcal M_P^{H^s})_{x_0},\overline g_s^{x_0})$ is 
a Riemannian submersion.  
Denote by $\overline{\nabla}^s$ the Levi-Civita connection of $\overline g_s$.  
Let $\mathcal V^{\mathcal A}$ be the vertical distribution of the Riemannian orbisubmersion $\pi_{\mathcal M_P}$ and 
$\mathcal H^{\mathcal A}$ the horizontal distribution of the Riemannian orbisubmersion $\pi_{\mathcal M_P}$.  
Also, let $\widetilde{\mathcal H}^{\mathcal A}$ the horizontal distribution of the Riemannian submersion $\pi_{(\mathcal M_P)_{x_0}}$.  
The vertical distribution $\mathcal V^{\mathcal A}$ of $\pi_{\mathcal M_P}$ is given by 
$$(\mathcal V^{\mathcal A})_{\omega}={\rm Im}\,d_{\omega}\quad\,\,(\omega\in\mathcal A_P^{H^s})\leqno{(8.1)}$$
and 
the horizontal distribution $\mathcal H^{\mathcal A}$ of $\pi_{\mathcal M_P}$ is given by 
%
%
$$\begin{array}{l}
\displaystyle{(\mathcal H^{\mathcal A})_{\omega}=\{\eta\in\Omega_{\mathcal T,1}^{H^s}(P,\mathfrak g)\,|\,\langle d_{\omega}(\rho),\eta\rangle_s^{\omega}=0\,\,\,\,
(\forall\,\rho\in\Omega_{\mathcal T,0}^{H^s}(P,\mathfrak g))\}}\\
\hspace{1.2truecm}\displaystyle{=\{\eta\in\Omega_{\mathcal T,1}^{H^s}(P,\mathfrak g)\,|\,\langle d_{\omega}(\rho),\square_{\omega}^s(\eta)\rangle_0=0\,\,\,\,
(\forall\,\rho\in\Omega_{\mathcal T,0}^{H^s}(P,\mathfrak g))\}}\\
\hspace{1.2truecm}\displaystyle{=\{\eta\in\Omega_{\mathcal T,1}^{H^s}(P,\mathfrak g)\,|\,\langle\rho,d_{\omega}^{\ast}(\square_{\omega}^s(\eta))\rangle_0=0\,\,\,\,
(\forall\,\rho\in\Omega_{\mathcal T,0}^{H^s}(P,\mathfrak g))\}}\\
\hspace{1.2truecm}\displaystyle{={\rm Ker}\,(d_{\omega}^{\ast}\circ\square_{\omega}^s).}
\end{array}\leqno{(8.2)}$$
Denote by $\widetilde{\mathcal V}^{\mathcal A}$ the vertical distribution of $\pi_{(\mathcal M_P)_{x_0}}$ and 
$\widetilde{\mathcal H}^{\mathcal A}$ the horizontal distribution of $\pi_{(\mathcal M_P)_{x_0}}$.  
The horizontal space $(\widetilde{\mathcal H}^{\mathcal A})_{\omega}$ is given by 
$$\begin{array}{l}
\displaystyle{(\widetilde{\mathcal H}^{\mathcal A})_{\omega}=\{\eta\in\Omega_{\mathcal T,1}^{H^s}(P,\mathfrak g)\,|\,\langle d_{\omega}(\rho),\eta\rangle_s^{\omega}=0}\\
\hspace{4.45truecm}\displaystyle{(\forall\,\rho\in\Omega_{\mathcal T,0}^{H^s}(P,\mathfrak g)\,\,{\rm s.t.}\,\,\rho|_{\pi^{-1}(x_0)}=0)\}}\\
\hspace{1.2truecm}\displaystyle{=\{\eta\in\Omega_{\mathcal T,1}^{H^s}(P,\mathfrak g)\,|\,\langle d_{\omega}(\rho),\square_{\omega}^s(\eta)\rangle_0=0}\\
\hspace{4.45truecm}\displaystyle{(\forall\,\rho\in\Omega_{\mathcal T,0}^{H^s}(P,\mathfrak g)\,\,{\rm s.t.}\,\,\rho|_{\pi^{-1}(x_0)}=0)\}}\\
\hspace{1.2truecm}\displaystyle{=\{\eta\in\Omega_{\mathcal T,1}^{H^s}(P,\mathfrak g)\,|\,\langle\rho,d_{\omega}^{\ast}(\square_{\omega}^s(\eta))\rangle_0=0}\\
\hspace{4.45truecm}\displaystyle{(\forall\,\rho\in\Omega_{\mathcal T,0}^{H^s}(P,\mathfrak g)\,\,{\rm s.t.}\,\,\rho|_{\pi^{-1}(x_0)}=0)\}}\\
\hspace{1.2truecm}\displaystyle{=\{\eta\in\Omega_{\mathcal T,1}^{H^s}(P,\mathfrak g)\,|\,\widehat{(d_{\omega}^{\ast}\circ\square_{\omega}^s)(\eta)}=\delta_{\vv}\,\,
{\rm for}\,\,{\rm some}\,\,\vv\in{\rm Ad}(P)_{x_0}\}}\\
\hspace{1.2truecm}\displaystyle{=\{\eta\in\Omega_{\mathcal T,1}^{H^s}(P,\mathfrak g)\,|\,(d_{\omega}^{\ast}\circ\square_{\omega}^s)(\eta)
=(d_{\omega}^{\ast}\circ d_{\omega})(G_{\Delta_{\omega}^0}^{\vv})}\\
\hspace{4.45truecm}\displaystyle{{\rm for}\,\,{\rm some}\,\,\vv\in{\rm Ad}(P)_{x_0}\}}\\
\hspace{1.2truecm}\displaystyle{=\{\eta\in\Omega_{\mathcal T,1}^{H^s}(P,\mathfrak g)\,|\,(d_{\omega}^{\ast}\circ\square_{\omega}^s)(\eta)
=(d_{\omega}^{\ast}\circ\square_{\omega}^s\circ\square_{\omega}^{-s}\circ d_{\omega})(G_{\Delta_{\omega}^0}^{\vv})}\\
\hspace{4.45truecm}\displaystyle{{\rm for}\,\,{\rm some}\,\,\vv\in{\rm Ad}(P)_{x_0}\}}\\
\hspace{1.2truecm}\displaystyle{=\{\eta\in\Omega_{\mathcal T,1}^{H^s}(P,\mathfrak g)\,|\,(d_{\omega}^{\ast}\circ\square_{\omega}^s)
(\,\eta-(\square_{\omega}^{-s}\circ d_{\omega})(G_{\Delta_{\omega}^0}^{\vv})\,)=0}\\
\hspace{4.45truecm}\displaystyle{{\rm for}\,\,{\rm some}\,\,\vv\in{\rm Ad}(P)_{x_0}\}}\\
\hspace{1.2truecm}\displaystyle{={\rm Ker}\,(d_{\omega}^{\ast}\circ\square_{\omega}^s)}\\
\hspace{1.7truecm}\displaystyle{\oplus{\rm Span}\{\,(\square_{\omega}^{-s}\circ d_{\omega})(G_{\Delta_{\omega}^0}^{\vv})\,|\,\vv\in{\rm Ad}(P)_{x_0}
\}}\\
\hspace{1.2truecm}\displaystyle{=(\mathcal H^{\mathcal A})_{\omega}
\oplus{\rm Span}\{\,(\square_{\omega}^{-s}\circ d_{\omega})(G_{\Delta_{\omega}^0}^{\vv})\,|\,\vv\in{\rm Ad}(P)_{x_0}
\},}
\end{array}\leqno{(8.3)}$$
where $G_{\Delta_{\omega}^0}^{\vv}$ is regarded as the element of $\Omega_{\mathcal T,0}^{H^s}(P,\mathfrak g)$.  
From $(8.1),\,(8.2)$ and $(8.3)$, we obtain 
$$T_{\omega}(\mathcal G_P^{H^{s+1}}\cdot\omega)=(\mathcal V^{\mathcal A})_{\omega}={\rm Im}\,d_{\omega}\,(\cong T_{{\rm id}}(\mathcal G_P^{H^{s+1}})/{\rm Ker}\,d_{\omega})\leqno{(8.4)}$$
and 
$$\begin{array}{l}
\hspace{0.5truecm}\displaystyle{T_{\omega}((\mathcal G_P^{H^{s+1}})_{x_0}\cdot\omega)}\\
\displaystyle{={\rm Im}\,d_{\omega}\cap\left({\rm Span}\{\,(\square_{\omega}^{-s}\circ d_{\omega})(G_{\Delta_{\omega}^0}^{\vv})\,|\,
\vv\in{\rm Ad}(P)_{x_0}
\,\}\right)^{\perp},}
\end{array}\leqno{(8.5)}$$
where we note that $T_{\rm id}\mathcal G_P^{H^{s+1}}$ and $T_{\omega}(\mathcal G_P^{H^{s+1}}\cdot\omega)$ are identified with 
$\Omega_{\mathcal T,0}^{H^{s+1}}(P,\mathfrak g)(=\Gamma^{H^{s+1}}({\rm Ad}(P)))$ 
and a subspace of $\Omega_{\mathcal T,1}^{H^s}(P,\mathfrak g)(=\Omega_1^{H^s}(B,{\rm Ad}(P)))$, 
respectively, and hence $d_{\omega}$ is regarded as a map from $T_{\rm id}\mathcal G_P^{H^{s+1}}$ to 
$T_{\omega}(\mathcal G_P^{H^{s+1}}\cdot\omega)=\mathcal V_{\omega}$.  
From Proposition 8.1 and $(8.4)$, we obtain the following fact.  

\vspace{0.25truecm}

\noindent
{\bf Proposition 8.3.} {\sl (i)\ \ If $\omega$ is an irreducible connection, then the gauge orbit $\mathcal G_P^{H^{s+1}}\cdot\omega$ is a regular orbit of the action 
$\mathcal G_P^{H^{s+1}}\curvearrowright\mathcal A_P^{H^s}$.  

(ii)\ \ If $\omega$ is a reducible connection, then the gauge orbit $\mathcal G_P^{H^{s+1}}\cdot\omega$ is a singular orbit of the action 
$\mathcal G_P^{H^{s+1}}\curvearrowright\mathcal A_P^{H^s}$ and the dimension of regular orbits of the slice representation of the action 
$\mathcal G_P^{H^{s+1}}\curvearrowright(\mathcal A_P^{H^s},g_s)$ at $\omega$ is equal to ${\rm dim}({\rm Ker}\,d_{\omega})$.  
}

\vspace{0.25truecm}


Denote by $\widehat{\mathcal G}$ the (full) isometry group of $(\mathcal A_P^{H^s},g_s)$.  
According to Proposition 8.2, the gauge transformation group $\mathcal G_P^{H^{s+1}}$ is regarded as a subgroup of $\widehat G$.  
Denote by ${\rm Aff}(\mathcal A_P^{H^s})$ the affine transformation group of the affine Hilbert space $\mathcal A_P^{H^s}$.  
The transformation ${\bf g}$ of $(\mathcal A_P^{H^s},g_s)$ is described as 
$${\bf g}\cdot\omega-\omega_0={\rm Ad}(\widehat{\bf g})(\omega-\omega_0)+{\rm Ad}(\widehat{\bf g})(\omega_0)
-(R_{\widehat{\bf g}(\cdot)})_{\ast}^{-1}\circ\widehat{\bf g}_{\ast\cdot}.$$
Thus ${\bf g}$ is an affine transformation of $\mathcal A_P^{H^s}$.  Hence $\mathcal G_P^{H^{s+1}}$ is regarded as a subgroup of 
$\widehat G\cap{\rm Aff}(\mathcal A_P^{H^s})$.  

Denote by ${\rm Diff}^{H^s}(\mathcal A_P^{H^s})$ the Hilbert Lie group of all $H^s$-diffeomorphisms of $\mathcal A_P^{H^s}$.  
From the definition of $g_s$, we can derive the following fact directly.  

\vspace{0.25truecm}

\noindent
{\bf Proposition 8.4.}\ {\sl Define a group $\widehat G'$ by 
\begin{align*}
\widehat G':=\{&F\in{\rm Diff}^{H^s}(\mathcal A_P^{H^s})\,\vert\,\\
&\langle(\widehat{F_{\ast\omega}(A_1)})_x,\,(\widehat{\square_{F(\omega)}^s(F_{\ast\omega}(A_2))})_x\rangle_{B,\mathfrak g}
=\langle(\widehat A_1)_x,\,(\widehat{\square_{\omega}^s(A_2)})_x\rangle_{B,\mathfrak g}\\
&\hspace{1truecm}(x\in B,\omega\in\mathcal A_P^{H^s},A_i\in T_{\omega}\mathcal A_P^{H^s}(\approx\Omega_{\mathcal T,1}^{H^s}(P,\mathfrak g))\,\,(i=1,2))\}.
\end{align*}
Then we have 
$$\mathcal G_P^{H^{s+1}}\subset\widehat G'\cap{\rm Aff}(\mathcal A_P^{H^s})\subset\widehat G'\subset\widehat G.$$
}

\section{The orbit geometry of the action $\mu_c^{-1}(\hat{\bf 0})\curvearrowright(\mathcal A_P^{H^s},g_s)$} 
In this section, we shall investigate the orbit geometry of the action $\mu_c^{-1}(\hat{\bf 0})\curvearrowright(\mathcal A_P^{H^s},g_s)$.  
Take $\omega\in\mathcal A_P^{H^s},\,\,A\in\Omega_{\mathcal T,1}^{H^s}(P,\mathfrak g)$ and $B\in\Omega_{\mathcal T,i}^{H^s}(P,\mathfrak g)$ ($i=0,1$).  
Define $\mathcal T_{\omega,A,B}\in\Gamma^{H^s}(T^{(0,i+1)}P\otimes\mathfrak g)$ by 
$$
\left\{\begin{array}{ll}
\displaystyle{(\mathcal T_{\omega,A,B})_u(\vv):=dB_u((\omega_u|_{(\mathcal V^{\mathcal A})_u})^{-1}(A_u(\vv_{\mathcal H^{\omega}_u}))} & (i=0)\\
\displaystyle{(\mathcal T_{\omega,A,B})_u(\vv_1,\vv_2):=dB_u((\omega_u|_{(\mathcal V^{\mathcal A})_u})^{-1}(A_u((\vv_1)_{\mathcal H^{\omega}_u})),\,(\vv_2)_{\mathcal H^{\omega}_u})} & 
(i=1),
\end{array}\right.$$
($u\in P,\,\,\vv,\vv_1,\vv_2\in T_uP$).  
Define a differential operator $\mathcal D_{\omega,A}:\Omega_{\mathcal T,i}(P,\mathfrak g)\to\Omega_{\mathcal T,i+1}(P,\mathfrak g)$ ($i=0,1$) of order one by 
$$\mathcal D_{\omega,A}(B):=\left\{\begin{array}{ll}
\displaystyle{-\int_0^1\mathcal T_{\omega+tA,A,B}\,dt} & (i=0)\\
\displaystyle{-2\int_0^1{\rm Alt}(\mathcal T_{\omega+tA,A,B})\,dt} & (i=1),
\end{array}\right.$$
where ${\rm Alt}(\cdot)$ denotes the alternization of $(\cdot)$.  
Denote by $\mathcal D_{\omega,A}^{\ast}$ the adjoint operator of $\mathcal D_{\omega,A}$.  
By using these operators, we define a differential operator $\Delta_{\omega,A}$ of order two by 
\begin{align*}
\Delta_{\omega,A}:=&(d_{\omega}\circ\mathcal D_{\omega,A}^{\ast}+\mathcal D_{\omega,A}^{\ast}\circ d_{\omega})
+(\mathcal D_{\omega,A}\circ d_{\omega}^{\ast}+d_{\omega}^{\ast}\circ\mathcal D_{\omega,A})\\
&+(\mathcal D_{\omega,A}\circ\mathcal D_{\omega,A}^{\ast}+\mathcal D_{\omega,A}^{\ast}\circ\mathcal D_{\omega,A}).
\end{align*}
For $A\in\Omega_{\mathcal T,1}^{H^s}(P,\mathfrak g)$, we define by a map $\tau_A$ by 
$$\tau_A(\omega):=\omega+A\quad\,\,(\omega\in\mathcal A_P^{H^s}).$$
The additive group $\Omega_{\mathcal T,1}^{H^s}(P,\mathfrak g)$ acts on $\mathcal A_P^{H^s}$ by 
$$A\cdot\omega:=\tau_A(\omega)\quad\,\,(A\in\Omega_{\mathcal T,1}^{H^s}(P,\mathfrak g)).$$
For the action $\Omega_{\mathcal T,1}^{H^s}(P,\mathfrak g)\curvearrowright(\mathcal A_P^{H^s},g_s)$, we have the following fact.  

\vspace{0.35truecm}

\noindent
{\bf Proposition 9.1.} {\sl For $\omega\in\mathcal A_P^{H^s}$ and $A\in\Omega_{\mathcal T,1}^{H^s}(P,\mathfrak g)$, we have 
$$\begin{array}{l}
\hspace{0.5truecm}\displaystyle{(\tau_A^{\ast}g_s)_{\omega}(B_1,B_2)-(g_s)_{\omega}(B_1,B_2)}\\
\displaystyle{=\langle B_1,((\square_{\omega}+\Delta_{\omega,A})^s-\square_{\omega}^s)(B_2)\rangle_0}
\end{array}\leqno{(9.1)}$$
($B_1,B_2\in\Omega_{\mathcal T,1}^{H^s}(P,\mathfrak g)$).  
} 

\vspace{0.35truecm}

\noindent
{\it Proof.}\ \ We consider the case of $i=1$.  
Take $\omega\in\mathcal A_P^{H^s}$, $A,B\in\Omega_{\mathcal T,1}^{H^s}(P,\mathfrak g)$, $u\in P$ and $\vv_1,\vv_2\in T_uP$.  
Then, from $(11.4)$ (see Section 11), we have 
\begin{align*}
&\left.\frac{d}{dt}\right|_{t=t_0}(d_{\omega+tA}B)_u(\vv_1,\vv_2)\\
=&-dB_u\left(((\omega+t_0A)_u|_{\mathcal V_u})^{-1}(A_u((\vv_1)_{\mathcal H^{\omega+t_0A}}),(\vv_2)_{\mathcal H^{\omega+t_0A}})\,\right)\\
&-dB_u\left((\vv_1)_{\mathcal H^{\omega+t_0A}},\,((\omega+t_0A)_u|_{\mathcal V_u})^{-1}(A_u((\vv_2)_{\mathcal H^{\omega+t_0A}})\,\right)
\end{align*}
($0\leq t_0\leq 1$).  By integrating the both sides of this relation from $0$ to $1$ with respect to $t_0$, we obtain 
$$(d_{\omega+A}B)_u(\vv_1,\vv_2)=(d_{\omega}B)_u(\vv_1,\vv_2)+(\mathcal D_{\omega,A}(B))_u(\vv_1,\vv_2).$$
Thus we obtain 
$$d_{\omega+A}=d_{\omega}+\mathcal D_{\omega,A}.\leqno{(9.2)}$$
Also, we obtain 
$$d^{\ast}_{\omega+A}=d^{\ast}_{\omega}+\mathcal D^{\ast}_{\omega,A},\leqno{(9.3)}$$
where $\mathcal D^{\ast}_{\omega,A}$ denotes the adjoint operator of $\mathcal D_{\omega,A}$.  
From $(9.2)$ and $(9.3)$, we can derive 
$$\Delta_{\omega+A}=\Delta_{\omega}+\Delta_{\omega,A},\leqno{(9.4)}$$
where $\Delta_{\omega}$ and $\Delta_{\omega+A}$ are the Hodge-de Rham Laplace operators with respect to $\omega$ and $\omega+A$, respectively.  
Similarly, in the case of $i=0$, we can show the relation $(9.4)$.  
Take $B_1,B_2\in\Omega_{\mathcal T,1}^{H^s}(P,\mathfrak g)$.  By using these relations, we can derive 
\begin{align*}
(\tau_A^{\ast}g_s)_{\omega}(B_1,B_2)=&(g_s)_{\omega+A}((\tau_A)_{\ast\omega}(B_1),(\tau_A)_{\ast\omega}(B_2))=(g_s)_{\omega+A}(B_1,B_2)\\
=&(g_s)_{\omega}(B_1,B_2)+\langle B_1,((\square_{\omega}+\Delta_{\omega,A})^s-\square_{\omega}^s)(B_2)\rangle_0.
\end{align*}
\qed

\vspace{0.35truecm}

The additive subgroup $\mu_c^{-1}(\hat{\bf 0})$ acts on $\mathcal A_P^{H^s}$ as a subaction of this action.  

\vspace{0.35truecm}

\noindent
{\bf Proposition 9.2.} {\sl 
The action $\mu_c^{-1}(\hat{\bf 0})\curvearrowright(\mathcal A_P^{H^s},g_s)$ is not isometric but horizontally isometric.}

\vspace{0.35truecm}

\noindent
{\it Proof.}\ \ According to Proposition 9.1, the action $\mu_c^{-1}(\hat{\bf 0})\curvearrowright(\mathcal A_P^{H^s},g_s)$ is not isometric.  
On the other hand, since $\mu_c$ is the orbit map of this action and 
$$\mu_c:(\mathcal A_P^{H^s},{\rm g}_s)\to(H^s([0,1],\mathfrak g),\langle\,\,,\,\,\rangle_s^o)$$
is a homothetic submersion by Proposition 6.2, we see that this action is horizontally isometric.  \qed

\section{Proof of Theorem B} 
In this section, we shall prove Theorem B stated in Introduction.  

\vspace{0.25truecm}

\noindent
{\it Proof of Theorem B.}\ \ 
According to Proposition 8.3, the subaction $\mathcal G_P^{H^{s+1}}\curvearrowright(\mathcal A_P^{H^s},g_s)$ is isometric.  
Also, according to Proposition 9.2, the subaction $\mu_c^{-1}(\hat{\bf 0})\curvearrowright(\mathcal A_P^{H^s},g_s)$ is horizontally isometric.  
Hence the action ${\widetilde{\mathcal G_P^{H^{s+1}}}}^c\curvearrowright(\mathcal A_P^{H^s},g_s)$ is horizontally isometric.  
Assume that $c$ is not a loop.  
Set $x_0:=c(0)$.  Take any $\omega_1,\,\omega_2\in\mathcal A_P^{H^s}$.  Set $u_i:=\mu_c(\omega_i)$ ($i=1,2$).  
Since $H^{s+1}([0,1],G)$ acts on $H^s([0,1],\mathfrak g)$ transitively, there exists $\bar{\bf g}\in H^{s+1}([0,1],G)$ with $\bar{\bf g}\cdot u_1=u_2$.  
It is clear that there exists ${\bf g}\in\mathcal G_P^{H^{s+1}}$ with $\lambda_c({\bf g})=\bar{\bf g}$.  
From $(6.7)$, we have 
$$\mu_c({\bf g}\cdot\omega_1)=\lambda_c({\bf g})\cdot\mu_c(\omega_1)=\bar{\bf g}\cdot u_1=u_2=\mu_c(\omega_2).$$
Thus, since the images of ${\bf g}\cdot\omega_1$ and $\omega_2$ by $\mu_c$ coincide, $A:=\overrightarrow{({\bf g}\cdot\omega_1)\,\omega_2}$ belongs to 
$\mu_c^{-1}(\hat{\bf 0})$.  
Hence we obtain $({\bf g},A)\cdot\omega_1=\omega_2$.  Therefore we see that ${\widetilde{\mathcal G_P^{H^{s+1}}}}^c$ acts on $\mathcal A_P^{H^s}$ transitively.  \qed

\section{Proof of Theorem C.} 
In this section, we shall prove Theorem C stated in Introduction.  
Denote by $\mathcal X(\mathcal A_P^{H^s})$ the space of all $C^{\infty}$-vector fields on $\mathcal A_P^{H^s}$.  
Note that $\mathcal X(\mathcal A_P^{H^s})$ is identified with the space $C^{\infty}(\mathcal A_P^{H^s},\Omega_{\mathcal T,1}^{H^s}(P,\mathfrak g))$ of all $C^{\infty}$-maps of 
$\mathcal A_P^{H^s}$ into $\Omega_{\mathcal T,1}^{H^s}(P,\mathfrak g)$.  
The correspondence 
$$d_{\bullet}:\omega\mapsto d_{\omega}(:\Omega_{\mathcal T,i}^{H^s}(P,\mathfrak g)\to\Omega_{\mathcal T,i+1}^{H^s}(P,\mathfrak g))\quad\,\,(\omega\in\mathcal A_P^{H^s})$$
is regarded as a $\Omega_{\mathcal T,i}^{H^s}(P,\mathfrak g)^{\ast}\otimes\Omega_{\mathcal T,i+1}^{H^s}(P,\mathfrak g)$-valued function over $\mathcal A_P^{H^s}$ 
and the correspondence 
$$d^{\ast}_{\bullet}:\omega\mapsto d^{\ast}_{\omega}(:\Omega_{\mathcal T,i}^{H^s}(P,\mathfrak g)\to\Omega_{\mathcal T,i-1}^{H^s}(P,\mathfrak g))\quad\,\,
(\omega\in\mathcal A_P^{H^s})$$
is regarded as a $\Omega_{\mathcal T,i}^{H^s}(P,\mathfrak g)^{\ast}\otimes\Omega_{\mathcal T,i-1}^{H^s}(P,\mathfrak g)$-valued function over $\mathcal A_P^{H^s}$.  
In particular, in the case of $i=1$, $d_{\bullet}$ (resp. $d^{\ast}_{\bullet}$) is regarded as 
a $\Omega_{\mathcal T,2}^{H^s}(P,\mathfrak g)$-valued (resp. $\Omega_{\mathcal T,0}^{H^s}(P,\mathfrak g)$-valued) $1$-form over $\mathcal A_P^{H^s}$ 
by identified $\Omega_{\mathcal T,1}^{H^s}(P,\mathfrak g)$ with $T_{\omega}\mathcal A_P^{H^s}$'s ($\omega\in\mathcal A_P^{H^s}$).  
Similarly the correspondence 
$$\square_{\bullet}:\omega\mapsto\square_{\omega}(:\Omega_{\mathcal T,1}^{H^s}(P,\mathfrak g)\to\Omega_{\mathcal T,1}^{H^s}(P,\mathfrak g))\quad\,\,
(\omega\in\mathcal A_P^{H^s})$$
also is regarded as a $\Omega_{\mathcal T,1}^{H^s}(P,\mathfrak g)$-valued $1$-form over $\mathcal A_P^{H^s}$.  

For $A_i\in\Omega_{\mathcal T,k_i}^{H^s}(P,\mathfrak g)$, we define a $\mathfrak g$-valued covariant tensor field $[A_1\otimes A_2]$ of degree $(k_1+k_2)$ on $P$ by 
$$\begin{array}{r}
\displaystyle{[A_1\otimes A_2]_u(\vv_1,\cdots,\vv_{k_1+k_2}):=[(A_1)_u(\vv_1,\cdots,\vv_{k_1}),A_2(\vv_{k_1+1},\cdots,\vv_{k_1+k_2})]}\\
\displaystyle{(u\in P,\,\,\vv_1,\cdots,\vv_{k_1+k_2}\in T_uP).}
\end{array}$$
Also, we define a $\mathfrak g$-valud $(k_1+k_2)$-form $[A_1\wedge A_2]$ on $P$ by 
$$\begin{array}{l}
\hspace{0.5truecm}\displaystyle{[A_1\wedge A_2]_u(\vv_1,\cdots,\vv_{k_1+k_2})}\\
\displaystyle{:=\frac{1}{(k_1+k_2)!}\,
\sum_{\sigma\in S_{k_1+k_2}}{\rm sgn}\,\sigma\,[(A_1)_u(\vv_{\sigma(1)},\cdots,\vv_{\sigma(k_1)}),A_2(\vv_{\sigma(k_1+1)},\cdots,\vv_{\sigma(k_1+k_2)})]}\\
\hspace{8truecm}\displaystyle{(u\in P,\,\,\vv_1,\cdots,\vv_{k_1+k_2}\in T_uP),}
\end{array}$$
where $S_{k_1+k_2}$ denotes the symmetric group of degree $(k_1+k_2)$.  

We can show the following relations for the directional derivatives of $d_{\bullet}$ and $d^{\ast}_{\bullet}$.  

\vspace{0.25truecm}

\noindent
{\bf Lemma 11.1.} {\sl Let $A,\,A_i\in T_{\omega}\mathcal A_P^{H^s}(=\Omega_{\mathcal T,1}^{H^s}(P,\mathfrak g))$ ($i=1,2$).  
Then the following statements {\rm(i), (ii)} and {\rm(iii)} hold.  

{\rm (i)}\ \ The directional derivative $A(d_{\bullet})$ of the $\Omega_{\mathcal T,0}^{H^s}(P,\mathfrak g)^{\ast}\otimes\Omega_{\mathcal T,1}^{H^s}(P,\mathfrak g)$-valued 
function $d_{\bullet}$ over $\mathcal A_P^{H^s}$ with respect to $A$ is given by 
$$A(d_{\bullet})(\eta)=-[\eta\otimes A]\qquad\,\,(\eta\in\Omega_{\mathcal T,0}^{H^s}(P,\mathfrak g)).$$

{\rm (ii)}\ \ The directional derivative $A_1(d_{\bullet})$ of the $\Omega_{\mathcal T,1}^{H^s}(P,\mathfrak g)^{\ast}\otimes\Omega_{\mathcal T,2}^{H^s}(P,\mathfrak g)$-valued 
function $d_{\bullet}$ over $\mathcal A_P^{H^s}$ with respect to $A_1$ is given by 
$$A_1(d_{\bullet})(A_2)=2[A_1\wedge A_2].$$

{\rm (iii)}\ \ The directional derivative $A_1(d^{\ast}_{\bullet})$ of 
the $\Omega_{\mathcal T,1}^{H^s}(P,\mathfrak g)^{\ast}\otimes\Omega_{\mathcal T,0}^{H^s}(P,\mathfrak g)$-valued function $d^{\ast}_{\bullet}$ over $\mathcal A_P^{H^s}$ 
with respect to $A_1$ is given by 
$$A_1(d^{\ast}_{\bullet})(A_2)=-{\rm Tr}_{g_S^{\omega}}^{\bullet}[A_1\otimes A_2](\bullet,\bullet).$$
where $g^{\omega}_S$ denotes the Sasaki metric of $P$ arising from $g_B$ and $\langle\,\,,\,\,\rangle_{\mathfrak g}$ and 
${\rm Tr}_{g^{\omega}_S}^{\bullet}$ denotes the trace in two covariant components $(\bullet)$'s with respect to $g^{\omega}_S$.  

{\rm (iv)}\ \ The directional derivative $A(d^{\ast}_{\bullet})$ of 
the $\Omega_{\mathcal T,2}^{H^s}(P,\mathfrak g)^{\ast}\otimes\Omega_{\mathcal T,1}^{H^s}(P,\mathfrak g)$-valued function $d^{\ast}_{\bullet}$ over $\mathcal A_P^{H^s}$ 
with respect to $A$ is given by 
$$A(d^{\ast}_{\bullet})(\varphi)={\bf 0}.$$

{\rm (v)}\ \ The directional derivative $A_1(\square_{\bullet})$ of 
the $\Omega_{\mathcal T,1}^{H^s}(P,\mathfrak g)^{\ast}\otimes\Omega_{\mathcal T,1}^{H^s}(P,\mathfrak g)$-valued function $\square_{\bullet}$ over $\mathcal A_P^{H^s}$ 
with respect to $A_1$ is given by 
$$\begin{array}{r}
\displaystyle{A_1(\square_{\bullet})(A_2)=-[d_{\omega}^{\ast}A_2\otimes A_1]-d_{\omega}\left({\rm Tr}_{g_S^{\omega}}^{\bullet}[A_1\otimes A_2](\bullet,\bullet)\right)
+2d_{\omega}^{\ast}([A_1\wedge A_2])}\\
\displaystyle{(A_1,A_2\in\Omega_{\mathcal T,1}^{H^s}(P,\mathfrak g)).}
\end{array}$$
}

\vspace{0.25truecm}

\noindent
{\it Proof.}\ \ 
Denote by $\mathcal V^P$ and $\mathcal H^{\omega}$ the vertical distribution and the horizontal distribution with respect to $\omega(\in\mathcal A_P^{H^s})$ of 
the $G$-bundle $\pi:P\to B$, respectively.  
Also, denote by $\nabla^B$ and $\nabla^{S,\omega}$ the Riemannian connections of $g_B$ and the Sasaki metric $g_S^{\omega}$ associated to $\omega$, respectively.  
Take $\vv_i\in T_uP$ ($i=1,2$).  Set $x:=\pi(u)$ and $\bar{\vv}_i:=\pi_{\ast u}(\vv_i)$.  
Let $\vV_i$ ($i=1,2$) be a local vector fields on a neighborhood of $x$ in $B$ with $(\nabla^B\vV_i)_x={\bf 0}$ and 
$\widetilde{\vV_i}^{\omega}$ be the horizontal lift of $\vV_i$ with respect to the connection $\omega$.  
Since $\pi:(P,g_S^{\omega})\to(B,g_B)$ is a Riemannian submersion and $\widetilde{\vV_i}^{\omega}$ is the horizontal lift of $\vV_i$ with respect to 
this Riemannian submersion, it is shown that $(\nabla^{S,\omega}_{\widetilde{\vV_1}^{\omega}}\widetilde{\vV_2}^{\omega})_{\mathcal H^{\omega}}$ equals to the horizontal lift 
of $\nabla^B_{\vV_i}\vV_j$ with respect to this Riemannian submersion, where $(\cdot)_{\mathcal H^{\omega}}$ denotes the horizontal component of $(\cdot)$ with respect to 
the connection $\omega$.  Hence, from $(\nabla^B_{\vV_1}\vV_2)_x={\bf 0}$, we have 
$$\left(\left(\nabla^{S,\omega}_{\widetilde{\vV_1}^{\omega}}\widetilde{\vV_2}^{\omega}\right)_{\mathcal H^{\omega}}\right)_u
=([\widetilde{\vV_1}^{\omega},\widetilde{\vV_2}^{\omega}]_{\mathcal H^{\omega}})_u={\bf 0}.\leqno{(11.1)}$$
First we shall show the relation in the statement (ii).  For $A_1(d_{\bullet})$, we have 
$$\begin{array}{l}
\hspace{0.5truecm}\displaystyle{((A_1(d_{\bullet}))\,(A_2))_u(\vv_1,\vv_2)}\\
\displaystyle{=\left.\frac{d}{dt}\right|_{t=0}(d_{\omega+tA_1}A_2)_u(\vv_1,\vv_2)}\\
\displaystyle{=\left.\frac{d}{dt}\right|_{t=0}(dA_2)_u((\vv_1)_{\mathcal H^{\omega+tA_1}},(\vv_2)_{\mathcal H^{\omega+tA_1}})}\\
\displaystyle{=(dA_2)_u\left(\left.\frac{d}{dt}\right|_{t=0}(\vv_1)_{\mathcal H^{\omega+tA_1}},(\vv_2)_{\mathcal H^{\omega}}\right)}\\
\hspace{0.5truecm}\displaystyle{+(dA_2)_u\left((\vv_1)_{\mathcal H^{\omega}},\left.\frac{d}{dt}\right|_{t=0}(\vv_2)_{\mathcal H^{\omega+tA_1}}\right)}\\
\displaystyle{=(dA_2)_u\left(\left.\frac{d}{dt}\right|_{t=0}\left((\widetilde{\vV}^{\omega}_1)_{\mathcal H^{\omega+tA_1}}\right)_u,\,
\left((\widetilde{\vV}^{\omega}_2)_{\mathcal H^{\omega}}\right)_u\right)}\\
\hspace{0.5truecm}\displaystyle{+(dA_2)_u\left(\left((\widetilde{\vV}^{\omega}_1)_{\mathcal H^{\omega}}\right)_u,\,
\left.\frac{d}{dt}\right|_{t=0}\left((\widetilde{\vV}^{\omega}_2)_{\mathcal H^{\omega+tA_1}}\right)_u,\,\right)}\\
\displaystyle{=\left(\left.\frac{d}{dt}\right|_{t=0}(\vv_1)_{\mathcal H^{\omega+tA_1}}\right)(\widetilde A_2(\widetilde{\vV}^{\omega}_2))}\\
\hspace{0.5truecm}\displaystyle{+(\vv_1)_{\mathcal H^{\omega}}\left(A_2\left(\left.\frac{d}{dt}\right|_{t=0}(\widetilde{\vV}^{\omega}_2)_{\mathcal H^{\omega+tA_1}}\right)\right)}\\
\hspace{0.5truecm}\displaystyle{-\left(\left.\frac{d}{dt}\right|_{t=0}(\vv_2)_{\mathcal H^{\omega+tA_1}}\right)(A_2(\widetilde{\vV}^{\omega}_1))}\\
\hspace{0.5truecm}\displaystyle{-(\vv_2)_{\mathcal H^{\omega}}\left(A_2\left(\left.\frac{d}{dt}\right|_{t=0}(\widetilde{\vV}^{\omega}_1)_{\mathcal H^{\omega+tA_1}}\right)
\right)}\\
\hspace{0.5truecm}\displaystyle{-(A_2)_u\left(\left[\left.\frac{d}{dt}\right|_{t=0}(\widetilde{\vV}^{\omega}_1)_{\mathcal H^{\omega+tA_1}},\,
(\widetilde{\vV}^{\omega}_2)_{\mathcal H^{\omega}}\right]_u\right)}\\
\hspace{0.5truecm}\displaystyle{
-(A_2)_u\left(\left[(\widetilde{\vV}^{\omega}_1)_{\mathcal H^{\omega}},\,\left.\frac{d}{dt}\right|_{t=0}(\widetilde{\vV}^{\omega}_2)_{\mathcal H^{\omega+tA_1}}\right]_u\right).}
\end{array}\leqno{(11.2)}$$
It is easy to show that $\displaystyle{\left.\frac{d}{dt}\right|_{t=0}(\widetilde{\vV}^{\omega}_i)_{\mathcal H^{\omega+tA_1}}}$ ($i=1,2$) are local sections of $\mathcal V^P$.  
Also, since $\widetilde{\vV}^{\omega}_j$ is $\pi$-projectable, 
$$\left[\left.\frac{d}{dt}\right|_{t=0}(\widetilde{\vV}^{\omega}_i)_{\mathcal H^{\omega+tA_1}},(\widetilde{\vV}^{\omega}_j)_{\mathcal H^{\omega}}\right]\,\,\,\,
((i,j)=(1,2)\,\,{\rm or}\,\,(2,1))$$
are local sections of $\mathcal V^P$.  Hence we obtain 
$$(A_2)_u\left(\left(\left.\frac{d}{dt}\right|_{t=0}(\widetilde{\vV}^{\omega}_i)_{\mathcal H^{\omega+tA_1}}\right)_u\right)={\bf 0}$$
and 
$$(A_2)_u\left(\left[\left.\frac{d}{dt}\right|_{t=0}(\widetilde{\vV}^{\omega}_i)_{\mathcal H^{\omega+tA_1}},\,(\widetilde{\vV}^{\omega}_j)_{\mathcal H^{\omega}}\right]_u\right)={\bf 0}.$$
From $(11.2)$ and these relations, we obtain 
$$\begin{array}{l}
\hspace{0.5truecm}\displaystyle{(A_1(d_{\bullet})\,(A_2))_u(\vv_1,\vv_2)}\\
\displaystyle{=\left(\left.\frac{d}{dt}\right|_{t=0}(\vv_1)_{\mathcal H^{\omega+tA_1}}\right)(A_2(\widetilde{\vV}^{\omega}_2))
-\left(\left.\frac{d}{dt}\right|_{t=0}(\vv_2)_{\mathcal H^{\omega+tA_1}}\right)(A_2(\widetilde{\vV}^{\omega}_1)).}
\end{array}\leqno{(11.3)}$$
Define $\Psi_u:G\to\pi^{-1}(x)$ by $\Psi_u(g):=ug$ ($g\in G$).  
Let $\alpha_i:(-\varepsilon,\varepsilon)\to\pi^{-1}(x)$ ($i=1,2$) be $C^{\infty}$-curves with 
$\displaystyle{\alpha_i'(0)=\left.\frac{d}{dt}\right|_{t=0}(\vv_i)_{\mathcal H^{\omega+tA_1}}}$ and 
$g_i:(-\varepsilon,\varepsilon)\to G$ be the $C^{\infty}$-curves with $\alpha_i(\bar t)=ug_i(\bar t)$ ($\bar t\in(-\varepsilon,\varepsilon)$).  
Then we have 
$$\begin{array}{l}
\hspace{0.5truecm}\displaystyle{\left(\left.\frac{d}{dt}\right|_{t=0}(\vv_1)_{\mathcal H^{\omega+tA_1}}\right)\left(A_2(\widetilde{\vV}^{\omega}_2)\right)}\\
\displaystyle{=\left.\frac{d}{d\bar t}\right|_{\bar t=0}(A_2)_{\alpha_1(\bar t)}\left((\widetilde{\vV}^{\omega}_2)_{\alpha_1(\bar t)}\right)}\\
\displaystyle{=\left.\frac{d}{d\bar t}\right|_{\bar t=0}(A_2)_{\alpha_1(\bar t)}\left((R_{{\bf g}_1(\bar t)})_{\ast}((\vv_2))_{\mathcal H^{\omega}})\right)}\\
\displaystyle{=\left.\frac{d}{d\bar t}\right|_{\bar t=0}{\rm Ad}(g_1(\bar t)^{-1})((A_2)_u(\vv_2))}\\
\displaystyle{=\left(\left.\frac{d}{d\bar t}\right|_{\bar t=0}{\rm Ad}(g_1(\bar t)^{-1})\right)((A_2)_u(\vv_2))}\\
\displaystyle{=-{\rm ad}(g_1'(0))((A_2)_u(\vv_2))}\\
\displaystyle{=-\left[(\Psi_u)_{\ast u}^{-1}\left(\left.\frac{d}{dt}\right|_{t=0}(\vv_1)_{\mathcal H^{\omega+tA_1}}\right),\,(A_2)_u(\vv_2)\right].}
\end{array}$$
Also, we have 
\begin{align*}
0=&\left.\frac{d}{dt}\right|_{t=0}\left((\omega+tA_1)_u\left((\vv_1)_{\mathcal H^{\omega+tA_1}}\right)\right)\\
=&\omega_u\left(\left.\frac{d}{dt}\right|_{t=0}(\vv_1)_{\mathcal H^{\omega+tA_1}}\right)+(A_1)_u(\vv_1).
\end{align*}
Hence, by noticing $(\Psi_u)_{\ast u}^{-1}=\omega_u$, we have 
$$(\Psi_u)_{\ast u}^{-1}\left(\left.\frac{d}{dt}\right|_{t=0}(\vv_1)_{\mathcal H^{\omega+tA_1}}\right)=-(A_1)_u(\vv_1).\leqno{(11.4)}$$
Therefore we obtain 
$$\left(\left.\frac{d}{dt}\right|_{t=0}(\vv_1)_{\mathcal H^{\omega+tA_1}}\right)(A_2(\widetilde{\vV}^{\omega}_2))
=\left[(A_1)_u(\vv_1),(A_2)_u(\vv_2)\right].\leqno{(11.5)}$$
Similarly, we have 
$$
\left(\left.\frac{d}{dt}\right|_{t=0}(\vv_2)_{\mathcal H^{\omega+tA_1}}\right)(A_2(\widetilde{\vV}^{\omega}_1))
=\left[(A_1)_u(\vv_2),(A_2)_u(\vv_1)\right].
\leqno{(11.6)}$$
From $(11.3),\,(11.5)$ and $(11.6)$, we obtain 
$$(A_1(d_{\bullet})(A_2))_u(\vv_1,\vv_2)=2[A_1\wedge A_2]_u(\vv_1,\vv_2).\leqno{(11.7)}$$
Thus the relation in the statement (ii) is derived.  
Similarly, for $\eta\in\Omega_{\mathcal T,0}^{H^s}(P,\mathfrak g)$, we can show 
$$(A_1(d_{\bullet})(\eta))_u(\vv)=-[\eta\otimes A_1]_u(\vv),\leqno{(11.8)}$$
that is, the relation in the statement (i).  

Take $\omega\in\mathcal A_P^{H^s}$, $x\in B$ and $u_x\in\pi^{-1}(x)$.  
Let $\vV_i$ ($i=1,2$) and $\widetilde{\vV_i}^{\omega}$ be as above.  
From $(11.1)$, we have 
$$(A_i)_{u_x}([\widetilde{\vV_1}^{\omega},\widetilde{\vV_2}^{\omega}]_{u_x})
=(A_i)_{u_x}((\nabla^{S,\omega}_{\widetilde{\vV_1}^{\omega}}\widetilde{\vV_2}^{\omega})_{u_x})={\bf 0}.\leqno{(11.9)}$$
Hence we have 
$$\begin{array}{l}
\displaystyle{(d_{\omega}A_i)_{u_x}((\widetilde{\vV_1}^{\omega})_{u_x},(\widetilde{\vV_2}^{\omega})_{u_x})
=(\widetilde{\vV_1}^{\omega})_{u_x}(A_i(\widetilde{\vV_2}^{\omega}))-(\widetilde{\vV_2}^{\omega})_{u_x}(A_i(\widetilde{\vV_1}^{\omega})).}
\end{array}\leqno{(11.10)}$$
From $(11.9)$ and $(11.10)$, we can derive 
$$\begin{array}{l}
\displaystyle{(dA_i)_{u_x}((\widetilde{\vV_1}^{\omega})_{u_x},(\widetilde{\vV_2}^{\omega})_{u_x})
=\left(\nabla^{S,\omega}_{(\widetilde{\vV_1}^{\omega})_{u_x}}A_i\right)((\widetilde{\vV_2}^{\omega})_{u_x})}\\
\hspace{4.95truecm}\displaystyle{-\left(\nabla^{S,\omega}_{(\widetilde{\vV_2}^{\omega})_{u_x}}A_i\right)((\widetilde{\vV_1}^{\omega})_{u_x}).}
\end{array}\leqno{(11.11)}$$

Let $(\ve_i)_{i=1}^n$ be a local orthonormal frame field on a neighborhood of $x$ in $B$ with $(\nabla^B\ve_i)_x={\bf 0}$ ($i=1,\cdots,n$) 
and $\widetilde{\ve_i}^{\omega}$ be the horizontal lift of $\ve_i$ with respect to the connection $\omega$.  
Note that $(11.1),\,(11.9)$ and $(11.11)$ hold for $\widetilde{\ve_i}^{\omega}$'s instead of $\widetilde{\vV_i}^{\omega}$'s.  
Let $\widetilde A_i$ be the constant vector fields on $\mathcal A_P^{H^s}$ defined by $(\widetilde A_i)_{\omega}:=A_i$ ($\omega\in\mathcal A_P^{H^s}$) and 
$\widetilde{\eta}$ be the $\Omega_{\mathcal T,0}^{H^s}(P,\mathfrak g)$-valued constant function on $\mathcal A_P^{H^s}$ defined by 
$\widetilde{\eta}_{\omega}:=\eta$ ($\omega\in\mathcal A_P^{H^s}$).  
Then, by using noticing $A_1(\widetilde A_2)={\bf 0}$ and $A_1(\widetilde{\eta})={\bf 0}$, we can derive 
$$\begin{array}{l}
\displaystyle{\left\langle A_1(d^{\ast}_{\bullet})(A_2),\,\eta\right\rangle_0
=A_1\left(\langle d^{\ast}_{\bullet}(\widetilde A_2),\widetilde{\eta}\rangle_0\right)-\langle d^{\ast}_{\omega}(A_2),A_1(\widetilde{\eta})\rangle_0}\\
\hspace{2.85truecm}\displaystyle{=A_1\left(\langle\widetilde{A_2},\,d_{\bullet}(\widetilde{\eta})\rangle_0\right)}\\
\hspace{2.85truecm}\displaystyle{=\langle A_1(\widetilde A_2),d_{\omega}(\eta)\rangle_0+\langle A_2,A_1(d_{\bullet})(\eta)\rangle_0}\\
\hspace{2.85truecm}\displaystyle{=\langle A_2,A_1(d_{\bullet})(\eta)\rangle_0.}
\end{array}\leqno{(11.12)}$$
Also, from $(11.8)$, we can derive 
$$\begin{array}{l}
\hspace{0.5truecm}\displaystyle{\left\langle(\widehat{A_2})_x,\,\,\widehat{A_1(d_{\bullet})(\eta)}_x\right\rangle_{B,\mathfrak g}}\\
\displaystyle{=\sum_{i=1}^n\langle(A_2)_{u_x}((\widetilde{\ve_i}^{\omega})_{u_x}),\,
(A_1(d_{\bullet})(\eta)))_{u_x}((\widetilde{\ve_i}^{\omega})_{u_x})\rangle_{\mathfrak g}}\\
\displaystyle{=-\sum_{i=1}^n\langle(A_2)_{u_x}((\widetilde{\ve_i}^{\omega})_{u_x}),\,[\eta\otimes A_1]_{u_x}((\widetilde{\ve_i}^{\omega})_{u_x})\rangle_{\mathfrak g}}\\
\displaystyle{=-\sum_{i=1}^n\langle[A_1\otimes A_2]_{u_x}((\widetilde{\ve}_i^{\omega})_{u_x},\,(\widetilde{\ve}_i^{\omega})_{u_x}),\eta_{u_x}\rangle_{\mathfrak g}}\\
\displaystyle{=-\langle{\rm Tr}_{g_S^{\omega}}^{\bullet}[A_1\otimes A_2]_{u_x}(\bullet,\bullet),\,\eta_{u_x}\rangle_{\mathfrak g}}\\
\displaystyle{=-\langle\left({\rm Tr}_{g_S^{\omega}}^{\bullet}[A_1\otimes A_2]_{u_x}(\bullet,\bullet)\right)^{\wedge}_x,\,\widehat{\eta}_x\rangle_{B,\mathfrak g}.}
\end{array}\leqno{(11.13)}$$
From $(11.12)$ and $(11.13)$, we can derive 
$$\langle A_1(d^{\ast}_{\bullet})(A_2),\eta\rangle_0=
-\langle{\rm Tr}_{g_S^{\omega}}^{\bullet}[A_1\otimes A_2](\bullet,\bullet),\,\eta\rangle_0$$
and hence 
$$A_1(d^{\ast}_{\bullet})(A_2)=-{\rm Tr}_{g_S^{\omega}}^{\bullet}[A_1\otimes A_2](\bullet,\bullet).$$
Thus the relation in the statement (iii) is derived.  

Let $\widetilde{\varphi}$ be the $\Omega_{\mathcal T,2}^{H^s}(P,\mathfrak g)$-valued constant function on $\mathcal A_P^{H^s}$ defined by 
$\widetilde{\varphi}_{\omega}:=\varphi$ ($\omega\in\mathcal A_P^{H^s}$).  
Then, by noticing $A_1(\widetilde A_2)={\bf 0}$ and $A_1(\widetilde{\varphi})={\bf 0}$, we can derive 
{\small
$$\begin{array}{l}
\displaystyle{\left\langle A_1(d^{\ast}_{\bullet})(\varphi),\,A_2\right\rangle_0
=A_1\left(\langle d^{\ast}_{\bullet}(\widetilde{\varphi}),\widetilde A_2\rangle_0\right)-\langle d^{\ast}_{\omega}(\varphi),A_1(\widetilde A_2)\rangle_0}\\
\hspace{2.85truecm}\displaystyle{=A_1\left(\langle\widetilde{\varphi},\,d_{\bullet}(\widetilde A_2)\rangle_0\right)}\\
\hspace{2.85truecm}\displaystyle{=\langle A_1(\widetilde{\varphi}),d_{\omega}(A_2)\rangle_0+\langle\varphi,A_1(d_{\bullet})(A_2)\rangle_0}\\
\hspace{2.85truecm}\displaystyle{=\langle\varphi,A_1(d_{\bullet})(A_2)\rangle_0.}
\end{array}
$$
}
Also, from $(11.7)$, we can derive 
$$\begin{array}{l}
\hspace{0.5truecm}\displaystyle{\left\langle\widehat{\varphi}_x,\,\,(A_1(d_{\bullet})(A_2))^{\wedge}_x\right\rangle_{B,\mathfrak g}}\\
\displaystyle{=\sum_{i=1}^n\sum_{j=1}^n\langle\varphi_{u_x}((\widetilde{\ve_i}^{\omega})_{u_x},(\widetilde{\ve_j}^{\omega})_{u_x}),\,
(A_1(d_{\bullet})(A_2))_{u_x}((\widetilde{\ve_i}^{\omega})_{u_x},(\widetilde{\ve_j}^{\omega})_{u_x})\rangle_{\mathfrak g}}\\
\displaystyle{=2\sum_{i=1}^n\sum_{j=1}^n\langle\varphi_{u_x}((\widetilde{\ve_i}^{\omega})_{u_x},(\widetilde{\ve_j}^{\omega})_{u_x}),\,}\\
\hspace{3.25truecm}\displaystyle{[A_1\wedge A_2]_{u_x}((\widetilde{\ve_i}^{\omega})_{u_x},(\widetilde{\ve_j}^{\omega})_{u_x})\rangle_{\mathfrak g}}\\
\displaystyle{={\bf 0}.}
\end{array}
$$
From these relations, we can derive 
$$\langle A_1(d^{\ast}_{\bullet})(\varphi),A_2\rangle_0={\bf 0}$$
and hence 
$$A_1(d^{\ast}_{\bullet})(\varphi)={\bf 0}.$$
Thus the relation in the statement (iv) is derived.  
From the relations in the statements ${\rm(i)}-{\rm(iv)}$ and 
$$A_1(\square_{\bullet})=A_1(d_{\bullet})\circ d^{\ast}_{\omega}+d_{\omega}\circ A_1(d^{\ast}_{\bullet})
+A_1(d^{\ast}_{\bullet})\circ d_{\omega}+d^{\ast}_{\omega}\circ A_1(d_{\bullet}),$$
we can derive the relation in the statement (v) directly.  \qed

\vspace{0.25truecm}

We can show the following fact for $d_{\omega}^{\ast}:\Omega_{\mathcal T,1}^{H^s}(P,\mathfrak g)\to\Omega_{\mathcal T,0}^{H^s}(P,\mathfrak g)$.  

\vspace{0.25truecm}

\noindent
{\bf Proposition 11.2.} {\sl For $A\in\Omega_{\mathcal T,1}^{H^s}(P,\mathfrak g)$, we have 
$$d_{\omega}^{\ast}A=-{\rm Tr}_{g_S^{\omega}}^{\bullet}(\nabla^{S,\omega}A)(\bullet,\bullet).$$
}

\vspace{0.25truecm}

\noindent
{\it Proof.}\ \ By using the divergence theorem, we can derive 
$$\begin{array}{l}
\hspace{0.5truecm}\displaystyle{\langle d_{\omega}^{\ast}A,\eta\rangle_0=\langle A,d_{\omega}\eta\rangle_0}\\
\hspace{0truecm}\displaystyle{=\int_{x\in B}\langle\widehat A_x,(\widehat{d_{\omega}\eta})_x\rangle_{B,\mathfrak g}\,dv_B}\\
\hspace{0truecm}\displaystyle{=\int_{x\in B}\sum_{i=1}^n\langle(A_{u_x}((\widetilde{\ve}_i^{\omega})_{u_x}),\,(\widetilde{\ve}_i^{\omega})_{u_x}(\eta)
\rangle_{\mathfrak g}\,dv_B}\\
\hspace{0truecm}\displaystyle{=\int_{x\in B}\sum_{i=1}^n(\widetilde{\ve}_i^{\omega})_{u_x}\left(\langle A(\widetilde{\ve}_i^{\omega}),\,\eta\rangle_{\mathfrak g}\right)\,dv_B}\\
\hspace{0.5truecm}\displaystyle{-\int_{x\in B}\sum_{i=1}^n\langle(\nabla^{S,\omega}_{(\widetilde{\ve}_i^{\omega})_{u_x}}A)((\widetilde{\ve}_i^{\omega})_{u_x}),\,
\eta_{u_x}\rangle_{\mathfrak g}\,dv_B}\\
\hspace{0.5truecm}\displaystyle{-\int_{x\in B}\sum_{i=1}^n\langle A_{u_x}(\nabla^{S,\omega}_{(\widetilde{\ve}_i^{\omega})_{u_x}}\widetilde{\ve}_i^{\omega}),\,\eta_{u_x}
\rangle_{\mathfrak g}\,dv_B}\\
\hspace{0truecm}\displaystyle{=\int_{x\in B}({\rm div}_{\nabla^B}\vX_{A,\eta})_x\,dv_B
-\int_{x\in B}\langle({\rm Tr}_{g_S^{\omega}}^{\bullet}(\nabla^{S,\omega}A)_{u_x}(\bullet,\bullet),\,\eta_{u_x}\rangle_{\mathfrak g}\,dv_B}\\
\displaystyle{=-\langle{\rm Tr}_{g_S^{\omega}}^{\bullet}(\nabla^{S,\omega}A)(\bullet,\bullet),\,\eta\rangle_0\hspace{2.5truecm}
(A\in T_{\omega}\mathcal A_P^{H^s}(\approx\Omega_{\mathcal T,1}^{H^s}(P,\mathfrak g))),}
\end{array}$$
where $u_x$ and $\widetilde{\ve}_i^{\omega}$ are as in the proof of Lemma 11.1, ${\rm div}_{\nabla^B}(\cdot)$ is the divergence of $(\cdot)$ with respect to 
the Riemannian connection $\nabla^B$ of $g_B$ and $\vX_{A,\eta}$ is a vector field on $B$ defined by 
$$(g_B)_x((\vX_{A,\eta})_x,\vv)=\langle A_{u_x}(\widetilde{\vv}^{\omega}_{u_x}),\eta_{u_x}\rangle_{\mathfrak g}\quad\,\,(x\in B,\,\vv\in T_xB).$$
Hence we obtain the desired relation.  \qed

\vspace{0.25truecm}

We can show the following fact for the $\Omega_{\mathcal T,1}^{H^s}(P,\mathfrak g)$-valued $1$-form $\square_{\bullet}$ on $\mathcal A_P^{H^s}$.  

\vspace{0.25truecm}

\noindent
{\bf Lemma 11.3.} {\sl The exterior derivative $d\square_{\bullet}$ of the $\Omega_{\mathcal T,1}^{H^s}(P,\mathfrak g)$-valued $1$-form $\square_{\bullet}$ 
on $\mathcal A_P^{H^s}$ is given by 
$$\begin{array}{l}
\displaystyle{(d\square_{\bullet})_{\omega}(A_1,A_2)=
-\sum_{\sigma\in S_2}{\rm sgn}\,\sigma\,{\rm Tr}_{g_S^{\omega}}^{\bullet}\left[(\nabla^{S,\omega}A_{\sigma(1)})\otimes A_{\sigma(2)}\right](\bullet,\bullet,\cdot)}\\
\hspace{3truecm}\displaystyle{
-2\sum_{\sigma\in S_2}{\rm sgn}\,\sigma\,{\rm Tr}_{g_S^{\omega}}^{\bullet}
\left[(\nabla^{S,\omega}A_{\sigma(1)})\otimes A_{\sigma(2)}\right](\cdot,\bullet,\bullet)}\\
\hspace{4.25truecm}\displaystyle{(\omega\in\mathcal A_P^{H^s},\,\,A_1,\,A_2\in T_{\omega}\mathcal A_P^{H^s}(\approx\Omega_{\mathcal T,1}^{H^s}(P,\mathfrak g))).}
\end{array}$$
}

\vspace{0.25truecm}

\noindent
{\it Proof.}\ \ Let $\widetilde A_i$ ($i=1,2$) be the constant vector fields on $\mathcal A_P^{H^s}$ defined by $(\widetilde A_i)_{\omega}=A_i$ ($\omega\in\mathcal A_P^{H^s}$).  
Then, by noticing that $A_i(\widetilde A_j)={\bf 0}$ ($\leq i,j\leq 2$) and $[\widetilde A_1,\widetilde A_2]={\bf 0}$ hold, we have 
$$\begin{array}{l}
\displaystyle{(d\square_{\bullet})_{\omega}(A_1,A_2)
=A_1(\square_{\bullet}(\widetilde A_2))-A_2(\square_{\bullet}(\widetilde A_1))-\square_{\omega}([\widetilde A_1,\widetilde A_2]_{\omega})}\\
\hspace{2.7truecm}\displaystyle{=A_1(\square_{\bullet})(A_2)-A_2(\square_{\bullet})(A_1).}
\end{array}$$
By using this relation and the relation in (v) of Lemma 11.1, we can derive 
$$\begin{array}{l}
\displaystyle{(d\square_{\bullet})_{\omega}(A_1,A_2)=\sum_{\sigma\in S_2}{\rm sgn}\,\sigma\left[d_{\omega}^{\ast}(A_{\sigma(1)})\otimes A_{\sigma(2)}\right]}\\
\hspace{3truecm}\displaystyle{-2d_{\omega}\left({\rm Tr}_{g_S^{\omega}}^{\bullet}\left[A_1\otimes A_2\right](\bullet,\bullet)\right)}\\
\hspace{4.25truecm}\displaystyle{(\omega\in\mathcal A_P^{H^s},\,\,A_1,\,A_2\in T_{\omega}\mathcal A_P^{H^s}(\approx\Omega_{\mathcal T,1}^{H^s}(P,\mathfrak g))).}
\end{array}$$
Also we have 
$$\begin{array}{l}
\hspace{0.5truecm}\displaystyle{d_{\omega}\left({\rm Tr}_{g_S^{\omega}}^{\bullet}\left[A_1\otimes A_2\right](\bullet,\bullet)\right)}\\
\displaystyle{=\sum_{\sigma\in S_2}{\rm sgn}\,\sigma\,{\rm Tr}_{g_S^{\omega}}^{\bullet}
\left[(\nabla^{S,\omega}A_{\sigma(1)})\otimes A_{\sigma(2)}\right](\cdot,\bullet,\bullet)}
\end{array}$$
From these relations and Proposition 11.2, we can derive the desired relation.  \qed

\vspace{0.25truecm}

By using Lemma 11.3, we prove Theorem C stated in Introduction.  

\vspace{0.25truecm}

\noindent
{\it Proof of Theorem C.}\ \ 
The statement (i) follows from Lemma 11.3 directly.  W e shall show the statement (ii).  
Let $(\mathfrak g,\mathfrak k)$ be a symmetric pair and $\theta$ be an involution of $\mathfrak g$ with $({\rm Fix}\,\theta)_0\subset\mathfrak k\subset{\rm Fix}\,\theta$, where 
${\rm Fix}\,\theta$ denotes the fixed point group of $\theta$ and $({\rm Fix}\,\theta)_0$ denotes its identity component.  
Set $\mathfrak p:={\rm Ker}(\theta+{\rm id})$.  Take a maximal abelian subspace $\mathfrak a$ of $\mathfrak p$ and a Cartan subalgebra $\mathfrak h$ 
of the complexification $\mathfrak g^{\mathbb C}$ of $\mathfrak g$ including $\mathfrak a$.  
Let $\widetilde{\triangle}(\subset\mathfrak h^{\ast})$ be the root system of $\mathfrak g^{\mathbb C}$ with respect to $\mathfrak h$ and $\mathfrak g^{\mathbb C}_{\alpha}$ 
be the root space for $\alpha\in\widetilde{\triangle}$, where $\mathfrak h^{\ast}$ denotes the (complex) dual space of $\mathfrak h$.  
Then we have the following root space decomposition:
$$\mathfrak g^{\mathbb C}=\mathfrak h\oplus\left(\mathop{\oplus}_{\alpha\in\widetilde{\triangle}}\mathfrak g^{\mathbb C}_{\alpha}\right).$$
Set 
$$\widetilde{\triangle}^{\bullet}:=\{\alpha\in\widetilde{\triangle}\,|\,\alpha|_{\mathfrak a}={\bf 0}\}$$ 
and 
$$\triangle:=\{\sqrt{-1}\alpha|_{\mathfrak a}\,|\,\alpha\in\widetilde{\triangle}\setminus\widetilde{\triangle}^{\bullet}\},$$
which is a root system.  Let $\triangle_+(\subset\triangle)$ be the positive root system under some lexicographical ordering of $\mathfrak a$ (with respect to some basis of 
$\mathfrak a$).  Let $\triangle_+=\{\lambda_1,\cdots,\lambda_l\}$.  
Set 
$$\mathfrak g_{\lambda_a}:=\mathop{\oplus}_{\alpha\in\widetilde{\triangle}\,\,{\rm s.t.}\,\,\sqrt{-1}\alpha|_{\mathfrak a}=\lambda_a}
((\mathfrak g^{\mathbb C}_{\alpha}\oplus\mathfrak g^{\mathbb C}_{-\alpha})\cap\mathfrak g)$$
and 
$$\mathfrak g_0:=\mathfrak a\oplus\left(\mathop{\oplus}_{\alpha\in\widetilde{\triangle}^{\bullet}}((\mathfrak g^{\mathbb C}_{\alpha}\oplus
\mathfrak g^{\mathbb C}_{-\alpha})\cap\mathfrak g)\right).$$
Here we note that $\mathop{\oplus}_{\alpha\in\widetilde{\triangle}^{\bullet}}((\mathfrak g^{\mathbb C}_{\alpha}\oplus\mathfrak g^{\mathbb C}_{-\alpha})\cap\mathfrak g)$ is equal 
to the centralizer $\mathfrak z_{\mathfrak k}(\mathfrak a)$ of $\mathfrak a$ in $\mathfrak k$.  
Then we have the following decomposition:
$$\mathfrak g=\mathfrak g_0\oplus\left(\mathop{\oplus}_{a=1}^l\mathfrak g_{\lambda_a}\right).\leqno{(11.14)}$$
Also, set $n_0:={\rm dim}\,\mathfrak g_0$ and $n_a:={\rm dim}\,\mathfrak g_{\lambda_a}$ ($a=1,\cdots,l$).  
Let $\{\overline{\ve}_i\}_{i=1}^N$ be an orthonormal basis of $\mathfrak g$ with respect to the $(-1)$-multiple of the Killing form of $\mathfrak g$ satisfying 
$$\overline{\ve}_i\in\left\{\begin{array}{ll}
\mathfrak g_0 & (1\leq i\leq n_0)\\
\mathfrak g_{\lambda_1} & (n_0+1\leq i\leq n_0+n_1)\\
\mathfrak g_{\lambda_2} & (n_0+n_1+1\leq i\leq n_0+n_1+n_2)\\
\hspace{0.15truecm}\vdots & \qquad\qquad\vdots\\
\mathfrak g_{\lambda_l} & (\sum\limits_{a=0}^{l-1}n_a+1\leq i\leq N).
\end{array}\right.$$
Also, let $c_{ij}^k$ ($1\leq i,j,k\leq N$) be the structure constants of $\mathfrak g$ with respect to $\{\overline{\ve}_i\}_{i=1}^N$, that is, 
$[\overline{\ve}_i,\overline{\ve}_j]=\sum\limits_{k=1}^Nc_{ij}^k\overline{\ve}_k$ ($1\leq i,j\leq N$).  
The bracket products of $\overline{\ve}_i$'s ($i=1,\cdots,N$) satisfy 
$$[\overline{\ve_i},\overline{\ve}_j]\in\left\{\begin{array}{ll}
\mathfrak g_0 & (1\leq i\leq j\leq n_0)\\
\mathfrak g_{\lambda_1}\setminus\{{\bf 0}\} & (1\leq i\leq n_0,\,\,n_0+1\leq n_0+n_1)\\
\mathfrak g_{\lambda_2}\setminus\{{\bf 0}\} & (1\leq i\leq n_0,\,\,n_0+n_1+1\leq n_0+n_1+n_2)\\
\hspace{0.5truecm}\vdots & \hspace{2.5truecm}\vdots\\
\mathfrak g_{\lambda_l}\setminus\{{\bf 0}\} & (1\leq i\leq n_0,\,\,\sum_{a=0}^{l-1}n_a+1\leq j\leq N)\\
\mathfrak g_{\lambda_{b_1}+\lambda_{b_2}} & (\sum\limits_{a=0}^{b_1-1}n_a+1\leq i\leq\sum\limits_{a=0}^{b_1}n_a,\,\,\,\sum\limits_{a=0}^{b_2-1}n_a+1\leq j\leq
\sum\limits_{a=0}^{b_2}n_a)\\
(\mathfrak g_0\oplus\mathfrak g_{2\lambda_b})\setminus\{{\bf 0}\} & (\sum\limits_{a=0}^{b-1}n_1+1\leq i\leq\sum\limits_{a=0}^bn_a,\,\,\,\,
\sum\limits_{a=0}^{b-1}n_a+1\leq j\leq\sum\limits_{a=0}^bn_a),
\end{array}\right.\leqno{(11.15)}$$
($1\leq b_1<b_2\leq l$, $1\leq b\leq l$), where $\mathfrak g_{\lambda_{b_1}+\lambda_{b_2}}$ (resp. $\mathfrak g_{2\lambda_b}$) implies $\{{\bf 0}\}$ in the case of 
$\lambda_{b_1}+\lambda_{b_2}\notin\triangle_+$ (resp. $2\lambda_b\notin\triangle_+$).  
This gives the detail datas of the matrices $(c_{ij}^a)$'s ($1\leq a\leq N$).  

Take $\omega\in\mathcal A_P^{H^s},\,\,A_1,A_2\in T_{\omega}\mathcal A_P^{H^s}(\approx\Omega_{\mathcal T,1}^{H^s}(P,\mathfrak g)),\,\,u\in P$ and $\vv\in T_uP$.  
Let $A_i=\sum_{j=1}^NA_i^j\otimes\overline{\ve}_j$ ($i=1,2$).  
Then, from Lemma 11.3, we can derive 
$$\begin{array}{l}
\hspace{0.5truecm}\displaystyle{((d\square_{\bullet})_{\omega}(A_1,A_2))_u(\vv)}\\
\displaystyle{=-\sum_{l=1}^N\sum_{j=1}^N\sum_{k=1}^N\sum_{\sigma\in S_2}{\rm sgn}\,\sigma\,
\left(({\rm div}_{\nabla^{S,\omega}}A_{\sigma(1)}^j)_u\cdot(A_{\sigma(2)}^k)_u(\vv)\right.}\\
\hspace{4.75truecm}\displaystyle{\left.+2\langle\nabla^{S,\omega}_{\vv}A_{\sigma(1)}^j,(A_{\sigma(2)}^k)_u\rangle\right)\,c_{jk}^l\overline{\ve}_l}\\
\displaystyle{=-\sum_{l=1}^N\sum_{1\leq j<k\leq N}\sum_{\sigma\in S_2}{\rm sgn}\,\sigma\,
\left(({\rm div}_{\nabla^{S,\omega}}A_{\sigma(1)}^j)_u\cdot(A_{\sigma(2)}^k)_u(\vv)\right.}\\
\hspace{4.75truecm}\displaystyle{-({\rm div}_{\nabla^{S,\omega}}A_{\sigma(1)}^k)_u\cdot(A_{\sigma(2)}^j)_u(\vv)}\\
\hspace{4.75truecm}\displaystyle{+2\langle\nabla^{S,\omega}_{\vv}A_{\sigma(1)}^j,(A_{\sigma(2)}^k)_u\rangle}\\
\hspace{4.75truecm}\displaystyle{\left.-2\langle\nabla^{S,\omega}_{\vv}A_{\sigma(1)}^k,(A_{\sigma(2)}^j)_u\rangle\right)\,c_{jk}^l\overline{\ve}_l}
\end{array}\leqno{(11.16)}$$
where $\langle\,\,,\,\,\rangle$ denotes the inner product of $T_u^{\ast}P$ given by $(g_S^{\omega})_u$.  
Hence $(d\square_{\bullet})_{\omega}(A_1,A_2)={\bf 0}$ holds if and only if 
$$\begin{array}{l}
\displaystyle{\sum_{1\leq j<k\leq N}\sum_{\sigma\in S_2}{\rm sgn}\,\sigma\,
\left(({\rm div}_{\nabla^{S,\omega}}A_{\sigma(1)}^j)_u\cdot(A_{\sigma(2)}^k)_u(\vv)\right.}\\
\hspace{3.5truecm}\displaystyle{-({\rm div}_{\nabla^{S,\omega}}A_{\sigma(1)}^k)_u\cdot(A_{\sigma(2)}^j)_u(\vv)}\\
\hspace{3.5truecm}\displaystyle{+2\langle\nabla^{S,\omega}_{\vv}A_{\sigma(1)}^j,(A_{\sigma(2)}^k)_u\rangle}\\
\hspace{3.5truecm}\displaystyle{\left.-2\langle\nabla^{S,\omega}_{\vv}A_{\sigma(1)}^k,(A_{\sigma(2)}^j)_u\rangle\right)\,c_{jk}^l={\bf 0}}\\
\hspace{7truecm}\displaystyle{(1\leq l\leq N)}
\end{array}\leqno{(11.17)}$$
hold for any $u\in P$ and any $\vv\in T_uP$.  
From the detail datas of the matrices $(c_{ij}^a)$'s ($a=1,\cdots,N$) obtained from $(11.15)$, we can find many pairs $(A_1,A_2)$'s of elements of 
$\Omega_{\mathcal T,1}^{H^s}(P,\mathfrak g)(\approx T_{\omega}\mathcal A_P^{H^s})$ which does not satisfy $(11.17)$.  
Therefore $\square_{\bullet}$ is not closed.  \qed

\section{Proofs of Theorems D and E.} 
In this section, we prove Theorems D and E stated in Introduction.  
Denote by $\pi_{\mathcal G_{K_1,K_2}}$ and $\pi_{\mathcal G_{K_1,K_2;c}}$ the orbit maps of the actions 
$(\mathcal G_P^{H^{s+1}})_{K_1,K_2}\curvearrowright\mathcal A_P^{H^s}$ and \newline
$(\mathcal G_P^{H^{s+1}})_{K_1,K_2;c}\curvearrowright\mathcal A_P^{H^s}$, respectively.  
Also, denote by $\mathcal H^{\mathcal A}_{K_1,K_2}$ and ${\mathcal H}^{\mathcal A}_{K_1,K_2;c}$ the horizontal distributions of 
the orbit maps $\pi_{\mathcal G_{K_1,K_2}}$ and $\pi_{\mathcal G_{K_1,K_2;c}}$, respectively.  
We give explicit descriptions of $\mathcal H^{\mathcal A}_{K_1,K_2}$ and $\mathcal H^{\mathcal A}_{K_1,K_2;c}$.  

\vspace{0.25truecm}

\noindent
{\bf Lemma 12.1.} {\sl {\rm (i)}\ For an irreducible connection $\omega(\in\mathcal A_P^{H^s})$, we have 
$$\begin{array}{l}
\hspace{0.5truecm}\displaystyle{(\mathcal H^{\mathcal A}_{K_1,K_2})_{\omega}}\\
\displaystyle{=(\mathcal H^{\mathcal A})_{\omega}
\oplus{\rm Span}\{\,(\square_{\omega}^{-s}\circ d_{\omega})(G_{\Delta_{\omega}^0}^{\sigma(0)\cdot\vv_1}+G_{\Delta_{\omega}^0}^{\sigma(1)\cdot\vv_2})\,|\,
\vv_i\in\mathfrak k_i^{\perp}\,\,\,(i=1,2)\},}
\end{array}\leqno{(12.1)}$$
where $\mathcal H^{\mathcal A}$ denotes the horizontal ditribution of $\pi_{\mathcal M_P}:\mathcal A_P^{H^s}\to\mathcal A_P^{H^s}/\mathcal G_P^{H^{s+1}}$ (see Section 8).  

{\rm (ii)}\ For an irreducible connection $\omega(\in\mathcal A_P^{H^s})$, we have 
{\small 
$$\begin{array}{l}
\hspace{0.5truecm}\displaystyle{(\mathcal H^{\mathcal A}_{K_1,K_2;c})_{\omega}}\\
\displaystyle{=\left\{\left.\square_{\omega}^{-s}\left(\int_0^1\breve{\delta}_{c'(t)^{\ast}\otimes\sigma(t)\cdot\vv(t)}\,dt\right)\,\right|\,\vv\in H^s([0,1],\mathfrak g)
\,\,\,\,{\rm s.t.}\,\,\,\,\vv(0)\in\mathfrak k_1^{\perp},\,\vv(1)\in\mathfrak k_2^{\perp}\right\}.}
\end{array}\leqno{(12.2)}$$
}
}

\vspace{0.25truecm}

\noindent
{\it Proof.}\ \ Take any irreducible connection $\omega(\in\mathcal A_P^{H^s})$.  
According to (ii) of Lemma 6.5, we have 
$${\rm hol}_c((\mathcal G_P^{H^{s+1}})_{K_1,K_2}\cdot\omega)=(K_1\times K_2)\cdot{\rm hol}_c(\omega).$$
Furthermore, from Proposition 6.1 and the linearity of $\mu_c$, we can derive 
$$(\mathcal G_P^{H^{s+1}})_{K_1,K_2;c}\cdot\omega={\rm hol}_c^{-1}((K_1\times K_2)\cdot{\rm hol}_c(\omega)).\leqno{(12.3)}$$
Also we have 
$$\begin{array}{l}
\hspace{0.5truecm}\displaystyle{(\mathcal H^{\mathcal A}_{K_1,K_2})_{\omega}}\\
\displaystyle{=\{\eta\in\Omega_{\mathcal T,1}^{H^s}(P,\mathfrak g)\,|\,\langle d_{\omega}(\rho),\eta\rangle_s^{\omega}=0}\\
\hspace{3.25truecm}\displaystyle{(\forall\,\rho\in\Omega_{\mathcal T,0}^{H^s}(P,\mathfrak g)\,\,{\rm s.t.}\,\,(\rho(\sigma(0)),\rho(\sigma(1)))\in\mathfrak k_1\times\mathfrak k_2
\}}\\
\displaystyle{=\{\eta\in\Omega_{\mathcal T,1}^{H^s}(P,\mathfrak g)\,|\,\langle d_{\omega}(\rho),\square_{\omega}^s(\eta)\rangle_0=0}\\
\hspace{3.25truecm}\displaystyle{(\forall\,\rho\in\Omega_{\mathcal T,0}^{H^s}(P,\mathfrak g)\,\,{\rm s.t.}\,\,(\rho(\sigma(0)),\rho(\sigma(1)))\in\mathfrak k_1\times\mathfrak k_2
\}}\\
\displaystyle{=\{\eta\in\Omega_{\mathcal T,1}^{H^s}(P,\mathfrak g)\,|\,\langle\rho,d_{\omega}^{\ast}(\square_{\omega}^s(\eta))\rangle_0=0}\\
\hspace{3.25truecm}\displaystyle{(\forall\,\rho\in\Omega_{\mathcal T,0}^{H^s}(P,\mathfrak g)\,\,{\rm s.t.}\,\,(\rho(\sigma(0)),\rho(\sigma(1)))\in\mathfrak k_1\times\mathfrak k_2
\}}\\
\displaystyle{=\{\eta\in\Omega_{\mathcal T,1}^{H^s}(P,\mathfrak g)\,|\,(d_{\omega}^{\ast}\widehat{\circ\square_{\omega}^s)}(\eta)
=\delta_{\sigma(0)\cdot\vv_1}+\delta_{\sigma(1)\cdot\vv_2}}\\
\hspace{3.25truecm}\displaystyle{{\rm for}\,\,{\rm some}\,\,\vv_i\in\mathfrak k_i^{\perp}\,\,\,(i=1,2)\}}\\
\displaystyle{=\{\eta\in\Omega_{\mathcal T,1}^{H^s}(P,\mathfrak g)\,|\,(d_{\omega}^{\ast}\circ\square_{\omega}^s)(\eta)
=(d_{\omega}^{\ast}\circ d_{\omega})(G_{\Delta_{\omega}^0}^{\sigma(0)\cdot\vv_1}+G_{\Delta_{\omega}^0}^{\sigma(1)\cdot\vv_2})}\\
\hspace{3.25truecm}\displaystyle{{\rm for}\,\,{\rm some}\,\,\vv_i\in\mathfrak k_i^{\perp}\,\,\,(i=1,2)\}}\\
\displaystyle{=\{\eta\in\Omega_{\mathcal T,1}^{H^s}(P,\mathfrak g)\,|\,(d_{\omega}^{\ast}\circ\square_{\omega}^s)
(\eta-(\square_{\omega}^{-s}\circ d_{\omega})(G_{\Delta_{\omega}^0}^{\sigma(0)\cdot\vv_1}+G_{\Delta_{\omega}^0}^{\sigma(1)\cdot\vv_2}))=0}\\
\hspace{3.25truecm}\displaystyle{{\rm for}\,\,{\rm some}\,\,\vv_i\in\mathfrak k_i^{\perp}\,\,\,(i=1,2)\}}\\
\displaystyle{={\rm Ker}\,(d_{\omega}^{\ast}\circ\square_{\omega}^s)}\\
\hspace{0.6truecm}\displaystyle{\oplus{\rm Span}\{\,(\square_{\omega}^{-s}\circ d_{\omega})(G_{\Delta_{\omega}^0}^{\sigma(0)\cdot\vv_1}
+G_{\Delta_{\omega}^0}^{\sigma(1)\cdot\vv_2})\,|\,
\vv_i\in\mathfrak k_i^{\perp}\,\,\,(i=1,2)\}}\\
\displaystyle{=(\mathcal H^{\mathcal A})_{\omega}
\oplus{\rm Span}\{\,(\square_{\omega}^{-s}\circ d_{\omega})(G_{\Delta_{\omega}^0}^{\sigma(0)\cdot\vv_1}+G_{\Delta_{\omega}^0}^{\sigma(1)\cdot\vv_2})\,|\,
\vv_i\in\mathfrak k_i^{\perp}\,\,\,(i=1,2)\}.}
\end{array}\leqno{(12.4)}$$
Thus we obtain the relation $(12.1)$ in the statement (i).  

For $\vv\in T_xB$, define $\vv^{\ast}\in T_x^{\ast}B$ by $\vv^{\ast}(\cdot):=(g_B)_x(\vv,\cdot)$.  
For $x\in B$ and $\hat a\in(T^{\ast}B\otimes{\rm Ad}(P))_x$, denote by $\delta_{\hat a}$ the element of $\Gamma^{H^s}(T^{\ast}B\otimes{\rm Ad}(P))$ such that 
$\langle\delta_{\hat a},\widehat A\rangle_0=\langle\hat a,\widehat A_x\rangle_{B,\mathfrak g}$ holds for any $\widehat A\in\Gamma^{H^s}(T^{\ast}B\otimes{\rm Ad}(P))$.  
Also, denote by $\breve{\delta}_{\hat a}$ the element of $\Omega_{\mathcal T,1}^{H^s}(P,\mathfrak g)$ corresponding to $\delta_{\hat a}$.  
Since $\mu_c^{-1}(\hat{\bf 0})$ is given by 
$$\mu_c^{-1}(\hat{\bf 0})=\{A\in\Omega_{\mathcal T,1}^{H^s}(P,\mathfrak g)\,|\,A_{\sigma(t)}(\sigma'(t))=0\,\,\,\,(\forall\,t\in[0,1])\},$$
the orthogonal complement $\mu_c^{-1}(\hat{\bf 0})^{\perp}$ of $\mu_c^{-1}(\hat{\bf 0})$ wih respect to $\langle\,\,,\,\,\rangle_s^{\omega}$ is given by 
$$\begin{array}{l}
\hspace{0.5truecm}\displaystyle{\mu_c^{-1}(\hat{\bf 0})^{\perp}=\left\{\left.\square_{\omega}^{-s}
\left(\int_0^1\breve{\delta}_{c'(t)^{\ast}\otimes\sigma(t)\cdot\vv(t)}\,dt\right)\,\right|\,
\vv\in H^s([0,1],\mathfrak g)\right\}}\\
\displaystyle{=\left\{\left.\int_0^1G_{\square_{\omega}^s}(\cdot,c(t))(c'(t)^{\ast}\otimes\sigma(t)\cdot\vv(t))\,dt\,\right|\,\vv\in H^s([0,1],\mathfrak g)\right\}.}
\end{array}\leqno{(12.5)}$$
According to $(12.3)$, $(\mathcal G_P^{H^{s+1}})_{K_1,K_2;c}\cdot\omega$ is $(\mathcal G_P^{H^{s+1}})_{x_0}$-invariant, 
we have $(\mathcal H^{\mathcal A}_{K_1,K_2;c})_{\omega}\subset\widetilde{\mathcal H}^{\mathcal A}_{\omega}$, where $\widetilde{\mathcal H}^{\mathcal A}$ denotes the horizontal 
distribution of $\pi_{(\mathcal M_P)_{x_0}}:\mathcal A_P^{H^s}\to\mathcal A_P^{H^s}/(\mathcal G_P^{H^{s+1}})_{x_0}$ (see Section 8).  
This together with $(12.4)$ and $(12.5)$ derive 
{\small
$$\begin{array}{l}
\hspace{0.5truecm}\displaystyle{(\mathcal H^{\mathcal A}_{K_1,K_2;c})_{\omega}=(\mathcal H^{\mathcal A}_{K_1,K_2})_{\omega}\cap\mu_c^{-1}(\hat{\bf 0})^{\perp}}\\
\hspace{0truecm}\displaystyle{=\left\{\left.\square_{\omega}^{-s}\left(\int_0^1\breve{\delta}_{c'(t)^{\ast}\otimes\sigma(t)\cdot\vv(t)}\,dt\right)\,\right|\,
\vv\in H^s([0,1],\mathfrak g)\,\,\,\,{\rm s.t.}\right.}\\
\hspace{0.5truecm}\displaystyle{\left.
d_{\omega}^{\ast}\left(\int_0^1\breve{\delta}_{c'(t)^{\ast}\otimes\sigma(t)\cdot\vv(t)}\,dt\right)=\breve{\delta}_{\sigma(0)\cdot\vw_1}+\breve{\delta}_{\sigma(1)\cdot\vw_2}\,\,
{\rm for}\,\,{\rm some}\,\,\vw_i\in\mathfrak k_i^{\perp}\,\,\,(i=1,2)\right\}.}
\end{array}$$
}
Also, we have 
$$d_{\omega}^{\ast}\left(\int_0^1\breve{\delta}_{c'(t)^{\ast}\otimes\sigma(t)\cdot\vv(t)}\,dt\right)
=\int_0^1\breve{\delta}_{\frac{d}{dt}(\sigma(t)\cdot\vv(t))}\,dt=\breve{\delta}_{\sigma(1)\cdot\vv(1)}-\breve{\delta}_{\sigma(0)\cdot\vv(0)}.$$
Hence we obtain the relation $(12.2)$ in the statement (ii).  \qed

\vspace{0.25truecm}


By using these lemmas, we shall prove Theorem D stated in Introduction.  

\vspace{0.25truecm}

\noindent
{\it Proof of Theorem D.}\ \ 
Denote by $\mathcal U$ the set of all irreducible $H^s$-connections of the bundle $P$.  
The set $\mathcal U$ coincides with the set of all regular points of the action $\mathcal G_P^{H^{s+1}}\curvearrowright\mathcal A_P^{H^s}$, 
where ``regular point'' means that the orbit of this action through the point is principal.  
Since $\mathcal H^{\mathcal A}_{K_1,K_2;c}$ is of constant dimension on $\mathcal U$, 
its restriction $\mathcal H^{\mathcal A}_{K_1,K_2;c}|_{\mathcal U}$ to $\mathcal U$ gives a distribution on $\mathcal U$.  
By using $(12.2)$, we shall show that $\mathcal H^{\mathcal A}_{K_1,K_2;c}|_{\mathcal U}$ is involutive (and hence integrable).  
Take $\vX_1,\vX_2\in\Gamma^{H^{s+1}}(\mathcal H^{\mathcal A}_{K_1,K_2;c}|_{\mathcal U})$.  
According to $(12.2)$, $\vX_i$ ($i=1,2$) are described as 
$$(\vX_i)_{\omega}=\square_{\omega}^{-s}\left(
\int_0^1\breve{\delta}_{c'(t)^{\ast}\otimes\sigma(t)\cdot\vv^i_{\omega}(t)}\,dt\right)\quad\,\,(\omega\in\mathcal U),$$
where $\vv^i_{\omega}$ ($i=1,2$) are elements of $H^s([0,1],\mathfrak g)$ satisfying 
$\vv^i_{\omega}(0)\in\mathfrak k_1^{\perp}$ and $\vv^i_{\omega}(1)\in\mathfrak k_2^{\perp}$.  
Define $\vY_i\in\Gamma^{H^{s+1}}(\mathcal H^{\mathcal A}_{K_1,K_2;c}|_{\mathcal U})$ ($i=1,2$) by 
$$(\vY_i)_{\omega}:=\int_0^1\breve{\delta}_{c'(t)^{\ast}\otimes\sigma(t)\cdot\vv^i_{\omega}(t)}\,dt\quad\,\,(\omega\in\mathcal U).$$
Since 
$$[\breve{\delta}_{c'(t_1)^{\ast}\otimes\sigma(t_1)\cdot\vv_{\omega}^1(t_1)},\,\,\breve{\delta}_{c'(t_2)^{\ast}\otimes\sigma(t_2)\cdot\vv_{\omega}^2(t_2)}]
=\breve{\delta}_{c'(t_1)^{\ast}\otimes\sigma(t_1)\cdot[\vv_{\omega}^1(t_1).\vv_{\omega}^2(t_1)]}\times\delta_{t_1t_2}$$
($\delta_{t_1t_2}\,:\,$ the Kronecker's delta) holds, we have 
$$[\vY_1,\vY_2]_{\omega}=\int_0^1\breve{\delta}_{c'(t)^{\ast}\otimes\sigma(t)\cdot[\vv_{\omega}^1(t),\vv_{\omega}^2(t)]}\,dt.\leqno{(12.6)}$$
Let $\sigma^{\omega}:[0,1]\to P$ be the horizontal lift of $c$ satisfying $\sigma^{\omega}(0)=\sigma(0)$ with respect to $\omega$ and 
${\bf g}^{\omega}:[0,1]\to G$ be the curve in $G$ defined by $\sigma^{\omega}(t)=\sigma(t){\bf g}^{\omega}(t)$.  
Then we have 
$$\sigma(t)\cdot\vv_{\omega}^i(t)=\sigma^{\omega}(t)\cdot{\rm Ad}({\bf g}^{\omega}(t)^{-1})(\vv_{\omega}^i(t)).$$
From the definition of $Y_i$ and this relation, we can derive 
$$\square_{\omega}((\vY_i)_{\omega})=\int_0^1\breve{\delta}_{c'(t)^{\ast}\otimes\sigma^{\omega}(t)\cdot(({\rm id}+\frac{1}{a^2}\,\frac{d^2}{dt^2})
({\rm Ad}({\bf g}^{\omega}(t)^{-1})(\vv_{\omega}^i(t))))}\,dt\quad\,\,(\omega\in\mathcal U).\leqno{(12.7)}$$
From $(12.6)$ and $(12.7)$, we can derive 
$$\begin{array}{l}
\displaystyle{[\square_{\bullet}\vY_1,\square_{\bullet}\vY_2]_{\omega}
=\int_0^1\breve{\delta}_{c'(t)^{\ast}\otimes\sigma^{\omega}(t)\cdot
(({\rm id}+\frac{1}{a^2}\,\frac{d^2}{dt^2})({\rm Ad}({\bf g}^{\omega}(t)^{-1})([\vv_{\omega}^1(t),\vv_{\omega}^2(t)])))}\,dt}\\
\hspace{2.6truecm}\displaystyle{=\square_{\omega}\left(
\int_0^1\breve{\delta}_{c'(t)^{\ast}\otimes\sigma(t)\cdot[\vv_{\omega}^1(t),\vv_{\omega}^2(t)]}\,dt\right)}\\
\hspace{2.6truecm}\displaystyle{=\square_{\omega}([\vY_1,\vY_2]_{\omega})\qquad\qquad\quad\,\,(\omega\in\mathcal U).}
\end{array}$$
Similarly we can show 
$$[\square_{\bullet}^{-i}\vY_1,\square_{\bullet}^{-i}\vY_2]_{\omega}=\square_{\omega}([\square_{\bullet}^{-i-1}\vY_1,\square_{\bullet}^{-i-1}\vY_2]_{\omega})
\qquad(i=1,\cdots,s-1)\quad\,\,(\omega\in\mathcal U).$$
Hence we obtain 
$$[\vX_1,\vX_2]_{\omega}=\square_{\omega}^{-s}([\vY_1,\vY_2]_{\omega}).\leqno{(12.8)}$$
Also, we have $[\vv_{\bullet}^1(0),\vv_{\bullet}^2(0)]\in\mathfrak k_1^{\perp}$ and $[\vv_{\bullet}^1(1),\vv_{\bullet}^2(1)]\in\mathfrak k_2^{\perp}$ hold.  
These facts together with $(12.2),\,(12.6)$ and $(12.8)$ imply that $[\vX_1,\vX_2]_{\omega}\in(\mathcal H^{\mathcal A}_{K_1,K_2;c})_{\omega}$ ($\omega\in\mathcal U$).  
Thus $\mathcal H^{\mathcal A}_{K_1,K_2;c}|_{\mathcal U}$ is involutive, that is, integrable.  
Take $\omega\in\mathcal U$.  Let $\mathcal L_{\omega}$ be the maximal integral manifold of $\mathcal H^{\mathcal A}_{K_1,K_2;c}|_{\mathcal U}$ through $\omega$.  
Since $K_i$ ($i=1,2$) are symmetric subgroups of $G$, the action $K_1\times K_2\curvearrowright G$ is hyperpolar action.  
Denote by $\Sigma_{{\rm hol}_c(\omega)}$ the section of this action through ${\rm hol}_c(\omega)$.  
From $(12.3)$, we see that $\mathcal L_{\omega}$ extends a complete Riemannian Hilbert submanifold.  Denote by $\widetilde{\Sigma}_{\omega}$ this complete Riemannian Hilbert 
submanifold.  Also, according to Theorem A in \cite{K3}, ${\rm hol}_c$ is a homothetic submersion of $(\mathcal A_P^{H^s},g_S)$ onto $(G,g_G)$.  
Hence ${\rm hol}_c|_{\widetilde{\Sigma}_{\omega}}$ is a homothetic covering of $\widetilde{\Sigma}_{\omega}$ onto $\Sigma_{{\rm hol}_c(\omega)}$.  
Since $\Sigma_{{\rm hol}_c(\omega)}$ meets all $(K_1\times K_2)$-orbits orthogonally, it follows from $(12.3)$ that $\widetilde{\Sigma}_{\omega}$ meets all 
$(\mathcal G_P^{H^{s+1}})_{K_1,K_2;c}$-orbits orthogonally.  
Furthermore, since $K_1\times K_2\curvearrowright G$ is hyperpolar, the induced metric on $\Sigma_{{\rm hol}_c(\omega)}$ is flat and hence the induced metric on 
$\widetilde{\Sigma}_{\omega}$ also is flat.  Also, according to Theorem B, the action $(\mathcal G_P^{H^{s+1}})_{K_1,K_2;c}\curvearrowright(\mathcal A_P^{H^s},g_s)$ 
is horizontally isometric.  Therefore, this action is horizontally hyperpolar.  \qed

\vspace{0.5truecm}

{\small 
\centerline{
\unitlength 0.1in
\begin{picture}( 67.2500, 15.0500)( -9.9000,-24.3500)
\put(40.7000,-16.6000){\makebox(0,0)[rt]{$\mathcal A_P^{H^s}$}}%
\put(53.5000,-16.6000){\makebox(0,0)[lt]{$H^s([0,1],\mathfrak g)$}}%
\put(34.2000,-12.2000){\makebox(0,0)[rb]{$(\mathcal G_P^{H^{s+1}})_{K_1,K_2;c}$}}%
\put(46.2000,-12.3000){\makebox(0,0)[lb]{$P^{H^{s+1}}(G,K_1\times K_2)$}}%
%
\special{pn 8}%
\special{pa 3620 1160}%
\special{pa 4510 1160}%
\special{fp}%
\special{sh 1}%
\special{pa 4510 1160}%
\special{pa 4444 1140}%
\special{pa 4458 1160}%
\special{pa 4444 1180}%
\special{pa 4510 1160}%
\special{fp}%
\put(39.6000,-11.0000){\makebox(0,0)[lb]{$\lambda_c$}}%
%
\special{pn 8}%
\special{pa 4230 1740}%
\special{pa 5180 1740}%
\special{fp}%
\special{sh 1}%
\special{pa 5180 1740}%
\special{pa 5114 1720}%
\special{pa 5128 1740}%
\special{pa 5114 1760}%
\special{pa 5180 1740}%
\special{fp}%
\put(45.5000,-16.9000){\makebox(0,0)[lb]{$\mu_c$}}%
%
\special{pn 13}%
\special{ar 3280 1620 590 450  5.0767644 5.9136541}%
%
\special{pn 13}%
\special{pa 3840 1460}%
\special{pa 3880 1540}%
\special{fp}%
\special{sh 1}%
\special{pa 3880 1540}%
\special{pa 3868 1472}%
\special{pa 3856 1492}%
\special{pa 3832 1490}%
\special{pa 3880 1540}%
\special{fp}%
%
\special{pn 13}%
\special{pa 5190 1260}%
\special{pa 5220 1270}%
\special{pa 5252 1280}%
\special{pa 5280 1294}%
\special{pa 5310 1306}%
\special{pa 5338 1320}%
\special{pa 5366 1336}%
\special{pa 5394 1354}%
\special{pa 5420 1372}%
\special{pa 5444 1392}%
\special{pa 5468 1414}%
\special{pa 5490 1436}%
\special{pa 5508 1462}%
\special{pa 5526 1488}%
\special{pa 5528 1492}%
\special{sp}%
%
\special{pn 13}%
\special{pa 5538 1496}%
\special{pa 5570 1560}%
\special{fp}%
\special{sh 1}%
\special{pa 5570 1560}%
\special{pa 5558 1492}%
\special{pa 5546 1512}%
\special{pa 5522 1510}%
\special{pa 5570 1560}%
\special{fp}%
%
\special{pn 8}%
\special{pa 3890 1920}%
\special{pa 3890 2290}%
\special{fp}%
\special{sh 1}%
\special{pa 3890 2290}%
\special{pa 3910 2224}%
\special{pa 3890 2238}%
\special{pa 3870 2224}%
\special{pa 3890 2290}%
\special{fp}%
%
\special{pn 8}%
\special{pa 5730 1890}%
\special{pa 5730 2260}%
\special{fp}%
\special{sh 1}%
\special{pa 5730 2260}%
\special{pa 5750 2194}%
\special{pa 5730 2208}%
\special{pa 5710 2194}%
\special{pa 5730 2260}%
\special{fp}%
\put(41.4000,-23.4000){\makebox(0,0)[rt]{$\mathcal A_P^{H^s}/(\mathcal G_P^{H^{s+1}})_{K_1,K_2;c}$}}%
\put(56.2000,-23.3000){\makebox(0,0)[lt]{$\mathcal A_P^{H^s}/(\mathcal G_P^{H^{s+1}})_{K_1,K_2;c}=K_1\setminus G/K_2$}}%
%
\special{pn 8}%
\special{pa 4330 2430}%
\special{pa 5480 2430}%
\special{fp}%
\special{sh 1}%
\special{pa 5480 2430}%
\special{pa 5414 2410}%
\special{pa 5428 2430}%
\special{pa 5414 2450}%
\special{pa 5480 2430}%
\special{fp}%
\put(47.6000,-24.0000){\makebox(0,0)[lb]{$\overline{\mu}_c$}}%
%
\special{pn 8}%
\special{pa 4200 1870}%
\special{pa 5450 2280}%
\special{fp}%
\special{sh 1}%
\special{pa 5450 2280}%
\special{pa 5394 2240}%
\special{pa 5400 2264}%
\special{pa 5380 2278}%
\special{pa 5450 2280}%
\special{fp}%
%
\special{pn 8}%
\special{ar 4312 2142 72 72  0.6528466 5.7088805}%
%
\special{pn 8}%
\special{pa 4372 2102}%
\special{pa 4392 2152}%
\special{fp}%
\special{sh 1}%
\special{pa 4392 2152}%
\special{pa 4386 2084}%
\special{pa 4372 2102}%
\special{pa 4350 2098}%
\special{pa 4392 2152}%
\special{fp}%
\put(37.9000,-21.5000){\makebox(0,0)[rb]{$\pi_{\mathcal G_{K_1,K_2;c}}$}}%
%
\special{pn 8}%
\special{ar 5412 2002 72 72  0.6528466 5.7088805}%
%
\special{pn 8}%
\special{pa 5472 1962}%
\special{pa 5492 2012}%
\special{fp}%
\special{sh 1}%
\special{pa 5492 2012}%
\special{pa 5486 1944}%
\special{pa 5472 1962}%
\special{pa 5450 1958}%
\special{pa 5492 2012}%
\special{fp}%
\end{picture}%
\hspace{12.5truecm}}
}

\vspace{0.35truecm}

\centerline{{\bf Diagram 12.1$\,\,:\,\,$ The orbit space of the action 
$(\mathcal G_P^{H^{s+1}})_{K_1,K_2;c}\curvearrowright\mathcal A_P^{H^s}$}}

\vspace{0.5truecm}

Next we prove Theorem E stated in Introduction.  

\vspace{0.25truecm}

\noindent
{\it Proof of Theorem E}\ \ First we show the statement (i).  Let $(\widetilde M,\widetilde g)$ be a Riemannian Hilbert manifold, 
$\mathcal G\curvearrowright(\widetilde M,\widetilde g)$ be a horizontally hyperpolar action and 
$M:=\mathcal G\cdot u_0$ be a principal orbit of this action.  
Since $M$ has flat section, it is curvature-invariant.  
Hence the normal Jacobi operators of $M$ are defined.  
Denote by $\pi_{\mathcal G}$ the orbit map the action $\mathcal G\curvearrowright(\widetilde M,\widetilde g)$ and $\mathcal U$ the set of all regular points 
of this action.  Here we note that $\mathcal U$ is an open dense subset of $\mathcal A_P^{H^s}$.  
Since $\pi_{\mathcal G}|_{\mathcal U}$ is a Riemannian submersion with integrable horizontal distribution, we can show that the $\pi_{\mathcal G}$-projecctable 
unit normal vector field $\widetilde{\xi}$'s of $M$ are parallel with respect to the normal connection, where we note that the $\pi_{\mathcal G}$-projectability 
of $\widetilde{\xi}$ means that $\widetilde{\xi}$ is the horizontal lift of a unit vector $\bar{\xi}$ of $\mathcal U/\mathcal G$ (equipped with a flat metric) 
at $\pi_{\mathcal G}(u_0)$.  Hence a normal frame field consisting of parallel unit normal vector fields of the principal orbit exists.  
Hence the normal holonomy group of $M$ is trivial.  
Again, since $\pi_{\mathcal G}|_{\mathcal U}$ is a Riemannian submersion with integrable horizontal distribution of $(\mathcal U,g_s|_{\mathcal U})$ onto 
$\mathcal U/\mathcal G$ equipped with a flat metric, the singular Riemannnian foliation $\mathcal F$ consisting of the $\mathcal G$-orbits is a polar foliation.  
Since $M$ is a principal orbit of this action, it is a regular leaf of this polar foliation $\mathcal F$.  Hence it is equifocal (see Theorem 2.1 of \cite{AT}), 
where we note that, even if singular Riemannian foliations on a finite dimensional Riemannian manifold are treated in \cite{AT}, the discussion in \cite{AT} is valid for our case because 
the action $\mathcal G\curvearrowright(\widetilde M, \widetilde g)$ is a proper action of finite cohomogeneity.  

Next we show the statement (ii).  
Assume that $\mathcal G\curvearrowright(\widetilde M,\widetilde g)$ be a hyperpolar CSJ-action.  
Let $\widetilde{\xi}$ be a parallel normal vector field of the principal orbit $M$.  Take ${\bf g}\in\mathcal G$.  
Since $\widetilde{\xi}$ is $\pi_{\mathcal G}$-projectable, $({\bf g})_{\ast u}(\widetilde{\xi}_u)=\widetilde{\xi}_{{\bf g}\cdot u}$ holds.  
Hence, since  ${\bf g}$ is an isometry of $(\widetilde M,\widetilde g)$, we have 
$A_{\widetilde{\xi}_{{\bf g}\cdot u}}={\bf g}_{\ast u}\circ A_{\widetilde{\xi}_u}\circ{\bf g}_{\ast u}^{-1}$ and 
$\widetilde R(\widetilde{\xi}_{{\bf g}\cdot u})={\bf g}_{\ast u}\circ\widetilde R(\widetilde{\xi}_u)\circ{\bf g}_{\ast u}^{-1}$.  
Thus $M$ satisfies the conditions (WI-ii) and (WI-iii).  Therefore it is weakly isoparametric.  

Assume that $(\widetilde M,\widetilde g)$ is locally symmetric and that $\mathcal G\curvearrowright(\widetilde M,\widetilde g)$ is hyperpolar regularizable action.  
Let $\mathcal G\cdot u_0$ be a principal orbit of this action.  Assume that $\mathcal G\cdot u_0$ is curvature-adapted and of bounded curvature.  
Then it follows from the statement (i) that $\mathcal G\cdot u_0$ is weakly isoparametric.  Furthermore, from (ii) of Theorem A, we can derive that $\mathcal G\cdot u_0$ 
is isoparametric.  \qed

\vspace{0.5truecm}

{\small 
\rightline{Department of Mathematics, Faculty of Science}
\rightline{Tokyo University of Science, 1-3 Kagurazaka}
\rightline{Shinjuku-ku, Tokyo 162-8601 Japan}
\rightline{(koike@rs.tus.ac.jp)}
}

\end{document}